\documentclass[12pt,oneside]{amsart}

\usepackage{amsmath} 
\usepackage{amssymb}
\usepackage{amsfonts}

\numberwithin{equation}{section}

\newtheorem{thm}{Theorem}[section]
\newtheorem{lem}[thm]{Lemma}
\newtheorem{cor}[thm]{Corollary}
\newtheorem{prop}[thm]{Proposition}
\newtheorem{rem}[thm]{Remark}

\theoremstyle{definition}
\newtheorem{defn}[thm]{Definition}

\theoremstyle{remark}
\newtheorem{claim}[thm]{Claim}

\newcommand{\br}{\mathbf R}
\newcommand{\cal}{\mathcal}
\newcommand{\ddt}{\partial_t }
\newcommand{\dtt}{{\textstyle\frac{d}{dt}}}

\newcommand{\im}{\operatorname{Im}}
\newcommand{\lip}{\operatorname{Lip}}

\newcommand{\op}{\operatorname{Op}}
\newcommand{\ol}{\overline}

\newcommand{\pb}[1]{\left\{\,#1\,\right\}}
\newcommand{\re}{\operatorname{Re}}
\newcommand{\restr}[1]{\big|_{#1}}
\newcommand{\set}[1]{\left\{\,#1\,\right\}}
\newcommand{\sing}{\operatorname{\rm sing}}
\newcommand{\sgn}{\operatorname{sgn}}
\newcommand{\sh}{\sharp}
\newcommand{\supp}{\operatorname{\rm supp}}
\newcommand{\mn}[1]{\Vert#1\Vert}
\newcommand{\w}[1]{\langle #1\rangle }
\newcommand{\wf}{\operatorname{WF}}
\newcommand{\wt}{\widetilde}

\begin{document} 

\textwidth 155 mm
\textheight 253 mm

\baselineskip 16pt 
\lineskip 2pt
\lineskiplimit 2pt

\title[The Nirenberg-Treves  conjecture]{The resolution of the
  Nirenberg-Treves conjecture} 
\author[NILS DENCKER]{{\textsc Nils Dencker}}
\address{Department of Mathematics, University of Lund, Box 118,
S-221 00 Lund, Sweden}
\email{dencker@maths.lth.se}
\maketitle

\section{Introduction}\label{intro}

In this paper we shall study the question of local solvability of a
classical pseudo-differential operator $P \in {\Psi}^m_{cl}(M)$ on a
$C^\infty$ manifold $M$. Thus, we assume that the symbol of $P$ is an
asymptotic sum of homogeneous terms, and that $p = {\sigma}(P)$ is the
homogeneous principal symbol of~ $P$. We shall also assume that $P$ is of
principal 
type, which means that the Hamilton vector field $H_p$ and the
radial vector field are linearly independent when $p=0$.

Local solvability of $P$ at  a compact set $K \subseteq M$ means that
the equation 
\begin{equation}\label{locsolv}
Pu = v 
\end{equation}
has a local solution $u \in \cal D'(M)$ in a neighborhood of $K$
for any $v\in C^\infty(M)$
in a set of finite codimension.  Local $L^2$ 
solvability at $x_0$ of a first order pseudo-differential operator $P$
means that~ \eqref{locsolv} has a local solution $u\in L_{loc}^2(M)$ 
in a neighborhood of $K$ for any
$v\in L^2_{loc}(M)$ in a set of finite codimension.
We can also define microlocal solvability at any compactly based cone
$K \subset T^*M$, see ~\cite[Definition~26.4.3]{ho:yellow}.

It was conjectured by Nirenberg and Treves \cite{nt} that condition
(${\Psi}$) was equivalent to local solvability of pseudo-differential
operators of principal type.
Condition (${\Psi}$) means that
\begin{multline}\label{psicond} \text{$\im (ap)$
    does not change sign from $-$ to $+$}\\ 
 \text{along the oriented
    bicharacteristics of $\re (ap)$}
\end{multline}
for any $0 \ne a \in C^\infty(T^*M)$; actually it suffices to check
this for some $a\in C^\infty(T^*M)$ such that $H_{\re (ap)} \ne 0$ by
\cite[Theorem~26.4.12]{ho:yellow}. By the oriented bicharacteristics
we mean the positive flow-out of the Hamilton vector field $H_{\re
(ap)} \ne 0$ on $\re (ap) =0$. These are also called
semi-bicharacteristics.  Condition~\eqref{psicond} is invariant under
conjugation with elliptic Fourier integral operators and
multiplication with elliptic pseudo-differential operators,
see \cite[Lemma~26.4.10]{ho:yellow}.

The necessity of (${\Psi}$) for local solvability of
pseudo-differential operators was proved by Moyer in two dimensions
and by H\"ormander in general, see Corollary 26.4.8 in
\cite{ho:yellow}. In the analytic category, the sufficiency of
condition (${\Psi}$) for solvability of microdifferential operators
acting on microfunctions was proved by Tr\'epreau \cite{trep} (see
also \cite[Chapter~VII]{ho:conv}). The sufficiency of condition
(${\Psi}$) for local $L^2$ solvability of first order
pseudo-differential operators in two dimensions was proved by Lerner
\cite{ln:2d}, leaving the higher dimensional case open.

For differential operators, condition (${\Psi}$) is equivalent to
condition ($P$), which rules out any sign changes of $\im (ap)$ along
the bicharacteristics of $\re (ap)$, $0\ne a \in C^\infty(T^*M)$. The
sufficiency of ($P$) for local $L^2$ solvability of first order
pseudo-differential operators was proved by Nirenberg and Treves
\cite{nt} in the case when the principal symbol is real analytic, and
by Beals and Fefferman \cite{bf} in the general case.

Lerner \cite{ln:ex} constructed counterexamples to the sufficiency of
(${\Psi}$) for local $L^2$ solvability of first order
pseudo-differential operators, raising doubts on whether the condition
really was sufficient for solvability. But it was proved by the author
\cite{de:ln} that Lerner's counterexamples are locally solvable with
loss of at most two derivatives (compared with the elliptic
case). Observe that local $L^2$ solvability of first order
pseudo-differential operators means loss of one derivative.
There are several other results giving local solvability 
under conditions stronger than~ (${\Psi}$), see
\cite{de:suff}, \cite{ho:solv}, \cite{ln:coh} and \cite{ln:fact}.

In this paper we shall prove local solvability of principal type
pseudo-differential operators $P\in {\Psi}^m_{cl}(M)$ satisfying
condition~(${\Psi}$). This resolves the Nirenberg-Treves conjecture.
To get local solvability we shall assume a strong form of the
non-trapping condition at ~$x_0$: that all semi-characteristics are
transversal to the fiber $T^*_{x_0}\br^n$, i.e., $p(x_0,{\xi}) = 0
\implies \partial_{\xi} p(x_0,{\xi}) \ne 0$.

\begin{thm}\label{mainthm}
If $P \in {\Psi}^m_{cl}(M)$ is of principal type satisfying
condition~$({\Psi})$ near $x_0 \in M$ and $\partial_{\xi} p(x_0,{\xi})
\ne 0$ when $ p(x_0,{\xi}) = 0$, then $P$ is locally solvable near $x_0$.
\end{thm}

It follows from the proof that we lose at most two derivatives in the
estimate of the adjoint, which is one more compared with the condition
~($P$) case.  Thus our result has the consequence that hypoelliptic
operators of principal type can lose at most two derivatives. In fact,
if the operator is hypoelliptic of principal type, then the adjoint is
solvable of principal type, thus satisfying condition $({\Psi})$ and
we obtain an estimate of the operator.

Theorem~\ref{mainthm} is going to be proved by the construction of a
pseudo-sign which will be used in a multiplier
estimate. This resembles the constructions by Lerner in~ \cite{ln:2d}
and~\cite{ln:coh}, but here the pseudo-sign is not $L^2$ bounded.  The
symbol of the pseudo-sign is, modulo elliptic factors, essentially a
perturbation of the signed homogeneous distance to the
sign changes of the imaginary part of the principal symbol.

Observe that Theorem~\ref{mainthm} can be
microlocalized: if condition ~(${\Psi}$) holds microlocally near
$(x_0,{\xi}_0) \in S^*(M)$ then $P$ is microlocally solvable near
~$(x_0,{\xi}_0)$, see Corollary~\ref{mikrocor}. Since we
lose two derivatives in the estimate this is not trivial, it is a
consequence of the special type of estimate.

Most of the earlier results on local solvability have relied on finding
a factorization of the imaginary part of the principal symbol, see for
example ~\cite{de:suff} and ~\cite{ln:fact}. We have not
been able to find a factorization in terms of sufficiently good symbol
classes to get local solvability. The best result seems to be given by
Lerner~\cite{ln:per}, where a factorization was made which proved that
every first order principal type pseudo-differential operator
satisfying condition (${\Psi}$) is a sum of a solvable operator and an
$L^2$ bounded operator. But the bounded perturbation is in a bad symbol
class, and the solvable operator is not $L^2$~ solvable.

This paper is a shortened and simplified version of ~\cite{de:psi}, and
the plan is as follows.  In Section~\ref{reduc} we reduce
the proof of Theorem~\ref{mainthm} to the estimate of
Proposition~\ref{mainprop} for a microlocal normal form for ~$P^*$.
This estimate follows from a general multiplier estimate given by
Proposition~\ref{psest} in Section~\ref{mult}, and it involves
a pseudo-sign with properties given by
Proposition~\ref{apsprop}. The main part of the paper consists of the
construction of the pseudo-sign, showing that it has
the required properties, and it will not be completed until
Section~\ref{lower}. We must define suitable symbol classes and
weights, the new symbol classes will be defined in Section~\ref{symb}
and the new weights in Section~\ref{weight}.  The construction relies
on the local properties of symbols satisfying condition (${\Psi}$),
which will be derived in Section~\ref{local}. In order to obtain
the pseudo-sign, we shall use the Wick quantization and
suitable norms defined in Section~\ref{norm}, but the actual construction
of the  pseudo-sign will be carried out in
Section~\ref{aps}.

The author would like to thank Lars H\"ormander and Nicolas Lerner for
valuable comments leading to corrections and improvements 
of the proof.

\section{Reduction to the multiplier estimate}\label{reduc}

In this section we shall reduce the proof of Theorem~\ref{mainthm} to
an estimate for a microlocal normal form of the
adjoint for the operator.
We shall consider operators on the form
\begin{equation}\label{pdef1}
P_0 = D_t + i F(t,x,D_x) 
\end{equation}
where $F  \in C(\br, {\Psi}^1_{cl}(T^*\br ^n))$ has real principal symbol 
${\sigma}(F) = f$. Observe that we do not assume that $t \mapsto
f(t,x,{\xi})$ is differentiable. 
We shall assume
condition ~($\overline{\Psi}$):
\begin{equation} \label{pcond} \text{
$t \mapsto f(t,x,{\xi})$ does not change sign from $+$
to $-$ with increasing $t$  for any $(x,{\xi})$.}
\end{equation}
This means that the adjoint operator $P_0^*$ satisfies
condition (${\Psi}$).

We shall use the Weyl quantization of symbols $a \in \cal
S'(T^*\br^n)$: 
\begin{equation*}
a^w(x,D_x)u(x) = (2{\pi})^{-n}\iint
\exp{(i\w{x-y,{\xi}})}a\!\left(\tfrac{x+y}{2},{\xi}\right)u(y)\,dyd{\xi}
\qquad u\in C_0^\infty(\br^n).
\end{equation*} 
For Weyl calculus notations and results, see \cite[Section
18.5]{ho:yellow}. Observe that $\re a^w = (\re a)^w$ is the
symmetric part and $i\im a^w = (i\im a)^w$ the
antisymmetric part of the operator $a^w$. Also, if $a \in
S^m_{1,0}(T^*\br ^n)$ then $a(x,D_x) \cong a^w(x,D_x)$ modulo
${\Psi}^{m-1}_{1,0}(T^*\br ^n)$. In the following, we shall 
denote $ S^m_{{\varrho},{\delta}}(T^*\br ^n)$ by
$S^m_{{\varrho},{\delta}}$, $0 \le {\delta} \le {\varrho} \le 1$.

\begin{defn} 
We say that the symbol $b(x,{\xi})$ is in $S^{m}_{1/2,1/2}$
of first order, if $b$ satisfies the estimates in
$S^{m}_{1/2,1/2}$ for derivatives of order $\ge 1$.
\end{defn}

This means that the homogeneous gradient $(\partial_xb, |{{\xi}}|
\partial_{\xi}b) \in S^{m+\frac{1}{2}}_{1/2,1/2} $, and
implies that the commutators of $b^w$ with operators in
${\Psi}^k_{1,0}$ are in ${\Psi}^{m+k-1/2}_{1/2,1/2}$.
Observe that this condition is preserved when multiplying with symbols
in $S^0_{1,0}$.

We are going to prove an estimate for operators $P_0$ satisfying
condition~\eqref{pcond}. The estimate is not an $L^2$ estimate,
it gives a loss of two derivatives compared with the elliptic
case, but it is still localizable.  Let $\mn{u}_{(s)}$ be the
usual Sobolev norm, let $\mn u = \mn u_{(0)}$ be the $L^2$ norm, and
$\w{u,v}$ the corresponding inner product.

\begin{prop}\label{mikroprop} 
Assume that $P = D_t + i F^w(t,x,D_x)$, with $F \in C(\br,
S^1_{cl})$ having real principal symbol $f$ satisfying
condition~\eqref{pcond}. Then there exists $T_0 > 0$ such that if\/ $0
< T \le T_0$ then we can choose a real valued symbol
$b_T(t,x,{\xi}) \in L^\infty\left(\br, S^{1/2}_{1/2,1/2}\right)$
uniformly, with the property that $b_T \in
S^{0}_{1/2,1/2}$ of first order uniformly, and 
\begin{equation}\label{propest1} 
\mn u^2_{(-1/4)} \le T \im\w{P_0u,b_T^w u}
\end{equation} 
for $u(t,x) \in C_0^\infty(\br\times \br^{n})$ having support where
$|t| \le T \le T_0$. 
\end{prop}

Note that we have to change the multiplier $b_T$ when we change ~$T$,
but that the multipliers are uniformly bounded in the symbol
classes. By the calculus, the conditions on $b_T$ are preserved when
composing $b_T^w$ with symmetric operators in $L^\infty(\br ,
{\Psi}^0_{1,0})$.  Since $b_T^w \in {\Psi}^{1/2}_{1/2,1/2}$ the
estimate ~\eqref{propest1} could depend on lower order terms in the
expansion of~ $P_0$, and it is not obvious that it is localizable.

\begin{rem}\label{microrem} 
The estimate ~\eqref{propest1} can be perturbed with terms in
$L^\infty(\br, S^0_{1,0})$ in the symbol of ~$P_0$ 
for small enough~$T$. Thus it can be microlocalized: if
${\phi}(x,{\xi}) \in S^0_{1,0}$ is real valued 
then we have
\begin{equation}\label{cutofferror}
\im\w{P_0{\phi}^w u, b_{T}^w{\phi}^w u} \le
\im\w{P_0u,{\phi}^w b_T^w{\phi}^w u} + C\mn u^2_{(-1/4)}
\end{equation}
where ${\phi}^wb_T^w{\phi}^w$ satisfies the same conditions as~$b_T^w$.
\end{rem}

In fact, assume that $P_0 = D_t +if^w(t,x,D_x) + r^w(t,x,D_x)$ with $r\in
L^{\infty}(\br, S^0_{1,0})$. By conjugation with $E^w(t,x,D_x)$ where 
 \begin{equation*}
E(t,x,{\xi}) = \exp\left (-\int_{0} ^t \im
  r(s,x,{\xi})\,ds\right) \in L^{\infty}(\br, S^0_{1,0}),
 \end{equation*}
we can reduce to the case when $\im r \in L^{\infty}(\br,
S^{-1}_{1,0})$. We find that 
$b_T^w$ is replaced with $B_{T}^w = E^wb^w_{T}E^w$, which is real and satisfies
the same conditions as $b_{T}^w$ since $E$ is real and $B_T^w$ is
symmetric.
Clearly, the estimate ~\eqref{propest1} can be perturbed with terms in
$L^\infty(\br,  S^{-1}_{1,0})$ in the symbol expansion of~ $P_0$, and 
if $a(t,x,{\xi}) \in
L^\infty(\br,  S^0_{1,0})$ is real valued, then
\begin{equation}\label{2.4}
 \im\w{a^w u, b_{T}^w u} = \frac{1}{2i}\w{[b_T^w,a^w]u,u} \le C \mn
 u^2_{(-1/4)}
\end{equation}
since $b_T \in S^{0}_{1/2,1/2}$ of first order, $\forall \, t$.  We
also find that $[P_0,{\phi}^w] \cong \pb {f,{\phi}} ^w$ modulo
$L^{\infty}(\br, {\Psi}^{-1}_{1,0})$ where $\pb {f,{\phi}} \in
L^{\infty}\left(\br, S^0_{1,0}\right)$ is real valued.  By
using~\eqref{2.4} with $a = \pb {f,{\phi}}$, we obtain that the
estimate ~\eqref{propest1} is localizable.

\begin{proof} [Proof of Theorem~\ref{mainthm}] 
Take $w_0 = (x_0,{\xi}_0) \in p^{-1}(0)$, then since
$\partial_{\xi}p(w_0) \ne 0$ we may use Darboux'
theorem and the Malgrange Preparation Theorem, to obtain
coordinates $(t,y) \in M = \br\times \br^{n-1}$ so that $w_0 =
(0, (0,{\eta}_0))$ and
\begin{equation} \label{mikrop}
P^* \cong A(D_{t} +
iF(t, y,D_{y})) = AP_0
\end{equation}
microlocally in a conical neighborhood ${\Gamma}_{w_0}$ of $
w_0$, where $F \in C^\infty(\br, {\Psi}^{1}_{cl})$ has real
principal symbol $f$ satisfying condition~\eqref{pcond} and
$A \in {\Psi}^{m-1}_{cl}$ is elliptic (see the proofs of Theorems
~21.3.6 and~26.4.7' in ~\cite{ho:yellow}).

Take a cut-off function $0 \le {\psi}(y,{\tau},{\eta}) \in S^0_{1,0}$
such that ${\psi}$ is constant in ~$t$, ${\psi} =1$ in a conical
neighborhood of ~$w_0$, and
$\supp {\psi}\bigcap \set{|t| < T}\subset {\Gamma}_{w_0}$ when $T$ is
small enough. As in Remark ~\ref{microrem} we find by using
~\eqref{propest1} that
\begin{equation}\label{propest0} 
\mn{{\psi}^w u}^2_{(-1/4)}\le T \left(
\im\w{P^*u, B_T^wu}+ C\mn u^2_{(-1/4)}\right)
\end{equation}
for $u(t,y) \in C_0^\infty(\br^{n})$ having support where $|t| \le T$
is small enough (with $C$ independent of $T$). Here we have $B_T^w =
A_0^*{\psi}^wb_T^w{\psi}^w \in {\Psi}^{3/2-m}_{1/2,1/2}$ for a
microlocal inverse $A_0 \in {\Psi}^{1-m}_{1,0}$ to ~~$A$ in
${\Gamma}_{w_0}$.

By using a partition of unity, the estimate~\eqref{propest0} and
the Cauchy-Schwarz inequality, we obtain $R_T \in S^{0}_{1,0}$ for
~$T \le T_0$, such that $x_0 \notin \sing \supp R_T$ and
\begin{equation}\label{locest1}
\mn u_{(-1/4)} \le C \mn{P^*u}_{(7/4-m)} +
 \mn{R_T^wu} _{(-1/4)}
\end{equation}
for $u(x) \in C_0^\infty(\br^{n})$ having support where $|x| \le T \le
T_0$ is small enough. In fact, outside
$p^{-1}(0)$ we can construct a microlocal inverse in
${\Psi}^{-m}_{1,0}$ to $P^*$. 
Now conjugation with $\w{D_x}^s$ does not change the principal symbol
of $P$. Thus, for any $s \in \br$ we obtain positive constants $C_s$
and $T_s$ and $R_{T,s} \in S^{0}_{1,0}$ such that $0 \notin \sing
\supp R_{T,s}$ and
\begin{equation}\label{corest} 
 \mn u_{(s)} \le C_s \mn{P^*u} _{(s+2-m)}  + \mn{R_{T,s}^wu}_{(s)}
\end{equation} 
for $u(x) \in C_0^\infty(\br^{n})$ having support where
$|x| \le T_s$. This gives the local solvability of ~$P$ with a loss
of at most two derivatives, and completes the proof of
Theorem~\ref{mainthm}.
\end{proof}

It follows from the proof of Theorem~\ref{mainthm} that we also get
microlocal solvability for ~$P$ (see Definition~26.4.3 in
\cite{ho:yellow}). In fact, we can write $P$ on the
form~\eqref{mikrop} microlocally near $(x_0,{\xi}_0) \in p^{-1}(0)$ with
${\sigma}(F) = f$ satisfying~\eqref{pcond}. By using
\eqref{propest0} 
and the Cauchy-Schwarz inequality we obtain
~\eqref{corest} with $(x_0,{\xi}_0) \notin \wf R_{T,s}$,
which gives microlocal solvability.

\begin{cor}\label{mikrocor} 
If $P \in {\Psi}^m_{cl}(M)$ is of principal type near
$(x_0,{\xi}_0)\in T^*M$, satisfying condition~$({\Psi})$
microlocally near $(x_0,{\xi}_0)$, then $P$ is microlocally solvable
near $(x_0,{\xi}_0)$.
\end{cor}

In order to prove  Proposition~\ref{mikroprop} we shall need to
make a ``second microlocalization'' using the specialized symbol
classes of the Weyl calculus (see  \cite[Section
18.5]{ho:yellow}).
Let us recall the definitions: let  $g_{x,\xi}(dx,d{\xi})$ be a metric
on $T^*\br^n$, then we say that $g$ is slowly varying if there exists $c
> 0$ such that 
\begin{equation*}
g_{x,\xi}(x-y,{\xi}-{\eta}) < c \implies cg_{x,{\xi}} \le g_{y,{\eta}}
\le g_{x,{\xi}}/c \qquad\text{$(x,{\xi})$, $(y,{\eta}) \in T^*\br^n$.}
\end{equation*} 
Let ${\sigma}$ be the standard symplectic form on $T^*\br^n$, and let
\begin{equation*}
 g^{\sigma}_{x,\xi}(y,{\eta}) =
 \sup_{(t,{\tau})}|{\sigma}((y,{\eta});(t,{\tau}))|^2/ g_{x,\xi}(t,{\tau}) 
\end{equation*} 
be the dual form
of $(y,{\eta}) \mapsto  g_{x,\xi}({\sigma}(y,{\eta}))$.
We say that $g$ is ${\sigma}$~ temperate if it is slowly varying and
there exist constants $C$ and $N$ such that
\begin{equation*}
 g_{y,{\eta}} \le C g_{x,{\xi}}( 1 +
 g^{\sigma}_{y,{\eta}}(y-x,{\eta}-{\xi}))^N  \qquad
\text{$(x,{\xi})$, $(y,{\eta}) \in T^*\br^n$.} 
\end{equation*}
A positive real valued function $m(x,{\xi})$ on $T^*\br^n$ is $g$~
continuous if there exists a constant $c$ so that
\begin{equation*}
 g_{x,\xi}(x-y,{\xi}-{\eta}) < c \implies cm({x,{\xi}}) \le
 m({y,{\eta}}) \le m({x,{\xi}})/c  \qquad\text{$(x,{\xi})$,
   $(y,{\eta}) \in T^*\br^n$.} 
\end{equation*}
We say that $m$ is ${\sigma}$, $g$~ temperate if it is $g$ ~continuous and 
there exist constants $C$ and $N$ such that
\begin{equation*}
 m({y,{\eta}}) \le C m({x,{\xi}})( 1 +
 g^{\sigma}_{y,{\eta}}(y-x,{\eta}-{\xi}))^N  \qquad\text{$(x,{\xi})$,
   $(y,{\eta}) \in T^*\br^n$.}  
\end{equation*}
Let $S(m,g)$ be the class
of symbols $a \in C^\infty(T^*\br^n)$ with the seminorms 
\begin{equation*} 
|a|^g_j(x,{\xi}) = \sup_{T_i\ne 0}
\frac{|a^{(j)}(x,{\xi},T_1,\dots,T_j)|}{\prod_1^j
g_{x,{\xi}}(T_i)^{1/2}}\le C_j
m(x,{\xi})\qquad\forall\,(x,{\xi})\text{ for $j \ge 0$.}
\end{equation*}
We shall use metrics which are conformal, they shall be on
the form $g_{x,\xi}(dx,d{\xi}) = H(x,{\xi}) g^\sh(dx,d{\xi})$ where $0
< H(x,{\xi}) \le 1$ and $g^\sh$ is a 
constant symplectic metric: $(g^\sh)^{\sigma} = g^\sh$.  
In the following, we say that $m>0$ is a weight for a metric
$g$ if $m$ is ${\sigma}$, $g$ ~temperate.

\begin{defn}\label{s+def}
Let $m$ be a weight for the ${\sigma}$~
temperate metric $g$.  We say that $a \in S^+(m,g)$ if $|a|_j^g \le
C_jm$ for $j \ge 1$.
\end{defn}

For example, $b \in S^+(1, g_{1/2,1/2})$, with $g_{1/2,1/2} =
\w{{\xi}}|dx|^2 + |d{\xi}|^2/\w{{\xi}}$ at $(x,{\xi})$, if and only if
$b \in S^{0}_{1/2,1/2}$ of first order.  After localization, we shall
consider operators of the type
\begin{equation}\label{pdef2}
P_0 = D_t + i f^w(t,x,D_x) 
\end{equation}
where $f \in C(\br, S(h^{-1},hg^{\sh}))$ is real, and
$h \cong (1 + |{\xi}|)^{-1} \le 1$ is constant.
After a microlocal change of coordinates, we find that $S^k_{1,0}$
corresponds to ~$S(h^{-k}, h g^\sh)$ and $S^k_{1/2,1/2}$
corresponds to ~$S(h^{-k}, g^\sh)$ microlocally.

\begin{prop}\label{mainprop} 
Assume that $P_0 = D_t + i f^w(t,x,D_x)$, with real valued $f
(t,x,{\xi}) \in C(\br, S(h^{-1}, hg^\sh))$ satisfying
condition~\eqref{pcond}, here $0 < h \le 1$ and $g^\sh =
(g^\sh)^{\sigma}$ are constant. Then there exists $T_0 > 0$, such that
if\/ $0 < T \le T_0$ there exists a uniformly ${\sigma}$ ~temperate
metric $G_{T} = H_Tg^\sh$, $h \le H_T(t,x,{\xi}) \le 1$, a real valued
symbol $b_T(t,x,{\xi}) \in L^\infty(\br, S(H_T^{-1/2}, g^\sh) \bigcap
S^+(1, g^\sh))$ uniformly, such that
\begin{equation}\label{propest}
 h^{1/2}\mn u^2\le  T \im\w{P_0u,b_T^wu}
\end{equation} 
for $u(t,x) \in C_0^\infty(\br\times \br^{n})$ having support where $|t|
\le T \le T_0$. Here $T_0$ and the seminorms of $b_T$ only depend 
on the seminorms of $f$ in $S(h^{-1}, hg^\sh)$ when $|t| \le 1$.
\end{prop}

We find that $H_T$ is a weight for $g^\sh$ since $G_T \le g^\sh$.
The conditions on $b_T$ means in $g^\sh$ ~orthonormal coordinates that
$|b_T| \le CH_T^{-1/2}$ and
$|\partial_x^{\alpha}\partial_{\xi}^{\beta} b_T| \le
C_{{\alpha}{\beta}}$ when $|{\alpha}| + |{\beta}| \ge 1$.
Note that we have to change the multiplier $b_T$
when we change $T$, but the multipliers are uniformly bounded in the
symbol classes.
Observe that when $f \equiv 0$ 
we obtain ~\eqref{propest} with $b_T = 2th^{1/2}/T$.
As before, 
the estimate ~\eqref{propest} can be perturbed with terms in
$L^\infty(\br, S(1,hg^\sh))$ in the symbol of ~$P_0$ for
small ~$T$  (with changed~$b_T$), and thus it can be microlocalized.

\begin{proof} [Proof of Proposition~\ref{mikroprop}] 
By Remark~\ref{microrem} we may assume that $F = {\sigma}(F) = f$.
Take a partition of unity $\set{{\phi}_j(x,{\xi})}_j \in S^0_{1,0}$
and $\set{{\psi}_j(x,{\xi})}_j $, $\set{{\Phi}_j(x,{\xi})}_j \in
S^0_{1,0}$ such that $\sum_{j} {\phi}_j^2 = 1$, ${\phi}_j \ge 0$,
${\psi}_j \ge 0$, ${\Phi}_j \ge 0$, ${\phi}_j{\psi}_j = {\phi}_j$ and
${\psi}_j{\Phi}_j = {\psi}_j$, $\forall\ j$.  We may assume that the
supports are small enough so that $\w {\xi} \approx \w {{\xi}_j}$ in
$\supp {\Phi}_j$ for some~ ${\xi}_j$. Then, after choosing new
symplectic coordinates
\begin{equation*}
\left\{
\begin{aligned}
& y = x \w{{\xi}_j}^{1/2}\\
& {\eta} = {\xi} \w{{\xi}_j}^{-1/2}
\end{aligned}
\right. 
\end{equation*}
we obtain that $S^k_{1,0} = S(h_j^{-k}, h_j g^\sh)$ and $S^k_{1/2,1/2}
= S(h_j^{-k}, g^\sh)$ in $\supp {\Phi}_j$, where $h_j = \w
{{\xi}_j}^{-1}\le 1$, and $g^\sh = |dy|^2 + |d{\eta}|^2$. 
Observe that ${\phi}_j$,  ${\psi}_j$ and ${\Phi}_j \in S(1, h_j
g^\sh)$ uniformly in the new coordinates. 

We find that 
$  {\phi}_j^wP_0 \cong{\phi}_j^w P_{0j}$ modulo $S(h_j, h_jg^\sh)$, 
where $P_{0j} = D_t + if_j^w(t,y,D_y)$ with real valued 
$$f_j(t,y,{\eta}) =
{\psi}_j(y,{\eta})f(t,y,{\eta}) \in C(\br,S(h_j^{-1}, 
h_j g^\sh))
$$ satisfying condition~\eqref{pcond} uniformly in $j$.  By using
Proposition~\ref{mainprop} we obtain $h_j \le H_{j, T} \le 1$ and real
valued symbols $b_{j,T} \in L^\infty(\br, S(H_{j, T}^{-1/2}, g^\sh)
\bigcap S^+(1, g^\sh))$ uniformly, such that
\begin{equation}\label{mmest}
h_j^{1/2} \mn u^2 \le T \im\w{P_{0j}u,b_{j,T}^wu}
\end{equation} 
for $u(t,y) \in C_0^\infty(\br\times \br^{n})$ having support where
$|t| \le T \le T_0$. 
Since ${\Phi}_jf_j \equiv f_j$, we may 
replace $b_{j,T}$ by $b_{j,T}{\Phi}_j + 2th^{1/2}(1-{\Phi}_j)/T$ for small
enough $T$.
By substituting ${\phi}_j^wu$ in~\eqref{mmest} and perturbing with
terms in $\op S(h_j,h_jg^\sh)$, using that
${\psi}_j {\phi}_j = {\phi}_j$ and  $h_j = \w
{{\xi}_j}^{-1}$, we find as in
Remark~\ref{microrem} that
\begin{equation*}
\sum_j\im\w{P_{0j}{\phi}^w_ju,b_{j,T}^w{\phi}_j^wu} \le \im\w{P_0u,
b_Tu} + C\mn u^2_{(-1/4)}
\end{equation*}
where $b_T^w = \sum_j {\phi}_j^wb_{j,T}^w{\phi}^w_j \in L^\infty\left(\br,
{\Psi}^{1/2}_{1/2,1/2}\right)$.  Here we consider $\set{{\phi}_j}_j \in
S^0_{1,0}$ as having values in $\ell^2$, $\set{f_j}_j \in S^1_{1,0}$
and $\set{b_{j,T}}_j \in S^{1/2}_{1/2,1/2}$ as having values in $\cal
L(\ell^2, \ell^2)$ (observe that $\partial_{y,{\eta}} b_{j,T} \equiv
0$ outside $\supp {\Phi}_j$).  Then, since $b_T^w$ is symmetric we
have $ \im b_T \equiv 0$ and the calculus gives that $b_T \in
S^0_{1/2,1/2}$ of first order (see ~\cite[p.\ 169]{ho:yellow}).  Since
\begin{equation*}
\mn{u}^2_{(-1/4)} \le C \sum_j h_j^{1/2}\mn{{\phi}_j^wu}^2, 
\end{equation*}
we obtain ~\eqref{propest0}.   
This completes the proof of
Proposition~\ref{mikroprop}.
\end{proof}

The proof of Proposition~\ref{mainprop} relies on the multiplier
estimate in Proposition~\ref{psest}, it will be given at the end of
Section~\ref{mult}.

\section{The Multiplier Estimate}\label{mult}

Let $\cal B = \cal B(L^2(\br^n) )$ be the set of bounded operators
$L^2(\br^n)\mapsto L^2(\br^n)$.  We say that $A(t) \in C(\br, \cal B)$
if $A(t) \in \cal B$ for all $t \in \br$ and $t \mapsto A(t) u \in
C(\br,L^2 (\br^n))$ for any $u \in L^2 (\br^n)$.  We shall consider
the operator
\begin{equation}\label{pdef}
P = D_t + iF(t)
\end{equation}
where $F(t) \in C(\br, \cal B)$. In the applications, we
will have $F(t) \in C(\br, \op S(h^{-1},hg^\sh))$ where $h$ is constant.
We shall use multipliers which are not continuous in $t$.

\begin{defn} 
If ${\mathbf R} \ni t
\mapsto A(t) \in \cal B$, then we say that $A(t)$ is
weakly measurable if $t \mapsto A(t)u$ is weakly measurable for every
$u \in \cal L^2(\br^n)$, i.e., $t \mapsto \w{A(t)u,v}$ is
measurable for any ~$u$, $v \in L^2(\br^n)$.  
\end{defn}

If $A(t)$ is weakly measurable and locally bounded in ~$\cal B$, then we say that
$A(t) \in L_{loc}^\infty(\br, \cal B)$. In that case, we find that
$t \mapsto \w{A(t)u,u}\in L_{loc}^\infty(\br)$ has weak
derivative $\w{\dtt A(\cdot)u,u} \in \cal D'(\br)$ for any $u \in \cal
S(\br^n)$, given by $\w{\dtt A(\cdot)u,u}({\phi}) =
-\int\w{A(t)u,u}{\phi}'(t)\,dt$, ${\phi}(t) \in C_0^\infty(\br)$.  We
also have that if $u(t)$, $v(t) \in C(\br,L^2 
(\br^n))$ and $A(t) \in \cal B$ is weakly measurable,
then $t \mapsto \w{A(t)u(t),v(t)}$ is measurable.

In the following, we let $\mn u(t)$ be the $L^2$ norm of ~$u(t,x)$ in
$\br^n$  for fixed  ~
$t$, and~ $\w{u,v}(t)$ the corresponding inner product.
We shall use the following multiplier estimate (see also~~ \cite{ln:2d}
and ~~ \cite{ln:coh} for similar estimates).

\begin{prop} \label{psest}
Let $P= D_t + iF(t)$ with $F(t) \in C(\br, \cal B)$. Assume that
$B(t)= B^*(t)  \in L_{loc}^\infty(\br, \cal B)$, such that
\begin{equation}\label{multestcond}
\re \w{\dtt B(t) u,u}  + 2\re\w{B(t)u,F(t)u} \ge \re
 \w{m(t)u,u}\quad\text{in $\cal D'(I)\qquad \forall\ u \in \cal
   S(\br^n)$} 
\end{equation}
for an open interval $I \subseteq \br$, where $m(t) \in L_{loc}^\infty(\br, \cal B)$.
Then we have
\begin{equation}\label{pest}
\int  \re \w{m(t)u(t),u(t)} \,dt \le 2\int \im\w{Pu(t),B(t)u(t)} \,dt
\end{equation}
for $u \in C_0^1(I, \cal S(\br^n))$.
\end{prop}

\begin{proof}
Since $B(t)\in \cal B$ is weakly measurable and locally
bounded, we may for $u \in \cal S(\br^n)$ define the
regularization
\begin{equation*}
\w{B_{\varepsilon}(t)u,u} = {\varepsilon}^{-1} \int
\w{B(s)u,u}{\phi}((t-s)/{\varepsilon})  \,ds
= \w{Bu,u}({\phi}_{{\varepsilon},t})\qquad{\varepsilon}>0
\end{equation*} 
where ${\phi}_{{\varepsilon},r}(s) = {\varepsilon}^{-1}
{\phi}((r-s)/{\varepsilon})$ with $0 \le {\phi} \in
C_0^{\infty}(\br)$ satisfying $\int {\phi}\,dt = 1$. Then $t \mapsto
\w{B_{\varepsilon}(t)u,u}$ is in $C^\infty(\br)$ with derivative at
$t=r$ equal to $\w{\dtt B_{\varepsilon}(r)u,u} = \dtt
\w{Bu,u}({\phi}_{{\varepsilon},r})$.  
Let $I_0$ be an open interval such that $I_0 \Subset I$.  Then for
small enough ${\varepsilon}>0$ we find from
condition~\eqref{multestcond} that
\begin{equation}\label{newestcond}
 \re \w{\dtt B_{\varepsilon}(t) u,u} +
 2\re\w{Bu,Fu}({\phi}_{{\varepsilon},t}) \ge \re
 \w{mu,u}({\phi}_{{\varepsilon},t})\qquad t \in I_0 \quad u \in
 \cal S(\br^n).
\end{equation}
In fact, ${\phi}_{{\varepsilon},t} \ge 0$ and $\supp
{\phi}_{{\varepsilon},t} \in C_0^\infty(I)$ for small enough
${\varepsilon}$ when $t \in I_0$.

Now we define for $u\in C_0^1(I_0,\cal S(\br^n))$ and small enough
${\varepsilon}>0$
\begin{equation}\label{eq:def1}
M_{\varepsilon,u}(t) = \re \w{B_{\varepsilon}u,u}(t) =
{\varepsilon}^{-1} \int
\w{B(s)u(t),u(t)}{\phi}((t-s)/{\varepsilon}) \,ds.
\end{equation}
By differentiating under the
integral sign we obtain that $M_{\varepsilon,u}(t) \in C_0^1(I_0)$,
with derivative $\dtt M_{\varepsilon,u} = \re \w{(\dtt B_{\varepsilon})u,u} +
2\re\w{B_{\varepsilon} u,\ddt u}$
since $B(t) \in L_{loc}^\infty(\br, \cal B)$. By
integrating with respect to $t$, we obtain the vanishing average
\begin{equation}\label{aveq}
0 =  \int M_{\varepsilon,u}(t)\,dt  =  \int  \re\w{(\dtt
  B_{\varepsilon})u,u}\,dt  + \int 2\re\w{B_{\varepsilon} u,\ddt u}\,dt
\end{equation}
when $u\in C_0^1(I_0,\cal S(\br^n))$.  Since $\ddt u = iPu + Fu$ we
obtain from ~\eqref{newestcond} and~\eqref{aveq} that
\begin{multline*}
0\ge \iint \big(\re\w{m(s) u(t),u(t)} + 2\re\w{B(s)u(t),iP u(t)}\\
+ \re\w{B(s)u(t),\big(F(t) -F(s)\big)u(t)}\big)
{\phi}_{{\varepsilon} ,t}(s)\,dsdt. 
\end{multline*}
By letting ${\varepsilon} \to 0$ we obtain by dominated convergence that 
\begin{equation*}
0 \ge \int  \re \w{m(t)u(t),u(t)} +2\re\w{B(t) u(t),iPu(t)} \,dt
\end{equation*}
since $F(t) \in C(\br, \cal B)$, $u\in C_0^1(I_0,\cal S(\br^n))$,
$m(t)$ and $B(t)$ are bounded in ~$\cal B$ when $t \in \supp
u$.
Now $ 2\re\w{Bu,iPu} =
-2\im\w{Pu,Bu}$, thus we obtain~\eqref{pest}
for $u\in C_0^1(I_0,\cal S(\br^n))$. Since $I_0$ is an arbitrary open
subinterval with compact closure in ~$I$, this completes the proof of
the Proposition.
\end{proof}

Now we can reduce the proof of Proposition~\ref{mainprop} to the
construction of a  pseudo-sign $B = b^w$ in a fixed interval.

\begin{prop}\label{apsprop}
Assume that $f \in C(\br, S(h^{-1}, hg^\sh))$ is a real valued symbol
satisfying condition ~~$(\ol {\Psi})$ given by~\eqref{pcond}, here $0
< h \le 1$ and $g^\sh = (g^\sh)^{\sigma}$ are constant. Then there
exist a positive constant $c_0$, uniformly ${\sigma}$ ~temperate
metric $G_1 = H_1g^\sh$, $h \le H_1(t,x,{\xi}) \le 1$, real valued
symbols\/ $b(t,x,{\xi}) \in L^\infty(\br, S(H_1^{-1/2}, g^\sh) +
S^+(1, g^\sh))$ and\/ ${\mu}(t,x,{\xi}) \in L^\infty(\br, S(1,
g^\sh))$ such that for $u(x) \in C_0^\infty(\br^n)$ we have
\begin{equation*}\left\{
\begin{aligned}
&\w{\ddt( b^w)u,u} \ge \w{{\mu}^wu,u} \ge c_0h^{1/2}\mn{u}^2 \\
&\re \w{b^wf^wu,u} \ge - \w{{\mu}^wu,u}/c_0 
\end{aligned}
\right.\qquad\text{
  in $\cal D'(\br)$ when $|t| < 1$.}\label{apsest}
\end{equation*}
Here $c_0$, and the seminorms of\/
$b$ and $m$ only depend on the seminorms of\/ $f$ in $S(h^{-1},
hg^\sh)$ for $|t| \le 1$.
\end{prop}

\begin{proof}[Proof of Proposition~\ref{mainprop}] 
By doing a dilation $s = t/T$, we find that $P\/$ transforms into
$T^{-1}P_T = T^{-1}(D_s + i\,Tf_T^w(s,x,D_x))$, where $f_T(s,x,{\xi}) =
f(Ts,x,{\xi})$ satisfies the conditions in Proposition~\ref{apsprop}
uniformly in $T$ when $0 < T \le 1$. Thus we obtain real $b_T$, ${\mu}_T$ and
$c_0$ such that when $|s| < 1$ we have
\begin{equation*}\left\{
\begin{aligned}
&\w{\partial_s( b_T^w)u,u} \ge \w{{\mu}_T^wu,u} \ge c_0h^{1/2}\mn{u}^2\\
&\re \w{b_T^wf_T^wu,u} \ge -\w{{\mu}_T^wu,u}/c_0
\end{aligned}\right.
\qquad \text{in $\cal D'(\br)$} 
\end{equation*} 
for $u \in
C_0^\infty(\br^n)$. This implies that
\begin{multline*}
\w{\partial_{s}b^w_T(s,x,D_x)u,u} +
2\re\w{Tf_T^w(s,x,D_x)u,b^w_{T}(s,x,D_x)u} \\\ge (1 -
2T/c_0)\w{{\mu}_T^w(s,x,D_x)u,u} \qquad\text{ in $\cal D'\big(]{-1},1[\big)$}
\end{multline*}
for $u \in C_0^\infty(\br^n)$. Thus, for $T \le c_0/4$ we obtain by
using Proposition~\ref{psest} with $P_T = D_s
+ i\,Tf_T^w(s,x,D_x)$, $B(s) = b_T^w(s,x, D_x)$ and $m(s) =
{\mu}_T^w(s,x,D_x)$ that
\begin{equation*}
c_0h^{1/2}\int\mn u^2 \,ds \le \int \w{{\mu}_T^wu,u} \,ds \le 4 \int
\im\w{P_Tu,b^w_Tu}(s)\,ds
\end{equation*}
if $u \in C_0^\infty(\br\times \br^{n})$ has support where $|s|
< 1$. Finally, we obtain that 
\begin{equation*}
c_0h^{1/2}\int\mn u^2 \,dt \le 4T \int
\im\w{Pu,\wt b^w_Tu}(t)\,dt 
\end{equation*} 
with $\wt b_T(t,x,{\xi}) = b_T(t/T,x,{\xi})$ for $u \in C_0^\infty(\br\times
\br^{n})$ has support where $|t| < T \le c_0/4$.
\end{proof}

It remains to prove Proposition~\ref{apsprop}, which will be done in
Section ~\ref{lower}. The proof involves construction of a pseudo-sign
$b$ and a suitable weight ${\mu}$, and it will occupy the remaining part
of the paper.

\section{Symbol Classes and Weights}\label{symb}

In this section we shall define the symbol classes we shall use. In
the following, we shall denote $ (x,{\xi})$ by $w \in
T^*\br^n$, and we shall assume
that $f \in C(\br, S(h^{-1}, hg^\sh))$ satisfies condition ~$(\ol
{\Psi})$ given by ~\eqref{pcond}, here $0 <
h \le 1$ and $g^\sh = (g^\sh)^{\sigma}$ are constant.
We shall only consider the values of $f(t,w)$ when $|t| \le 1$, thus
for simplicity we let $f(t,w) = f(1,w)$ when $t \ge 1$ and  $f(t,w) =
f(-1,w)$ when $t \le -1$.

First, we shall define the signed distance function~ ${\delta}_0(t,w)$
in $T^*\br^n$ for fixed $t \in \br$, with the property that $t \mapsto
{\delta}_0(t,w)$ is non-decreasing and ${\delta}_0f \ge 0$.  Let
\begin{align}\label{xpmdef}
&X_+ = \set{(t,w) \in \br \times T^*\br^n : \exists\,s\le t,\ f(s,w) >0}\\
&X_- = \set{(t,w) \in \br \times T^*\br^n : \exists\,s\ge
 t,\ f(s,w) < 0}.\label{xpmdef1} 
\end{align}
We have that $X_\pm$ are open in $\br \times T^*\br^n$, and
by condition $(\ol{\Psi})$
we obtain that $X_- \bigcap X_+ = \emptyset$ and $\pm f \ge 0$ on $X_\pm$.
Let $X_0 = \br \times T^*\br^n\setminus \left(X_+\bigcup
  X_-\right)$, which is closed in
$\br \times T^*\br^n$, by the definition of $X_\pm$ we have
$f=0$ on $X_0$. 
Let
\begin{equation*}
d_0(t_0,w_0) = \inf\set{g^\sh(w_0-z)^{1/2}:\ (t_0,z) \in X_0}
\end{equation*}
which is the $g^\sh$ distance in $T^*\br^n$ to $ X_0$ for fixed $t_0$,
it could be equal to $+\infty$ in the case that 
$X_0\bigcap \set{t=t_0} = \emptyset$.  
Observe that $d_0$ is equal to the $g^\sh$ distance to $\complement
X_\pm$ on $X_\pm$.

\begin{defn}
We say that $w \mapsto a(w)$  is Lipschitz continuous on $T^*\br^n$
with respect to the metric ~$g^\sh$ if 
\begin{equation*}
\sup_{w \ne z \in  T^*\br^n}|a(w)- a(z)|/g^\sh(w-z)^{1/2} = C < \infty
\end{equation*} 
and $C$ is the Lipschitz constant of $a$. We shall denote by
$\lip(T^*\br^n)$ the Lipschitz continuous functions on $T^*\br^n$ with
respect to the metric $g^\sh$.
\end{defn}

Now if we have a bounded family $\set{a_j(w)}_{j \in J}$ of
Lipschitz continuous functions with uniformly bounded Lipschitz
constant $C$, then the infimum $A(w) = \inf_{j \in J} a_j(w)$ is
Lipschitz with the same constant (so also the supremum). 
By the triangle
inequality we find that $w \mapsto g^\sh(w-z)^{1/2}$ is  Lipschitz
continuous with respect to the metric $g^\sh$  with Lipschitz constant
equal to 1. By taking the infimum over ~$z$ we find that
$w \mapsto d_0(t,w)$ is Lipschitz
continuous with respect to the metric $g^\sh$ for those ~$t$ when it
is not equal to $\infty$, with Lipschitz constant equal to 1.
Since $X_0$ is closed we find that $t \mapsto d_0(t,w)$ is lower
semicontinuous, but we shall not use this fact.

\begin{defn} \label{signdef}
We define the sign of ~$f$ by
\begin{equation}\label{fsign}
\sgn (f) = 
\left\{
\begin{alignedat}{2}
 \pm 1& &\quad &\text{on  $X_\pm$}\\
 0& &\quad &\text{on  $X_0$}
\end{alignedat}
\right.
\end{equation}
then $\sgn(f)\cdot f \ge 0$.
\end{defn}

\begin{defn}\label{d0deforig}
We define the signed  distance function ${\delta}_0$ by
\begin{equation}\label{delta0def}
{\delta}_0(t,w) =  \sgn(f)(t,w)\min(d_0(t,w),h^{-1/2}).
\end{equation}
\end{defn}

By the definition we have that $|{\delta}_0| \le h^{-1/2}$ and
$|{\delta}_0| = d_0$ when  $|{\delta}_0| < h^{-1/2}$.
The signed distance function has the following properties.

\begin{rem}
The signed  distance function $w \mapsto {\delta}_0(t,w)$ given by
Definition~\ref{d0deforig} is Lipschitz 
continuous with respect to the metric $g^\sh$ with Lipschitz constant
equal to 1. We also find that 
${\delta}_0(t,w)f(t,w) \ge 0$ and $
t \mapsto {\delta}_0(t,w)$ is non-decreasing.
\end{rem}

In fact, since ${\delta}_0 =0$ on $X_0$ it suffices to show the
Lipschitz continuity of $w \mapsto {\delta}_0(t,w)$ on ~$X_+$ and
~$X_-$, which follows from the Lipschitz continuity of $w \mapsto
d_0(t,w)$. Since $(t,w) \in X_+$ implies $(s,w) \in X_+ $ for $s \ge 
t$ and $(t,w) \in X_-$ implies $(s,w) \in X_-$ for $s \le t$, we
find that $ t \mapsto {\delta}_0(t,w)$ is non-decreasing.

Since $t \mapsto {\delta}_0(t,w)$ is non-decreasing and bounded, it is
a regulated function. This means that the left and right limits
${\delta}_0(t\pm,w) = \lim_{0 < {\varepsilon} \to 0}{\delta}_0(t\pm
{\varepsilon},w)$ exist for any $(t,w)$ (see ~\cite{dieu}).  Since $t
\mapsto |{\delta}_0(t,w)|$ is lower semicontinuous, and $t
\mapsto{\delta}_0(t,w)$ is non-decreasing such that $\sgn({\delta}_0) =
\pm 1$ on $X_\pm$, we find that $t \mapsto {\delta}_0(t,w)$ is continuous
from left in ~$X_+$ and from right in ~$X_-$, but we shall not use
this fact.

In the following, we shall treat $t$ as a parameter, and denote $f' =
\partial_wf$ and $f'' = f^{(2)}$, where the differentiation is in the $w$
variables only. We shall also in the following assume that we have
choosen $g^\sh$ orthonormal coordinates so that $g^\sharp(dw) =
|dw|^2$. We shall use the norms $|f'|_{g^\sh} = |f'|$ and
$\mn{f''}_{g^\sh} = \mn {f''}$, but we shall omit the index ~$g^\sh$.
Next, we shall define the metric we are going to use.

\begin{defn}\label{g1def}
Define the weight
\begin{equation}\label{h2def}
H_1^{-1/2} = 1 + |{\delta}_0| + \frac{|f'|}{\mn {f''} + h^{1/4}|f'|^{1/2} +
  h^{1/2}}
\end{equation}
and the corresponding metric
$G_1 = H_1 g^\sh$.
\end{defn}

Since $|f'|/(\mn {f''} + h^{1/4}|f'|^{1/2} + h^{1/2})$ is continuous
in $(t,w)$ we find that $t \mapsto H_1^{1/2}(t,w)$ is a regulated
function.  We also have that
\begin{equation} \label{H1hest}
1 \le H_1^{-1/2}\le 1 + |{\delta}_0| + h^{-1/4} |f'|^{1/2}\le Ch^{-1/2}
\end{equation}
since $|f'| \le C_1 h^{-1/2}$ and $|{\delta}_0| \le
h^{-1/2}$. Moreover, $|f'| \le H_1^{-1/2}(\mn{f''} + h^{1/4}|f'|^{1/2}
+ h^{1/2})$ so the Cauchy-Schwarz inequality gives
\begin{equation}\label{dfest0}
|f'| \le 2 \mn{f''} H_1^{-1/2} + 3h^{1/2}H_1^{-1} \le CH_1^{-1/2}. 
\end{equation}
If $1 + |{\delta}_0| \ll H_1^{-1/2}$, then we find
from Definition~\ref{g1def} that $|f'| \gg h^{1/2}$ and thus
$H_1^{1/2} \cong \mn{f''}/|f'| + h^{1/4}/|f'|^{1/2}$.

By Proposition~\ref{g1prop} below the metric $G_1$ is ${\sigma}$ temperate.
The denominator 
\begin{equation} 
D = \mn {f''} + h^{1/4}|f'|^{1/2} + h^{1/2}
\end{equation}
in ~\eqref{h2def} may seem strange, but it has 
the following natural explanation which we owe to Nicolas
Lerner~\cite{ln:private}.

\begin{rem} \label{H1rem}
Let $F = h^{-1/2}f \in S(h^{-3/2}, hg^\sh)$, then
the largest $H_2 \le 1$ for which $F \in S(H_2^{-3/2}, H_2g^\sh)$ 
is given by
\begin{multline}\label{H2def}
 H_2^{-1/2} \cong 1 + |F|^{1/3} +
|F'|^{1/2} + \mn{F''} \\ = 1 + h^{-1/6}|f|^{1/3} + h^{-1/4}| f'|^{1/2} +
h^{-1/2} \mn{f''} \le Ch^{-1/2}
\end{multline}
modulo uniformly elliptic factors. When $ |{\delta}_0| \ll H_2^{-1/2}$
we obtain that $D \cong H_2^{-1/2}h^{1/2}$,
\begin{equation}\label{H2H1comp}
 H_1^{-1/2} \cong 1 +  |{\delta}_0| +  |F'|H_2^{1/2} \le CH_2^{-1/2}
\end{equation}
and by the Cauchy-Schwarz inequality we find that
$ 
H_2^{-1/2} \cong H_1^{-1/2} +
\mn{F''}.
$
Thus
\begin{equation}\label{fbisscomp}
H_1^{-1/2} \cong H_2^{-1/2}\quad \iff \quad \mn{f''} \le C h^{1/2}
H_1^{-1/2}\qquad \text{when $|{\delta}_0| \ll H_2^{-1/2}$}, 
\end{equation}
and then $H_2^{-1/2} \cong 1 + |F'|^{1/2}$ so $D \cong h^{1/4}|f'|^{1/2}
+ h^{1/2}$.
\end{rem}

Since $G_2$ is ${\sigma}$ temperate by Taylor's formula, we find that
~\eqref{H2H1comp} implies that $D$ is a weight for $G_1$ when $
|{\delta}_0| \ll H_1^{-1/2}$. But when $|{\delta}_0(w_0)| \cong
H_1^{-1/2}(w_0)$ we could have $D = h^{1/2}$ at $w_0$ with the
variation ${\Delta} D \ge c_0 h|w-w_0|^2 \gg h^{1/2}$ when $h^{-1/4}
\ll |w-w_0| \ll H_1^{-1/2}$. Thus $D$ is not a weight for $G_1$,
but one can show that $D + |{\delta}_0|H_1^{1/2}$ is.

The
advantage of using the metric $G_1$ is that when  $H_1 \ll 1$ in a
$G_1$ neighborhood of the sign changes we find that
$|f'| \ge c h^{1/2}$ is a
weight for $G_1$, ${\delta}_0 \in
S(H_1^{-1/2},G_1)$ and  the curvature 
of $f^{-1}(0)$ is bounded by $CH_1^{1/2}$ by Remark~\ref{dfsymbol} and
Proposition~\ref{ffactprop}. It follows from
Proposition~\ref{mproplem} that $f \in S(H_1^{-1}, G_1)$, but that
will not suffice for the proof of the conjecture.
We shall define the weight we are going to use, and in the following
we shall omit the parameter ~$t$.

\begin{defn} \label{Mdef}
Let
\begin{equation*}
M = |f| +|f'| H_1^{-1/2} + \mn {f''} H_1^{-1} + h^{1/2}H_1^{-3/2}
\end{equation*}
then we have $h^{1/2} \le M \le c h^{-1}$.
\end{defn}

\begin{prop}\label{g1prop}
We find that $G_1$ is ${\sigma}$ ~temperate  such that $G_1 =
H_1^2G_1^{\sigma}$ and
\begin{equation}\label{tempest}
H_1(w) \le C_0 H_1(w_0)(1 + H_1(w)g^\sh(w-w_0)) \le C_0 H_1(w_0)(1 +
g^\sh(w-w_0)). 
\end{equation}
We also have that $M$ is a weight for $G_1$ such that 
\begin{equation}\label{mtemp}
M(w) \le C_1 M(w_0)(1 + H_1(w_0)g^\sh(w-w_0))^{3/2}  \le C_1M(w_0) 
(1 + g^\sh(w-w_0))^{3/2}
\end{equation}
and we find that $f \in S(M,G_1)$. The constants $C_0$, $C_1$ and the
seminorms of $f$ in $S(M,G_1)$ only depend on the
seminorms of $f$ in $S(h^{-1},hg^\sh)$.
\end{prop}

\begin{proof}
Observe that since $G_1 \le g^\sh \le G_1^{\sigma}$ we find that the
conditions \eqref{tempest}--\eqref{mtemp} are stronger than the
property of being ${\sigma}$~temperate.  If $H_1(w_0) g^\sh(w-w_0) \ge
c> 0$ then we immediately obtain~\eqref{tempest} with $C_0 = c^{-1}$.
Thus, in order to prove~\eqref{tempest} it suffices to prove that $H_1(w) \le
C_0 H_1(w_0)$ when 
$H_1(w_0) g^\sh(w-w_0) \ll 1$, i.e., that $G_1 $ is slowly
varying.

First we consider the case $1 + |{\delta}_0(w_0)| \ge
H_1^{-1/2}(w_0)/2$. Then we 
find by the uniform Lipschitz continuity of $w \mapsto |{\delta}_0(w)|$ that
\begin{equation*}
H_1^{-1/2}(w) \ge 1 + |{\delta}_0(w)| \ge 1 + |{\delta}_0(w_0)| -
H_1^{-1/2}(w_0)/6 \ge  H_1^{-1/2}(w_0)/3
\end{equation*}
when $|w-w_0| \le H_1^{-1/2}(w_0)/6$, which gives
the slow variation in this case with $C_0 = 9$.

In the case $1 + |{\delta}_0(w_0)| \le H_1^{-1/2}(w_0)/2$ we have
\begin{equation} \label{h1dfest}
H_1^{-1/2}(w_0) \le 2|f'(w_0)|/(\mn {f''(w_0)} +
h^{1/4}|f'(w_0)|^{1/2} + h^{1/2}),
\end{equation}
thus we find
\begin{align}\label{4.17}
&\mn{f''(w_0)} \le 2H_1^{1/2}(w_0)|f'(w_0)| \qquad\text{ and}\\ 
&h^{1/2} \le 4H_1(w_0)|f'(w_0)|.\label{4.18}
\end{align}
By Taylor's formula we find
\begin{multline}\label{1.19}
|f'(w) -f'(w_0)| \le {\varepsilon}H_1^{-1/2}(w_0)\mn{f''(w_0)} + C_3
{\varepsilon}^2h^{1/2}H_1^{-1} (w_0) \\ \le (2{\varepsilon} + 4
C_3{\varepsilon}^2 )|f'(w_0)|\qquad\text{$|w-w_0| \le
{\varepsilon}H_1^{-1/2}(w_0)$}
\end{multline}
and when $2{\varepsilon} + 4
C_3{\varepsilon}^2 \le 1/2$ we obtain that
\begin{equation}\label{dfvar0}
1/2 \le |f'(w)|/|f'(w_0| \le 3/2 \qquad |w-w_0| \le
{\varepsilon}H_1^{-1/2}(w_0).
\end{equation}
Taylor's formula and ~\eqref{4.18} gives that
\begin{equation*}
\mn{f''(w)-f''(w_0)} \le C_3{\varepsilon}h^{1/2} H_1^{-1/2}(w_0)
\le 4 C_3{\varepsilon}H_1^{1/2}(w_0)|f'(w_0)|
\end{equation*}
when $|w-w_0| \le {\varepsilon}H_1^{-1/2}(w_0)$. We find that
$\mn{f''(w)}/|f'(w_0)|\le (2 +
4C{\varepsilon})H_1^{1/2}(w_0) 
$ when $|w-w_0| \le {\varepsilon}H_1^{-1/2}(w_0)$ by~\eqref{4.17} .
We then obtain from ~\eqref{h1dfest} and~\eqref{dfvar0} that
\begin{equation*}
 H_1^{1/2}(w) \le \mn{f''(w)}|f'(w)|^{-1} + h^{1/4}|f'(w)|^{-1/2} +
 h^{1/2}|f'(w)|^{-1} \le CH_1^{1/2}(w_0)
\end{equation*} 
when $|w-w_0| \le {\varepsilon}H_1^{-1/2}(w_0)$ and ${\varepsilon} \ll
1$, which gives the slow variation and ~\eqref{tempest}.

Next, we prove ~\eqref{mtemp}.
Taylor's formula gives as before that
\begin{equation}\label{taylorest0}
 \mn{f^{(k)}(w)} \le C\left(\sum_{j=0}^{2-k}
   \mn{f^{(k+j)}(w_0)}|w-w_0|^{j} + 
 h^{1/2}|w-w_0|^{3-k}\right) \qquad 0 \le k \le 2.
\end{equation}
By~\eqref{Mdef} we obtain that
\begin{equation*}
M(w) \le C\sum_{k=0}^{2}\mn{f^{(k)}(w_0)}(|w-w_0| + H_1^{-1/2}(w))^k +
Ch^{1/2} (|w-w_0| + H_1^{-1/2}(w))^3.
\end{equation*}
We obtain from ~\eqref{tempest} that $H_1^{-1/2}(w) \le C
H_1^{-1/2}(w) \le C (H_1^{-1/2}(w_0) + |w-w_0|)$. Thus we find
\begin{multline*}
 M(w) \le C\sum_{k=0}^{2}\mn{f^{(k)}(w_0)}H_1^{-k/2}(w_0) (1 +
 H_1^{1/2}(w_0)|w-w_0|)^k  \\ + 
Ch^{1/2}H_1^{-3/2}(w_0) (1 + H_1^{1/2}(w_0)|w-w_0|)^3 \le C'M(w_0) (1 +
H_1^{1/2}(w_0)|w-w_0|)^3 
\end{multline*}
which gives ~\eqref{mtemp}.

It is clear from the definition
of $M$ that  
$ 
\mn{f^{(k)}} \le MH_1^{k/2}$ when $k \le 2 
$, 
and when $k \ge 3$ we have
\begin{equation*}
 \mn{f^{(k)}} \le C_kh^{\frac{k-2}{2}}\le
 C'_kh^{1/2}H_1^{\frac{k-3}{2}}\le C'_k M H_1^{\frac{k}{2}}
\end{equation*}
since $h \le CH_1$ by~\eqref{H1hest} and $ h^{1/2}H_1^{-3/2} \le M$. 
This completes the proof.
\end{proof}

Observe that $f \in S(M,H_1g^\sh)$ for any choice of $H_1 \ge
ch$ in Definition~\ref{Mdef}, we do not have to use the other
properties of $H_1$. 
It follows from the
proof of Proposition~\ref{g1prop} that $|f'|$ is a weight for $G_1$
when $H_1^{1/2} \ll 1$ in a $G_1$ neighborhood of the sign changes.

\begin{rem}\label{dfsymbol}
When $1 + |{\delta}_0(w_0)| \le
H_1^{-1/2}(w_0)/2$ we have $|f'(w_0)| \ge  h^{1/2}/4$ and
\begin{equation}\label{dfvar}
1/C \le |f'(w)|/|f'(w_0| \le C \qquad \text{for\/ $|w-w_0| \le
{\varepsilon}H_1^{-1/2}(w_0)$.}
\end{equation}
We also have that $f' \in S(|f'|,G_1)$, i.e.,
\begin{equation} \label{dfsymbolest}
|f^{(k)}(w_0)| \le C_k |f'(w_0)|H_1^{\frac{k-1} {2}}(w_0)\qquad\text{for $k
\ge 1$}
\end{equation} 
when  $1 + |{\delta}_0(w_0)| \le
H_1^{-1/2}(w_0)/2$.
\end{rem}

In fact, \eqref{dfsymbolest} is trivial if $k =1$, follows from
~\eqref{4.17} for $k = 2$, and when $k \ge 3$ we have
$$|f^{(k)}(w_0)| \le C_kh^{\frac {k-2}2} \le 4C_k |f'|H_1 h^{\frac
{k-3}2}\le C'_k |f'| H_1^{\frac {k-1}2}
$$ 
by~\eqref{H1hest} and~\eqref{4.18}.

\begin{prop}\label{ffactprop}
Let $H_1^{-1/2}$ be given by Definition~\ref{g1def} for $f \in
S(h^{-1},hg^\sh)$. There exists ${\kappa_1}>0$ so that if\/
$|{\delta}_0(w_0)| \le {\kappa}_1H_1^{-1/2}(w_0)$, $H_1^{1/2}(w_0) \le
{\kappa_1}$ and  
\begin{equation}\label{4.19}
 \partial_{w_1}f(w_0) \ge c_0|f'(w_0)|
\end{equation}
for some $c_0>0$, then there exists $c_1 > 0$ such that
\begin{align}\label{ffactor}
&f(w) = {\alpha}_1(w)(w_1- {\beta}(w'))\\
&{\delta}_0(w) = {\alpha}_0(w)(w_1 - {\beta}(w'))\label{dfactor} 
\end{align}
when $|w- w_0| \le
 c_1H_1^{-1/2}(w_0)$.
Here\/ $0 < c_1 \le {\alpha}_0 \in S(1, G_1)$,
$c_1|f'| \le {\alpha}_1 \in S(|f'|, G_1)$ and
${\beta} \in S(H_1^{-1/2}, G_{1})$ only depends on $w'$, $w =
(w_1,w')$. 
The constant $c_1$ and the seminorms of ${\alpha}_j$, $j=1$, $2$, and
${\beta}$ only depend on $c_0$ and the seminorms of $f$ in $S(h^{-1},hg^\sh)$.
\end{prop}

\begin{proof}
We choose coordinates so that $w_0 = 0$, and put $H_1^{1/2} =
H_1^{1/2}(0)$.  Since $ H_1^{1/2} \le {\kappa_1}$ and $|{\delta}_0(0)| \le
{\kappa}_1H_1^{-1/2}$ we find from the Lipschitz continuity of
${\delta}_0$ and the slow variation that $1 + |{\delta}_0(w)| \le
H_1^{-1/2}(w)/2$ when $|w| \le {\varepsilon}H_1^{-1/2}$ for
sufficiently small ~${\varepsilon}$ and ~${\kappa_1}$. By
Remark~\ref{dfsymbol} we find 
that $|f'(w)| \le C|f'(0)|$ when $|w| \le
{\varepsilon}H_1^{-1/2}$ for small ${\varepsilon}$. 
We find from ~\eqref{dfsymbolest} and~\eqref{4.19} that
\begin{equation} \label{4.20}
\partial_{w_1}f(w) \ge
\partial_{w_1}f(0) - C'{\varepsilon}|f'(0)| \ge c_0|f'(0)|/2 \ge
c_0h^{1/2}/2C^2
\end{equation}
when $|w| \le {\varepsilon}H_1^{-1/2}$  for
${\varepsilon} \le c_0/2C'$.  Since $|{\delta}_0(0)| \le
{\kappa}_1H_1^{-1/2} \le C{\kappa}_1h^{-1/2}$ we find from
Definition~\ref{d0deforig} that $f(w)= 0$ for some $|w| \le
{\kappa}_1H_1^{-1/2}$ when $C{\kappa}_1 < 1$.  Thus, by the implicit
function theorem we can solve
\begin{equation*}
f(w) = 0 \iff w_1 = {\beta}(w')\qquad\text{when $|w|\le {\varepsilon}
  H_1^{-1/2}$} 
\end{equation*} 
for sufficiently small ${\varepsilon}$ and ${\kappa}_1$.
We find that $|{\beta}| =\cal O(H_1^{-1/2})$, $|{\beta}'| = 
|\partial_{w'}f|/|\partial_{w_1}f| = \cal O(1)$  and it follows from 
~\eqref{dfsymbolest} that
\begin{equation*}
|{\beta}''| \le C(|\partial_{w_1}^2f| |{\beta}'|^2 +
|\partial_{w'}\partial_{w_1}f||{\beta}'| +
\mn{\partial_{w'}^2f})/|\partial_{w_1}f| = \cal O(H_1^{1/2})
\end{equation*}
when $w_1 = {\beta}(w')$ and $|w|\le {\varepsilon} H_1^{-1/2}$.

Now we have by  Remark~\ref{dfsymbol} that $\partial_{w_j}f \in
S(|f'|,G_1)$ when $|w|\le {\varepsilon} H_1^{-1/2}$ for small enough
${\varepsilon}$, and by ~\eqref{4.20} we find $|f'|  \le
C\partial_{w_1}f$ when $|w|\le {\varepsilon} H_1^{-1/2}$. 
Thus we find $|\partial_w^{\alpha}f| \le C_{\alpha}
\partial_{w_1}fH_1^{({|{\alpha}| - 1})/{2}}$. 
Assume by induction that
$|{\beta}^{({\alpha})}| \le C_{\alpha} H_1^{({|{\alpha}| - 1})/{2}}$ when
$|w|\le {\varepsilon} H_1^{-1/2}$ for $|{\alpha}| < N$, where $N
\ge 3$. Then by differentiating the equation $f({\beta}(w'),w') = 0$
we find for $|{\alpha}| = N$
\begin{equation*} 
\partial^{\alpha}{\beta} = - \big(\sum c_{\gamma}
\partial^k_{w_1}\partial^{{\gamma}_0}_{w'}f
\prod_{j=1}^{k}\partial_{w'}^{{\gamma}_j}{\beta} +
\partial^{{\alpha}}_{w'}f\big)/\partial_{w_1}f\qquad\text{when $w_1 =
{\beta}(w')$ and $|w|\le {\varepsilon} H_1^{-1/2}$}
\end{equation*}
where the sum is over $k \ge 2$ and $
\sum_{j=1}^{k}{\gamma}_j + {\gamma}_0= {\alpha}$; or $k = 1$,
${\gamma}_0 \ne 0$ and ${\gamma}_0 + {\gamma}_1 ={\alpha}$.
In any case, we obtain that $|{\gamma}_j| < |{\alpha}|$ in the summation.
Since $f' \in S(|\partial_{w_1}f|, G_1)$, we obtain that
$|\partial^{{\alpha}}_{w'}f/\partial_{w_1}f| \le C_{\alpha}
H_1^{(|{\alpha}|-1)/2}$. By the induction hypothesis we find that  
\begin{equation*}
|\partial^k_{w_1}\partial^{{\gamma}_0}_{w'}f
\prod_{j=1}^{k}\partial_{w'}^{{\gamma}_j}{\beta}/\partial_{w_1}f| \le
C_{k,\gamma}H_1^{(k+|{\gamma}_0|-1)/2} H_1^{({\sum_{j = 1}^k|{\gamma}_j|-
    k})/{2}} = C_{k,\gamma}H_1^{(|{\alpha}|-1)/2} 
\end{equation*}
which completes the induction argument.

Now by using Taylor's formula we find $f(w) =
{\alpha}(w)(w_1-{\beta}(w'))$ where
\begin{equation*}
 {\alpha}(w) = \int_0^{1}
 \partial_{w_1}f({\theta}w_1 +(1-{\theta}){\beta}(w'),w')\,d{\theta} \qquad
 \text{for }|w| \le {\varepsilon}H_1^{-1/2}.
\end{equation*}
Thus ${\alpha}(w) \cong |f'(0)|$ when $|w| \ll
H_1^{-1/2}$ and ${\kappa}_1 \ll 1$, since then
$|{\beta}(w')| \ll H_1^{-1/2}$. Now
$\partial_{w_1}f \in S(|f'|, G_1)$  when
$|w| \le {\varepsilon}H_1^{-1/2}$, so we obtain that
${\alpha}(w) = f_0(w,{\beta}(w'))$ where $f_0 \in S(|f'|,G_1)$. Since
${\beta} \in 
S(H_1^{-1/2},G_1)$, differentiation gives
\begin{equation*} 
|\partial^{\gamma}{\alpha}| \le C_{\gamma}\sum_{\sum_{j=1}^{k} {\gamma}_j +
{\gamma}_0 = {\gamma}}  \left | \partial_{w_1}^{k}\partial_{w'}^{{\gamma}_0}f_0
\prod_{j=1}^k \partial^{{\gamma}_j}{\beta}\right | \le
C'_{\gamma}|f'| H_1^{({k +
|{\gamma}_0|}+ {\sum_{j=1}^k |{\gamma}_j|-k})/{2}} \le
C''_{\gamma}|f'| H_1^{|{\gamma}|/2}
\end{equation*}
which proves~\eqref{ffactor}. 

It remains to prove the statements about ${\delta}_0(w)$. It suffices
to prove that ${\delta}(w) = H_1^{1/2}(0){\delta}_0(w) \in
S(1,G_{1,0})$, here ${\delta}(w)$ is 
the signed $G_{1,0} = H_1(0)g^\sh$ distance to $X_0$. By choosing
\begin{equation*}
\left\{
\begin{aligned}
&z_1 = H_1^{1/2}(0)(w_1 - {\beta}(w'))\\
&z' = H_1^{1/2}(0) w' 
\end{aligned}\right. 
\end{equation*}
as new coordinates, then we find that $G_{1,0}$ transforms to a
uniformly bounded $C^\infty$ metric in a fixed neighborhood of the
origin. Now ${\delta}_1(z)$ is $\sgn(z_1)$ times the distance to $z_1
=0$ with respect to this metric, and this is a $C^\infty$ function in
a sufficiently small neighborhood of the origin. Clearly,
$|\partial_{z_1}{\delta}_1| \ge c > 0$ in a fixed neighborhood of the
origin, so Taylor's formula gives ${\delta}_1(z) = {\alpha}_0(z)z_1$,
where $c/2 \le {\alpha}_0 \in C^\infty$ in a smaller neighborhood.
This completes the proof of the Proposition.
\end{proof}

We shall compare our metric with the Beals--Fefferman metric $G 
= Hg^\sh$ on $ T^*\br^n$, where
\begin{equation}\label{hdef}
H^{-1} = 1 + |f| + |f'|^2 \le C h^{-1}.
\end{equation}
This metric is continuous in $t$,
 ${\sigma}$ ~temperate on $T^*\br^n$ 
and $\sup G/G^{{\sigma}} = H^{2}\le 1$.  We also have
$f \in S(H^{-1},G)$ (see for example the proof of
Lemma~26.10.2 in \cite{ho:yellow}).

\begin{prop}\label{mproplem}
We have that $H^{-1} \le CH_1^{-1}$ and $M \le C H_1^{-1}$,
which gives that $f \in S(M,G_1) \subseteq S(H_1^{-1},
G_1)$. When $|{\delta}_0| \le {\kappa}_0
H_1^{-1/2}$ for\/ $0 < {\kappa}_0$ sufficiently small we find that 
\begin{equation}\label{Mcomp}
 1/C \le M / (\mn{f''}H_1^{-1} + h^{1/2}H_1^{-3/2}) \le C.
\end{equation}
When  $|{\delta}_0| \le {\kappa}_0
H_1^{-1/2}$ and $H_1^{1/2} \le {\kappa}_0$ for\/ $0 < {\kappa}_0$
sufficiently small, we find 
\begin{equation} 
1/C_1 \le M/|f'|H_1^{-1/2} \le C_1.
\end{equation}
The constants only depend on the seminorms of
$f$ in $S(h^{-1},hg^\sh)$.
\end{prop}

Thus, the metric $G_1$ gives a coarser localization than the
Beals-Fefferman metric $G$ and smaller localization errors.

\begin{proof}
First note that by the Cauchy-Schwarz inequality we have
$$ M = |f| +|f'|H_1^{-1/2} + \mn {f''} H_1^{-1} +
h^{1/2}H_1^{-3/2} \le C(H^{-1} + H_1^{-1}).
$$ 
Thus we obtain $M \le CH_1^{-1}$  if we show that $H^{-1} \le
C H_1^{-1}$. Observe that we only have to prove this when
$|{\delta}_0| \ll H^{-1/2}$, since else
$H^{-1/2} \le C|{\delta}_0| \le C H_1^{-1/2}$.

If $|{\delta}_0(w_0)| \le {\kappa}H^{-1/2}(w_0) \le C{\kappa}h^{-1/2}$
for $C{\kappa}< 1$, then there exists $w \in f^{-1}(0)$ such that
$|w-w_0| \le {\kappa}H^{-1/2}(w_0) $. For
sufficiently small ${\kappa}$ we find from Taylor's formula and the
slow variation that $|f(w_0)| \le C{\kappa}H^{-1}(w_0)$. We obtain for
small enough ${\kappa}$ that
$$H^{-1}(w_0) \le
(1-C{\kappa})^{-1}(1 +|f'(w_0)|^2) \le C'H_1^{-1}(w_0)$$
which proves that $H^{-1} \le C H_1^{-1}$ and thus  $M \le CH_1^{-1}$.

Next, assume that $|{\delta}_0(w_0)| \le {\kappa}_0 H_1^{-1/2}(w_0)$.
When ${\kappa}_0 $ is small enough, we find as before that there
exists $w \in f^{-1}(0)$ such that $|w-w_0| = |{\delta}_0(w_0)|$. By
Remark~\ref{dfsymbol} we obtain that $|f'|$ only varies with a fixed
factor in $|w-w_0| \le {\kappa}_0 H_1^{-1/2}(w_0)$ for ${\kappa}_0 \ll
1$. Since $f(w)=0$, Taylor's formula gives that $|f(w_0)| \le
C{\kappa}_0|f'(w_0)|H_1^{-1/2}(w_0)$. Then we find that $M \le
C(|f'|H_1^{-1/2} + \mn {f''} H_1^{-1} + h^{1/2}H_1^{-3/2}) \le C'(\mn
{f''} H_1^{-1} + h^{1/2}H_1^{-3/2})$ by ~\eqref{dfest0}, which gives
~\eqref{Mcomp}.
When also $H_1^{1/2}(w_0) \le {\kappa}_0 \ll 1$  we find from
\eqref{4.17}--\eqref{4.18} that 
$
M(w_0) \le C(|f'|H_1^{-1/2} + \mn {f''} H_1^{-1} + h^{1/2}H_1^{-3/2})
\le  C'|f'(w_0)|H_1^{-1/2}(w_0).
$
This completes
the proof of the Proposition. 
\end{proof}

When ${\delta}_0  \in S(H_1^{-1/2},G_1)$ we find from the Fefferman-Phong
inequality that $({\delta}_0f)^w \ge r^w$ for some $r \in
S(H_1^{1/2},G_1)$ since ${\delta}_0f \ge 0$. But this will not suffice
for the proof, instead we shall use that $r \in S(MH_1^{3/2},G_1)$.  
In the next section, we shall estimate
the term $MH_1^{3/2} \cong \mn{f''}H_1^{1/2} + h^{1/2}$ near
the sign changes (see Proposition~\ref{symbest}).

\section{Local Properties of the Symbol}\label{local}

In this section we shall study the local properties of the symbol near
the sign changes. We start with a one dimensional result.

\begin{lem} \label{lem1}
Assume that $f(t) \in C^3(\br)$ such that
$\mn{f^{(3)}} \le h^{1/2}$ is bounded. If
\begin{equation}\label{fsgncond0}
\sgn(t)f(t) \ge 0\qquad\text{when ${\varrho}_0 \le |t|
\le {\varrho}_1$}
\end{equation} 
for ${\varrho}_1 \ge 3{\varrho}_0 >0$, then we find
\begin{align}
&|f(0)| \le \frac{3}{2}\left(f'(0){\varrho}_0 + 
  h^{1/2}{\varrho}_0^3/2\right)\label{f0est}\\
&|f''(0)| \le f'(0)/{\varrho}_0 + 7
  h^{1/2}{\varrho}_0/6.\label{d2fest}
\end{align}
\end{lem}

We obtain ~\eqref{f0est} since we are close to the zeroes of $f$, and
~\eqref{d2fest} since we have a lower bound on $f'$ in the domain.

\begin{proof}
By Taylor's formula we have
\begin{equation*}
0 \le \sgn({t})f({t}) = |{t}| f'(0) +
\sgn({t})(f(0) +f''(0){t}^2/2) + R({t}) \qquad
{\varrho}_0 \le |{t}| \le {\varrho}_1
\end{equation*} 
where $|R({t})| \le  h^{1/2}|{t}|^3/6$. We obtain that
\begin{equation}\label{5.2a}
\left|f(0)+ t^2 f''(0)/2\right| \le
f'(0)|{t}| + h^{1/2}|{t}|^3/6
\end{equation} 
for any $|{t}| \in [{\varrho}_0, {\varrho}_1]$. 
By choosing
$|{t}| = {\varrho}_0$ and $|{t}| = 3{\varrho}_0$, we obtain that
$$4{\varrho}_0^2 |f''(0)| \le 
4 f'(0){\varrho}_0 + 28 h^{1/2}{\varrho}_0^3/6
$$ 
which gives~\eqref{d2fest}. By letting $|t| ={\varrho}_0$  in
~\eqref{5.2a} and substituting ~\eqref{d2fest}, we obtain~\eqref{f0est}. 
\end{proof}

Next, we study functions $f(w) \in S(h^{-1}, hg^\sh)$ where ~~$h $ is
constant, then $\mn{f^{(3)}} \le C_3h^{1/2}$.  In the following we
shall assume that $H_0^{1/2}$ is a parameter, but later on it will be
given by the metric $G_1$, see~\eqref{h0small}.

\begin{prop}\label{newestprop}
Let $f \in S(h^{-1}, hg^\sh)$ for constant\/ $h$ and let\/ $H_1^{1/2}$
be given by 
by Definition~\ref{g1def}. Assume that there exist positive
constants $C_0$  and $H_0^{1/2}$ such that  $H_0^{1/2} \ge h^{1/2}$ and 
\begin{equation}\label{newlemass}
\sgn(w_1)f(w) \ge 0\quad\text{when $(1 + H_0^{1/2}|w'|^2)/C_0 \le |w_1| \le
  C_0H_0^{-1/2}$ and $|w'| \le C_0 H_0^{-1/2}$}
\end{equation}
where $w= (w_1,w')$. If\/
$H_0^{1/2} \le {\kappa}_0$ is sufficiently small, then
there exist $c_1$ and $C_1$ such that
\begin{align}\label{fdfest}
&|f(0)|  \le \partial_{w_1}f(0){\varrho}_0 +  C_1h^{1/2}{\varrho}_0^3 \\
\label{ddfest}
&\mn{f''(0)} \le \partial_{w_1}f(0)/{\varrho}_0 +
C_1h^{1/2}{\varrho}_0  
\end{align}
for any $1 \le {\varrho}_0 \le c_1H_0^{-1/2}$.
Here $c_1$, $C_1$ and ${\kappa}_0$ only depend on $C_0$ and the
seminorms of\/ ~$f$ in $S(h^{-1},hg^\sh)$.
\end{prop}

\begin{proof}
We shall consider the function $t \mapsto f(t,w')$ which satisfies
~\eqref{fsgncond0} for fixed $w'$ with 
$$(1 + H_0^{1/2}|w'|^2)/C_0 = {\varrho}_0(w') \le
 {\varrho}_1(w')/3 = C_0H_0^{-1/2}/3 $$ 
if $H_0^{1/2} (1 + H_0^{1/2}|w'|^2) \le C_0^2/3$. We obtain this
when $|w'| \le {\kappa}_0H_0^{-1/2}$ and $H_0^{1/2} \le {\kappa}_0$
if ${\kappa}_0 + {\kappa}_0^2 \le C_0^2/3$, which we assume in what 
follows. Then \eqref{f0est} and~\eqref{d2fest} gives that
\begin{align}
\label{f0est0}
&|f(0,w')| \le   \frac{3}{2}\partial_{w_1}f(0,w'){\varrho}_0 +
 Ch^{1/2}{\varrho}_0^3\\
\label{d2fest0}
&|\partial_{w_1}^2f(0,w')| \le \partial_{w_1}f(0,w')/{\varrho}_0 +
 Ch^{1/2}{\varrho}_0\qquad \text{}
\end{align}
for ${\varrho}_0(w') \le
{\varrho}_0 \le  {\varrho}_1(w')/3$ and $|w'| \le {\kappa}_0H_0^{-1/2}$.
By letting $w' = 0$ we obtain 
~\eqref{fdfest} from~\eqref{f0est0}, and by taking $w' =0$
in~\eqref{d2fest0} we find that
\begin{equation}\label{d2w1fest}
|\partial_{w_1}^2f(0)| \le \partial_{w_1}f(0)/{\varrho}_0  +
Ch^{1/2}{\varrho}_0
\end{equation} 
for $1 \le {\varrho}_0 \le c_1H_0^{-1/2}$ and $c_1 \le C_0/3$. 
By letting ${\varrho}_0 = {\varrho}_0(w')$ in~\eqref{f0est0}, dividing
with $3{\varrho}_0(w')/2$ and using Taylor's formula
for $w' \mapsto \partial_{w_1}f(0,w')$, we
find that there exists $C > 0$ so that
\begin{equation*}
0 \le \partial_{w_1}f(0) + \w{w',\partial_{w'}(\partial_{w_1}f)(0)} +
Ch^{1/2}(1 +|w'|^2)  \qquad\text{} 
\end{equation*}
when $|w'| \le
{\kappa}_0H_0^{-1/2}$ since then ${\varrho}_0(w') \le (1+
{\kappa}_0|w'|)C_0^{-1}$. Thus 
by optimizing over fixed ~$|w'|$, we  
obtain that  
\begin{equation}\label{dd1fest}
|w'||\partial_{w'}(\partial_{w_1}f)(0)| \le  \partial_{w_1}f(0) +
Ch^{1/2}(1 +|w'|^2)  \qquad\text{when $|w'|
  \le {\kappa}_0 H_0^{-1/2}.$} 
\end{equation}
By again letting ${\varrho}_0 = {\varrho}_0(w')$ in~\eqref{f0est0},
using Taylor's formula for $w' \mapsto \partial_{w_1}f(0,w')$ 
and substituting ~\eqref{dd1fest}, we obtain 
\begin{equation*}
 |f(0,w')| \le  \partial_{w_1}f(0){\varrho}_0(w')  +
C h^{1/2}(1 + |w'|^3))\qquad\text{when $|w'| \le
   {\kappa}_0 H_0^{-1/2}$.} 
\end{equation*}
By considering the odd and even terms in Taylor's formula for $w'
\mapsto f(0,w')$, optimizing over fixed ~$|w'|$ using~\eqref{fdfest}
with ${\varrho}_0 = 1$, we obtain that
\begin{equation}\label{dfddfest}
 \left|\partial_{w'}f(0)\right||w'| + \mn{\partial_{w'}^2f(0)}|w'|^2/2 
\le \partial_{w_1}f(0){\varrho}_0(w')  +
C h^{1/2}(1 + |w'|^3))
\end{equation}
when $|w'| \le {\kappa}_0 H_0^{-1/2}$.  Now we have
$C_0^{-1}\le {\varrho}_0(w') \le C_0^{-1}(1 + {\kappa}^2_0)$ if $|w'| \le
{\kappa}_0H_0^{-1/4}$.  Thus we obtain ~\eqref{ddfest} by taking $|w'|
= {\varrho}_0 \in [1, {\kappa}_0H_0^{-1/2}]$ in ~\eqref{dd1fest},
$|w'|^2 = {\varrho}_0 \in [1, {\kappa}_0^2H_0^{-1/2}]$ in
~\eqref{dfddfest} and using~\eqref{d2w1fest}. This finishes the proof
of the Proposition.
\end{proof}

We find from ~\eqref{dfddfest} that we also get an estimate on
$|\partial_{w'}f(0)|$. But we shall need an estimate on the direction
of the gradient of $f$.

\begin{rem}\label{estdfrem}
Assume that $f \in S(h^{-1},hg^\sh)$ for constant\/ $h$,  $f(0) = 0$ and
\begin{equation}\label{semifest}
f(w) \ge 0\quad\text{when $c_1 \le w_1 \le |w| \le C_1$}
\end{equation}
for some $c_1 < C_1$. There exist
${\kappa}_0 > 0$ and $C_0 > 0 $ such that if $H_1^{1/2}(0) \le
{\kappa}_0$ then
\begin{equation}\label{semidfest}
|\partial_{w'}f(0)| \le C_0|\partial_{w_1}f(0)|.
\end{equation}
The values of~ ${\kappa}_0$ and ~$C_0$ only
depend on~$c_1$, $C_1$ and the seminorms of  $f \in S(h^{-1},hg^\sh)$.
\end{rem}

In fact, by Taylor's  formula for $w_1 \mapsto f(w_1,w')$ we find
from~\eqref{semifest} that 
\begin{equation*}
 -f(0,w') \le  C(|\partial_{w_1}f(0,w')| + \mn{f''(0,w')} +
 h^{1/2})\qquad \text{when $|w'| \le \sqrt{C_1^2 -c_1^2} $}.
\end{equation*}
By using the Taylor expansion of $w' \mapsto
f(0,w')$ and of  $w' \mapsto \partial_{w_1} f(0,w')$ we find that
\begin{equation*} 
 \w{w',d_{w'}f(0)} \le C'(|\partial_{w_1}f(0)| + \mn{f''(0)} +
 h^{1/2})\qquad \text{when $|w'| \le \sqrt{C_1^2 -c_1^2} $}.
\end{equation*}
By optimizing over fixed $|w'|> 0$ we find 
\begin{equation*}
 |\partial_{w'}f(0)| \le C''(|\partial_{w_1}f(0)| + \mn{f''(0)}+ h^{1/2}).
\end{equation*}
When $H_1^{1/2}(0) \ll 1$ we obtain
from~\eqref{4.17} and~\eqref{4.18} that $ \mn{f''(0)}+ h^{1/2} \ll
|f'(0)|$ since $f(0) = 0$, which gives ~\eqref{semidfest}.
Next, we shall estimate the weight $M$ near the sign
changes, this will be important for the lower bounds in Section~\ref{lower}.

\begin{prop}\label{symbest}
Let $f \in S(h^{-1}, h g^\sh)$ for constant\/ $h$ and let\/
$H_1^{1/2}$ be given by 
by Definition~\ref{g1def}. Assume that there exist positive
constants  $C_0$ and $H_0^{1/2}$ such that
\begin{equation}\label{newlemass0}
\sgn(w_1)f(w) \ge 0\quad\text{when $(1 + H_0^{1/2}|w'|^2)/C_0 \le |w_1| \le
  C_0H_0^{-1/2}$ and $|w'| \le C_0 H_0^{-1/2}$}
\end{equation}
where $w= (w_1,w')$, and $h^{1/2}\le H_0^{1/2} \le 1$. 
Then we obtain that
\begin{equation}
M(0)H_1^{3/2}(0) \le C_1H_0^{1/2}.\label{Mest}
\end{equation}
Here $C_1$ only depends on $C_0$ and the seminorms of $f$ in
$S(h^{-1},hg^\sh)$.
\end{prop}

\begin{proof}
From the definition of $M$ we find
$ 
MH_1^{3/2} = |f| H_1^{3/2} + |f'| H_1 + \mn{f''} H_1^{1/2} + h^{1/2}. 
$ 
In the following we shall denote $H_1^{1/2} = H_1^{1/2}(0)$.  
First we observe that if $H_1^{1/2} \le C_0H_0^{1/2}$ then
$M(0)H_1^{3/2} \le CH_1^{1/2} \le CC_0 H_0^{1/2}$ by
Proposition~\ref{mproplem}. Thus, in
the following we shall assume $H_0^{1/2} \le {\kappa}_0H_1^{1/2}
\le {\kappa}_0$
for some ${\kappa}_0> 0$ to be determined later.
For small enough ~${\kappa}_0$ we find from  ~ \eqref{fdfest} in 
Proposition~\ref{newestprop} that
\begin{equation}\label{fdfest0}
|f(0)|  \le C(|f'(0)| + h^{1/2}).
\end{equation}
Since $H_1 \le 1$ it suffices to estimate $\mn{f^{(k)}}$ for $k=1$, $2$.
We obtain from ~\eqref{dfest0} that
\begin{equation}\label{dfdddfest}
|f'(0)|  H_1 \le 2\mn{f''(0)}H_1^{1/2} + 3h^{1/2}.
\end{equation}
Thus, it remains to estimate
$\mn{f''(0)}H_1^{1/2}$
in order to obtain~\eqref{Mest}.
Now $H_1^{-1/2} \le {\kappa}_0H_0^{-1/2}$,
and by ~ \eqref{ddfest} we have 
\begin{equation*}
 \mn{f''(0)}H_1^{1/2} \le C_1
H_1^{1/2}(|f'(0)|/{\varrho}_0 + h^{1/2}{\varrho}_0)
\end{equation*}
for $1 \le {\varrho}_0 \le c_1H_0^{-1/2}$.
Thus, if ${\kappa}_0 \le c_1/4C_1$, we can choose ${\varrho}_0
= 4C_1H_1^{-1/2} \le 4C_1{\kappa}_0 H_0^{-1/2} \le c_1H_0^{-1/2}$
which gives
\begin{equation}\label{5.20} 
\mn{f''(0)}H_1^{1/2} 
\le
\frac{1}{4}|f'(0)|H_1 + Ch^{1/2} \le \frac{1}{2}
\mn{f''(0)}H_1^{1/2} + C_2h^{1/2}.
\end{equation}
by \eqref{dfdddfest}. This gives
$
 \mn{f''(0)}H_1^{1/2} \le 2 C_2h^{1/2} \le 2 C_2C_0H_0^{1/2}
$
and completes the proof.
\end{proof}

It follows from the proof that we obtain the bound $MH^{3/2} \le
Ch^{1/2}$ in the case when $H_0^{-1/2} \gg H_1^{-1/2}$. 
Next, we shall investigate the conditions we need in order to obtain
~\eqref{newlemass0}. 

\begin{lem}\label{betalemma}
Let $f \in S(h^{-1},hg^\sh)$ satisfy condition $(\ol {\Psi})$ given by
~\eqref{pcond}, let\/ ${\delta}_0$ be given by
Definition~\ref{d0deforig} and\/ $H_1^{1/2}$ by
Definition~\ref{g1def}.  Assume that there exist\/ $t'\le t_0 \le t''$
such that
\begin{align}
&{\Delta} = \max_{t = t',\, t''} |{\delta}_0(t,w_0) -{\delta}_0(t_0,w_0)| \le 1
\label{d0var}\\ 
&|{\delta}_0(t,w_0)|  \le{\varrho}_0 H_0^{-1/2}\qquad\text{for $t =
  t'$ and $t''$,} 
\label{h0h1est1}
\end{align}
where 
\begin{equation} \label{h0small}
H_0^{1/2} = \max(H_1^{1/2}(t',w_0),H_1^{1/2}(t'',w_0))
\le {\varrho}_0.
\end{equation} 
If ${\varrho}_0$ is sufficiently small in ~\eqref{h0h1est1}--\eqref{h0small}
then there exist $g^\sh$ orthonormal coordinates $w=
(w_1,w')$ so that $w_0 = (z_1,0)$  with $z_1 =
{\delta}_0(t',w_0)$ and 
\begin{equation} \label{betapmdef}
\left\{
\begin{aligned} 
&\sgn f(t',w) = \sgn(w_1-{\beta}_-(w'))\\
&\sgn f(t'',w) = \sgn(w_1-{\beta}_+(w'))
\end{aligned}\right.\qquad \text{for $|w|\le
c_0 H_0^{-1/2}$}
\end{equation}
where ${\beta}_\pm(w') \in S(H_0^{-1/2},H_0g^\sh).$ We also find
that ${\beta}_-(0) = |{\beta}'_-(0)|
=0$, $|{\beta}_+(0)| \le C {\Delta}$ and $|{\beta}'_+(0)| \le C
H_0^{1/4}{\Delta}^{1/2}$.  Here ${\varrho}_0$, $c_0$ and $C$ only depend
on the seminorms of $f$ in $S(h^{-1},hg^\sh)$.
\end{lem}

\begin{proof}
Since $|{\delta}_0(t',w_0)| \le {\varrho}_0 H_0^{-1/2}\le
C{\varrho}_0h^{-1/2}$ by~\eqref{h0h1est1} and ~\eqref{H1hest}, we find
by ~\eqref{delta0def} there exists $(t',\ol w) \in X_0$, so
that $|w_0 -\ol w| = |{\delta}_0(t',w_0)|$ if $C{\varrho}_0 <1$,
which we shall assume in what follows.  By the slow variation, we find
from ~\eqref{h0h1est1} and~\eqref{h0small} that $H_1^{1/2}(t,\ol
w)  \le C H_0^{1/2}  \le C {\varrho}_0 \ll 1 $ for $t = t'$ and
$t''$ when ${\varrho}_0 \ll 1$. We may 
choose $g^\sh$ orthonormal coordinates so that $\ol w = 0$ and $w_0 = (z_1,0)$
with $z_1 = {\delta}_0(t',w_0)$.
 
Since $H_1^{1/2} (t',0) \le C{\varrho}_0 \ll 1$ and
${\delta}_0(t',0) =0$, we find from Remark~\ref{dfsymbol} that
$|f'(t',0 )| \ne 0$. Since $f(t',0 ) = 0$ and
${\delta}_0(t',(z_1,0))= z_1$, we find $\partial_{w'}f(t',0 ) =
0$. By using Proposition~\ref{ffactprop} at $(t',0)$ we obtain
~\eqref{betapmdef} for $t = t'$ when 
${\varrho}_0$ are small enough. Here
${\beta}_-(w') \in S(H_1^{-1/2}(t',0),H_{1}(t',0)g^\sh)$ when $|w'|
\le c_1 H_1^{-1/2}(t',0)$, and ${\beta}_-(0) = |{\beta}'_-(0)| =0$.  We
find that $\mn{{\beta}_-''} \le C_2H_0^{1/2}$, thus by using Taylor's
formula we obtain that ${\beta}_-(w') \in S(H_0^{-1/2}, H_0 g^\sh)$ when
$|w'| \le c_0 H_0^{-1/2}$ since $H_0^{-1/2} \le C H_1^{-1/2}(t',0)$. It
is clear that we may also obtain 
~\eqref{betapmdef} for $t= t''$ but not necessarily with the same
coordinates.

\begin{claim}\label{lemclaim} 
When ${\varrho}_0 > 0$ is small enough in ~\eqref{h0h1est1}, there
exists $\wt w$ such that 
$(t'',\wt w) \in X_0$, 
\begin{equation}\label{t+zeroeq} 
|\wt w'| \le C_0{\varrho}_0^{1/2}{\Delta}^{1/2} H_0^{-1/4}\quad\text{
  and }\quad|\wt w_1| \le C_0{\Delta}
\end{equation}
which gives $|\wt w|\le C_0{\varrho}_0 H_0^{-1/2}$ since ${\Delta} \le
1 \le {\varrho}_0 H_0^{-1/2}$.
\end{claim}

\begin{proof}[Proof of Claim~\ref{lemclaim}] 
We shall consider the cases when ${\delta}_0(t',w_0)$ and
${\delta}_0(t'',w_0)$ have opposite or the same sign.  If
${\delta}_0(t',w_0)$ and ${\delta}_0(t'',w_0)$ have the opposite
sign (including one or more being zero) then we obtain that
$|{\delta}_0(t,w_0)|\le |{\delta}_0(t'',w_0) - {\delta}_0(t',w_0)|
\le 2{\Delta}$ by~\eqref{d0var} for $t = 
t'$ and $t''$. We find from ~\eqref{h0h1est1} that
$|{\delta}_0(t'',w_0)| 
\le {\varrho}_0 H_0^{-1/2} \ll h^{-1/2}$ when ${\varrho}_0 \ll 1$,
and then there must exist a point $(t'',\wt w) \in X_0$ such
that $|\wt w -w_0| \le 2{\Delta}$. Since $|w_0| =
|{\delta}_0(t',w_0)|$ we find $|\wt w|
\le |w_0| + |\wt w -w_0 | \le 4{\Delta}$, thus we
obtain~\eqref{t+zeroeq} in this case.

When ${\delta}_0(t',w_0) = z_1$ and ${\delta}_0(t'',w_0)$ have the same sign
(and are non-zero), we shall first consider the case when they both
are negative.  Then we have that $f(t'',(z_1,0))  < 0$,
and since ${\beta}_-(0) = 0$ we find that
$f(t',(w_1,0)) >0$ for $0 < w_1 \le c_0 H_0^{-1/2}$, which implies
that $f(t'',(w_1,0)) \ge 0$ by condition~ ($\ol{\Psi}$). We find that
$f$ must have a zero at $(t'',\wt w_1,0)$ for some $0 \ge \wt w_1 \ge
z_1 ={\delta}_0(t',w_0)$.  Since $|{\delta}_0(t'',w_0) -
{\delta}_0(t',w_0)| \le 2{\Delta}$, we find $z_1 -\wt w_1 \le
{\delta}_0(t'',w_0) \le z_1 +2 {\Delta}$. Thus $0 \ge\wt w_1 \ge
-2{\Delta}$, which gives ~\eqref{t+zeroeq} in this case with $\wt w'
=0$.

Finally, we consider the case when both ${\delta}_0(t',w_0)= z_1$ and
${\delta}_0(t'',w_0)$ are positive.  By condition~ ($\ol
{\Psi}$) the sign changes of $f$ when $t = t''$ and $|w| \le c_0
H_0^{-1/2}$ are located where $w_1 \le {\beta}_-(w')$. Since
$|{\delta}_0(t'',w_0) - {\delta}_0(t',w_0)| \le 2 {\Delta}$, we
find that ${\delta}_0(t'',w_0) \le z_1 +
2{\Delta}$.  Now let ${\Gamma}_0$ be a circle with radius $r_1
={\delta}_0(t'',w_0) \le {\varrho}_0H_0^{-1/2}$ 
centered at $w_0 = (z_1,0)$, by shrinking ${\varrho}_0$ we may assume
that ${\varrho}_0 \le c_0$. We find that $(w_1,w') \in {\Gamma}_0$
implies that $w_1 = z_1 \pm \sqrt{r_1^2 -|w'|^2} \ge z_1 - r_1 +
|w'|^2/2r_1$ and $|w'| \le r_1$, thus
$$ w_1  \ge 
-2 {\Delta} + \frac{1}{2{\varrho}_0}
H_0^{1/2}|w'|^2\quad\text{and}\quad |w'| \le r_1 \le  {\varrho}_0
H_0^{-1/2}
$$ since $z_1 -r_1 \ge -2{\Delta}$ and $\sqrt{1- t^2} \le 1-t^2/2$ for $|t|
\le 1$. Now 
$|{\beta}_-(w')| \le
C_2 H_0^{1/2}|w'|^2$ for $|w'| \le c_0 H_0^{-1/2}$,
thus ${\Gamma}_0$ will intersect the set $\set {w_1\le {\beta}_-(w')}$
only at points $\wt w = (\wt w_1,\wt w')$ where
\begin{equation*}
-2 {\Delta}  + \frac{1}{2{\varrho}_0}
H_0^{1/2}|w'|^2 \le  C_2 H_0^{1/2}|w'|^2\quad \text{and}\quad |w'| \le
  {\varrho}_0H_0^{-1/2}. 
\end{equation*}
When ${\varrho}_0 \le
(4C_2)^{-1}$ is small enough we find that $|\wt w'| \le \sqrt{8{\varrho}_0}
{\Delta}^{1/2}H_0^{-1/4}$ and   
$|\wt w_1| = |{\beta}_-(\wt w')| \le 8
C_2{\varrho}_0 {\Delta}$. Since there must be sign changes of $w
\mapsto f(t'',w)$ in this set, we find that there exists $(t'', \ol w)
\in X_0$ such that $|\ol w'| \le  \sqrt{8{\varrho}_0}
  {\Delta}^{1/2} H_0^{-1/4}$ and $-2{\Delta}\le \ol w_1
\le |{\beta}_-(\ol w')| \le
8 C_2 {\varrho}_0{\Delta}$. Thus we obtain~\eqref{t+zeroeq} in this
case, which completes the proof of Claim~\ref{lemclaim}.
\end{proof}

It remains to finish the Proof of Lemma~\ref{betalemma}.  By the slow
variation and ~\eqref{t+zeroeq} we find
$ 
H_1^{1/2}(t'',\wt w) \le C^2H_0^{1/2} \le C^2{\varrho}_0
$ 
for ${\varrho}_0$ small enough.
Since $f(t',w)$ has the same sign as $w_1- {\beta}_-(w')$ and $t'' \ge t'$, we
find from condition ($\ol {\Psi}$) that
$ f(t'',w) \ge 0 $ when $w_1 \ge |{\beta}_-(w')|$ and $|w| \le
c_1 H_0^{-1/2}$. By ~\eqref{t+zeroeq}
we find that $|{\beta}_-(w')| \le C_2 H_0^{1/2}|w'|^2 \le
C_3({\varrho}_0{\Delta}+ H_0^{1/2}) \le C'$ when $|w'- \wt w'| \le C$, thus 
$$f(t'',w_1,\wt w') \ge 0\qquad\text{when $C' \le w_1 \ll H_0^{-1/2}$
and $|w'- \wt w'| \le C$.} 
$$ Since $f(t'',\wt w) = 0$ we find by Remark~\ref{estdfrem} that
$|\partial_{w_1}f(t'',\wt w)| \ge c'|\partial_w f(t'',\wt w)|$ for
some $c' >0$.  By using Proposition~\ref{ffactprop} as before at
$(t'',\wt w)$ for small enough ${\varrho}_0$, we obtain
~\eqref{betapmdef} for $t = t''$ with ${\beta}_+(w') \in
S(H_1^{-1/2}(t'', \wt w), H_1(t'', \wt w)g^\sh)$ when $|w'- \wt w'|
\le c_1H_1^{-1/2}(t'',\wt w)$, which contains the neighborhood
$\set{|w'| \le c_0 H_{0}^{-1/2}}$ when $c_0$ and ${\varrho}_0$ are
small enough. Observe that $|{\beta}_+(\wt w')| = |\wt w_1| \le C_0
{\Delta}$ by ~\eqref{t+zeroeq}.

It remains to prove the estimates on ${\beta}_+$. By condition $(\ol
{\Psi})$ we find $\exists\, c >0$ such that
$${\beta}_-(w') -{\beta}_+(w')\ge 0\qquad\text{in $|w'- \wt w'| \le c
H_{0}^{-1/2}$}
$$ 
when ${\varrho}_0 \ll 1$. We also have $\mn{{\beta}_\pm''} \le
C_2H_0^{1/2}$ and $|{\beta}_\pm(\wt w')| \le C{\Delta}$ by
~\eqref{t+zeroeq}, thus we find from Lemma~7.7.2 in ~\cite{ho:yellow}
that 
$$|{\beta}'_-(\wt w') -{\beta}'_+(\wt w')| \le C H_0^{1/4}{\Delta}^{1/2}.
$$ Since ${\beta}_-'(0) =0$ we find that $|{\beta}'_-(\wt w')| \le
C_2C_0{\varrho}_0^{1/2} H_0^{1/4}{\Delta}^{1/2}$ by ~\eqref{t+zeroeq},
and thus $|{\beta}'_+(\wt w')| \le C''
H_0^{1/4}{\Delta}^{1/2}$. By using Taylor's formula, we find that
$|{\beta}'_+(0)| \le C_0 H_0^{1/4}{\Delta}^{1/2}$, $|{\beta}_+(0)|
\le C_0{\Delta}$ and ${\beta}_+ \in S(H_0^{-1/2}, H_0g^\sh)$
when $|w'| \le c_0 H_0^{-1/2}$, for some $C_0 > 0$. This completes the
proof of Lemma~\ref{betalemma}.
\end{proof}

We will also need the following geometrical Lemma  in
Section~\ref{weight}. It tells how far away the minimum of the
distance is attained.

\begin{lem}\label{minlem}
Let ${\Sigma} = \set{w_1 = {\beta}(w')}$, $(w_1,w') \in \br^{1 + n}$,
where ${\beta}(w') \in C^2(\br^n)$, $|{\beta}'(0)| =
{\gamma}H_0^{1/2}$ and $\mn{{\beta}''(w')}\le  H_0^{1/2}$,
$\forall\,w'$, for some positive constant $H_0^{1/2}$, and let
${\delta}(w)$ be the minimal euclidean distance from $w= (w_1,w')$ to
${\Sigma}$.  Then if $|w| \le
{\lambda}H_0^{-1/2}$ and $|{\beta}(0)| \le {\lambda}H_0^{-1/2}$ for
${\lambda} \le 1/6$, then we find ${\delta}(w) = |w-z|$ where
$z =({\beta}(z'),z') \in {\Sigma}$ satisfies
\begin{equation}\label{zwest}
|z'-w'| \le 6{\lambda}({\gamma} + |w'|) \le {\gamma} + |w'|. 
\end{equation}
\end{lem}

\begin{proof}[Proof of Lemma~\ref{minlem}]
Since $|{\beta}(0)| \le {\lambda}H_0^{-1/2}$ we find that
${\delta}(w)\le 2{\lambda}H_0^{-1/2} $  when $|w| \le
{\lambda}H_0^{-1/2}$. We find that
there exists $z = ({\beta}(z'),z')$ such that $|w-z| ={\delta}(w)$ and
$|z'| \le |w| + {\delta}(w) \le 3{\lambda}H_0^{-1/2}$.

Take the right-angled triangle with corners at $w$, $({\beta}(z'),z')$
and $({\beta}(z'),w')$.  By Pythagoras' theorem we find that
\begin{equation}\label{cordeq}
|z' - w'|^2 = {\delta}^2(w) - |w_1 - {\beta}(z')|^2 = ({\delta}(w) +
|w_1 - {\beta}(z')|)({\delta}(w) -|w_1 - {\beta}(z')|).
\end{equation}
Since $|w_1 - {\beta}(z')|\le {\delta}(w) \le |w_1- {\beta}(w')|$ we
find by the triangle inequality
$$
|{\delta}(w) -|w_1 - {\beta}(z')|| \le | {\beta}(z') -{\beta}(w')|
$$ 
By Taylor's formula we obtain 
\begin{multline*}
 |{\beta}(z') - {\beta}(w')| \le |{\beta}'(w')||z'-w'| +
 |z'-w'|^2H_0^{1/2}/2 \\\le |z'-w'|({\gamma} + |w'| +
 |z'-w'|/2 ) H_0^{1/2},
\end{multline*}
since $|{\beta}'(w')| \le |{\beta}'(0)| + |w'|  H_0^{1/2}
\le  ({\gamma} + |w'|) H_0^{1/2}$, thus
$$
|{\delta}(w) -|w_1 - {\beta}(z')||
\le ({\gamma} + |w'| + |z'-w'|/2 )|z'-w'| H_0^{1/2}.
$$ 
Since $|w_1 - {\beta}(z')|\le {\delta}(w) \le 2{\lambda}
H_0^{-1/2}$ we obtain from~ \eqref{cordeq} that
\begin{equation*}
|z'- w'|^2 \le 4{\lambda}({\gamma} + |w'| +
|z'-w'|/2 )|z'-w'|.
\end{equation*}
Thus, for  ${\lambda} \le 1/6$ we
find 
$|z'-w' | \le 6{\lambda}({\gamma} + |w'|)$
which proves the Lemma.
\end{proof}

\section{The Weight function}\label{weight}

In this section, we shall define the weight $m_{{\varrho}}$ we shall
use, it will depend on a parameter $0 <{\varrho}\le 1$.  The weight
will essentially measure how much $t \mapsto {\delta}_0(t,w)$ changes
between the minima of $t \mapsto H_1^{1/2}(t,w)$.  Since $H_1^{1/2}$
gives an upper bound on the curvature of the zero set when $H_1^{1/2}
\ll 1$, the weight will give a bound on the sign changes of the symbol
similar to condition~\eqref{newlemass0} in suitable coordinates.

Recall that $t \mapsto {\delta}_0(t,w)$ and $t \mapsto H_1^{1/2}(t,w)$
are regulated functions, and as before we shall assume that they are
constant when $|t| \ge 1$. 
In the
following, we let $\w{s} = 1 + |s|$.

\begin{defn}\label{h0def}
For $0 < {\varrho} \le 1$ and  $(t_0,w_0) \in \br\times T^*\br^n$ we
define 
\begin{align}\label{m-def}
&m_{-,{\varrho}}(t_0,w_0) = \inf_{t' \le t_0}
\set{{\varrho}^2({\delta}_0(t_0,w_0) -{\delta}_0(t',w_0)) +
  H_1^{1/2}(t',w_0)\w{{\varrho} {\delta}_0(t',w_0)}}\\
&m_{+,{\varrho}}(t_0,w_0) = \inf_{t_0 \le t''}\label{m+def}
\set{{\varrho}^2({\delta}_0(t'',w_0) -{\delta}_0(t_0,w_0)) +
  H_1^{1/2}(t'',w_0)\w{{\varrho} {\delta}_0(t'',w_0)}}
\end{align}
and
\begin{equation} \label{mphodef}
m_{{\varrho}} = \min(\max(m_{+,{\varrho}},
m_{-,{\varrho}}),{\varrho}^2).
\end{equation}
\end{defn}

Then we have $ch^{1/2} \le m_{\pm,{\varrho}} \le H_1^{1/2}
\w{{\varrho} {\delta}_0} \le 1$ 
by \eqref{h2def}--\eqref{H1hest}. 
We find that $m_1 = \max(m_{+,1},m_{-,1})$ and 
\begin{equation} \label{hhhest} 
\min(c h^{1/2}, {\varrho}^2) \le m_{{\varrho}} \le
\min(H_1^{1/2}\w{{\varrho}{\delta}_0}, {\varrho}^2).
\end{equation} 
Now we have 
\begin{multline*}
 m_1(t_0,w_0) \cong \inf_{t' \le t_0 \le
    t''}\Big\{\,{\delta}_0(t'',w_0) -{\delta}_0(t',w_0)\\ +
    H_1^{1/2}(t',w_0)\w{{\delta}_0(t',w_0) }  +
    H_1^{1/2}(t'',w_0)\w{{\delta}_0(t'',w_0) }\,\Big\}
\end{multline*}
and thus $m_1(t_0,w_0) \cong 1$ when $|{\delta}_0(t,w_0)|
\cong H_1^{-1/2}(t,w_0)$ for $t \ge t_0$ or $t \le t_0$.  When $t
\mapsto {\delta}_0(t,w_0)$ is constant, we find that $m_{{\varrho}}$
is proportional to the quasi-convex hull of $t \mapsto
H_1^{1/2}(t,w_0)$ (i.e., it is convex with respect to the constant
functions). The weight also has the ``convexity property'' given by
Proposition~\ref{qmaxpropo}: if $\max_I m_{1} \gg \min_I m_{1}$ on $I
= \set{(t,w): \ a \le t \le b}$, then $\exists\ c > 0$ so that the
variation in $t$ of ${\delta}_0$ on ~~$I$ is bounded from below:
$|{\Delta}_I{\delta}_0| \ge c\max_I m_{1}$.  We shall use the
parameter ${\varrho}$ to obtain suitable norms in
Section~\ref{norm}, but this is just a technicality: all $m_{\varrho}$
are equivalent according to Proposition~\ref{weightprop}.  
Next, we shall show that the conditions in Lemma
~\ref{betalemma} are obtained for small enough $m_{\varrho}$.

\begin{prop}\label{weightprop}
If ${\varrho}=1$ or\/  $m_{{\varrho}}(t_0,w_0) <
{\varrho}^2 < 1$, then there exist $t' \le
t_0 \le t''$ such that
\begin{align}\label{6.6}
&|{\delta}_0(t,w_0) -{\delta}_0(t_0,w_0)| < {\varrho}^{-2}
m_{{\varrho}}(t_0,w_0) \le 1\\ 
&H^{1/2}_1(t,w_0)\w{{\varrho}{\delta}_0(t,w_0)} <
2 m_{{\varrho}}(t_0,w_0) \le 2{\varrho}^2.
\label{6.7}
\end{align}
for $t = t'$ and $t''$. 
The function $t \mapsto m_{{\varrho}}(t,w)$ is
regulated such that
\begin{equation}\label{h0eq}
{\varrho}_1^2/{\varrho}_2^2 \le m_{{\varrho}_1}(t,w) /
m_{{\varrho}_2}(t,w) \le 1
\end{equation}
when $0 <{\varrho}_1 \le {\varrho}_2 \le 1$.
\end{prop}

We obtain from Proposition~\ref{weightprop} that $H^{1/2}_1(t,w_0) <
2{\varrho}^2$ and $|{\delta}_0(t,w_0)| \le 2{\varrho}
H_1^{-1/2}(t,w_0)$ when $t = t'$, $t''$. Thus, when
$m_{{\varrho}}(t_0,w_0) < {\varrho}^2 \ll 1$ we may use
Proposition~\ref{ffactprop} at $(t',w_0)$ and
$(t'',w_0)$. 
Observe that when $m_{{\varrho}}(t_0,w_0) <
{\varrho}^2$ (or ${\varrho}=1$) we obtain from ~\eqref{6.6} that
\begin{equation}\label{6.11a}
1/2 \le \w{{\varrho}{\delta}_0(
  t,w_0)}/\w{{\varrho}{\delta}_0(t_{0},w_0)} \le 2  \qquad t =
t',\ t''
\end{equation}
which gives
$ 
 H^{1/2}_1(t,w_0)\le 4H^{1/2}_1(t_{0},w_0)$ for $ t =
t',\ t''
$, 
since then $H^{1/2}_1(t,w_0)\w{{\varrho}{\delta}_0(t,w_0)} <
2m_{{\varrho}}(t_0,w_0)\le
2H^{1/2}_1(t_0,w_0)\w{{\varrho}{\delta}_0(t_0,w_0)}$ by ~\eqref{hhhest}
and ~\eqref{6.7}. 
Thus if
$m_{{\varrho}}(t_0,w_0) < {\varrho}^2$ (or ${\varrho}=1$) we have
\begin{equation}\label{H1H0est}
 H^{-1/2}_1(t_0,w_0)\le 4 \min(H^{-1/2}_1(t',w_0),\,
 H_1^{-1/2}(t'',w_0)).
\end{equation}
As in ~\eqref{h0small} we shall in the following 
denote $H_0^{1/2} = \max(H^{1/2}_1(t',w_0),\, H_1^{1/2}(t'',w_0))$,
see for example ~\eqref{h0deforig}.

\begin{proof}
We have that  $m_{\pm,{\varrho}} \le m_{{\varrho}}$ when 
$m_{{\varrho}} < {\varrho}^2< 1$ or when ${\varrho}=1$.
By approximating the limit, we may choose $t'' \ge t_0$ so that
\begin{equation} \label{limapprox}
{\varrho}^2({\delta}_0(t'',w_0) -{\delta}_0(t_0,w_0)) + H_1^{1/2}(
t'',w_0) \w{{\varrho}{\delta}_0(t'',w_0)}
< m_{+,{\varrho}}(t_0,w_0) + ch^{1/2}
\end{equation} 
where $c$ is chosen as in ~\eqref{hhhest}. Then we find $
{\varrho}^{2}({\delta}_0(t'',w_0) -{\delta}_0(t_0,w_0)) <
m_{+,{\varrho}}(t_0,w_0) $ and
$H_1^{1/2}(t'',w_0)\w{{\varrho}{\delta}_0(t'',w_0)} <
m_{+,{\varrho}}(t_0,w_0) + ch^{1/2} \le 2 m_{+,{\varrho}}(t_0,w_0)$.
We similarly obtain this estimate for $m_{-,{\varrho}}$ with $t' \le t_0$,
which gives ~\eqref{6.6}--\eqref{6.7}.

To prove ~\eqref{h0eq} we let $F_{\varrho}(t,s,w) =
{\varrho}^2|{\delta}_0(s,w) -{\delta}_0(t,w)| + H_1^{1/2}(s,w)
\w{{\varrho}{\delta}_0(s,w)}$. Then we have $F_{\varrho_1} \le
F_{\varrho_2}$ and 
${\varrho}_1^2F_{{\varrho}_2} \le {\varrho}_2^2 F_{{\varrho}_1}$ when
${\varrho}_1 \le {\varrho}_2$. Since these estimates are preserved
when taking infimum and supremum, we obtain ~\eqref{h0eq} for
$m_{\pm,{\varrho}_j}$ and  $m_{{\varrho}_j}$, $j =1$, $2$.

To prove that $t \mapsto m_{\varrho}(t,w)$ is a regulated function, it
suffices to prove that $t \mapsto m_{\pm,\varrho}(t,w)$ is a regulated
function since this property is preserved when taking maximum and
minimum.  We note that
$$t \mapsto m_{+,{\varrho}}(t,w_0) = \inf_{t \le t''}
\set{{\varrho}^2{\delta}_0(t'',w_0) +
  H_1^{1/2}(t'',w_0)\w{{\varrho}{\delta}_0(t'',w_0)}}
-{\varrho}^2{\delta}_0(t,w_0)$$
and since the infimum is non-decreasing and bounded, we find that this
gives a regulated function in ~$t$. A similar argument works for
$m_{-,{\varrho}}$, which proves the result.
\end{proof}

In the following we shall assume the coordinates chosen so that
$g^\sh(w) = |w|^2$. Observe that $m_{\varrho}$ is not a weight for
$G_1$, but the following Proposition shows that it is a weight for
$g_{\varrho} = {\varrho}^2g^\sh$ uniformly in ${\varrho}$.

\begin{prop}\label{h0slow}
We find that there exists $C>0$  such that
\begin{equation}\label{wtemp}
m_{{\varrho}}(t,w) \le C m_{{\varrho}}(t,w_0)(1 +
{\varrho}^2g^\sh(w-w_0))\qquad \forall\,t
\end{equation}
uniformly when $0 <{\varrho}\le 1$, which implies
that $m_{{\varrho}}$ is a weight
for $g_{\varrho} = {\varrho}^2g^\sh$ since $g_{\varrho}^{\sigma} =
g^\sh/{\varrho}^{2}$.  The constant $C$
only depends on the seminorms of $f$ in $S(h^{-1},hg^\sh)$.
\end{prop}

\begin{proof}
Since $m_{{\varrho}} \le {\varrho}^2$ we only have to
consider the case when 
\begin{equation} \label{h0rhoest}
m_{{\varrho}}(t_0,w_0) < {\varrho}^2.
\end{equation} 
Now, it suffices to show that
\begin{equation}\label{6.25}
m_{{\varrho}}(t_0,w)/m_{{\varrho}}(t_0,w_0) \le
 C(1 + {\varrho}^2|w-w_0|^2)\qquad\text{when $|w-w_0|\le {\varrho}
   m_{{\varrho}}^{-1}(t_0,w_0)$}
\end{equation} 
uniformly in $0 <{\varrho} \le 1$. In fact, when $|w-w_0| > {\varrho}
m_{{\varrho}}^{-1}(t_0,w_0)$ we obtain that ${\varrho}^2|w-w_0|^2 >
{\varrho}^4 m_{{\varrho}}^{-2}(t_0,w_0) > m_{{\varrho}}(t_0,w)/
m_{{\varrho}}(t_0,w_0)$ by~\eqref{h0rhoest}. Thus ~\eqref{wtemp} is
trivially satisfied with $C=1$ when $|w-w_0| > {\varrho}
m_{{\varrho}}^{-1}(t_0,w_0)$. Thus, in the following we shall only
consider $w$ such that $|w-w_0|\le {\varrho}
m_{{\varrho}}^{-1}(t_0,w_0)$.

Since $m_{{\varrho}}$ are equivalent when ${\varrho} \ge {\varrho}_0 >0$
by~\eqref{h0eq}, it suffices to consider ${\varrho} \le {\varrho}_0
\ll 1$.  In fact, if ~\eqref{wtemp}
holds for $m_{{\varrho}_0}$ then it holds for $m_{{\varrho}}$ when
${\varrho}_0 \le {\varrho} \le 1$, with
$C$ replaced by $C/{\varrho}_0^2$. 
In the following we shall assume $0 < {\varrho} \le {\varrho}_0$,
where ${\varrho}_0$ shall be determined later.
Since we assume~\eqref{h0rhoest} we may use 
Proposition~\ref{weightprop} to obtain $t' \le t_0 \le t''$
such that \eqref{6.6}--\eqref{6.7} hold. 
By~\eqref{6.7} and~\eqref{6.11a} we obtain that
\begin{equation}\label{delta0ref}
|{\delta}_0(t, w_0)| \le 8{\varrho} H_0^{-1/2} \qquad \text{for
  $t = t'$, $t''$}
\end{equation} 
where 
\begin{equation}\label{h0deforig}
H_0^{1/2} =
\max(H_1^{1/2}(t',w_0),H_1^{1/2}(t'',w_0)) <
2m_{{\varrho}}(t_0,w_0) < 2 {\varrho}^2,
\end{equation}
which gives ${\varrho{}}m_{{\varrho}}^{-1}(t_0,w_0) < 
2{\varrho{}}H_0^{-1/2}$. Since
\begin{equation} \label{Deltadef}
{\Delta} = \max_{t = t',\ t''} |{\delta}_0(t,w_0) -{\delta}_0(t_0,w_0)|
 < {\varrho}^{-2}m_{{\varrho}}(t_0,w_0) < 1
\end{equation} 
we may use Lemma~\ref{betalemma} when $2{\varrho} < 8{\varrho}^2 \ll 1$
to obtain $g^\sh$ ~orthonormal coordinates $w= (w_1,w')$ so that $w_0
= (z_1,0)$, $z_1 = {\delta}_0(t',w_0)$ and
\begin{equation*}
\left\{
\begin{aligned} 
&\sgn f(t',w) = \sgn(w_1-{\beta}_-(w'))\\
&\sgn f(t'',w) = \sgn(w_1-{\beta}_+(w'))
\end{aligned}\right.\qquad \text{for $|w|\le
c_0 H_0^{-1/2}$}
\end{equation*}
where $w' \mapsto {\beta}_\pm(w') \in S(H_0^{-1/2},H_0g^\sh)$. We also find
${\beta}_-(0) = |{\beta}'_-(0)| =0$, and since ${\Delta} <
{\varrho}^{-2}m_{{\varrho}}(t_0,w_0)$ by ~\eqref{Deltadef} we obtain that
\begin{equation}\label{beta+est} 
|{\beta}_+(0)| \le C{\varrho}^{-2}m_{{\varrho}}(t_0,w_0)\qquad \text{and}
\qquad|{\beta}'_+(0)| \le
C{\varrho}^{-1}H_0^{1/4}m_{{\varrho}}^{1/2}(t_0,w_0).
\end{equation}
Since $|w_0| = |{\delta}_0(t',w_0)| \le 8{\varrho} H_0^{-1/2}$
by ~\eqref{delta0ref} and $ |w-w_0| \le
{\varrho{}}m_{{\varrho}}^{-1}(t_0,w_0) <
2{\varrho{}}H_0^{-1/2}$ by ~\eqref{h0deforig}, we have $|w| \le
10{\varrho} H_0^{-1/2}$ and 
$$|{\delta}_0(t,w)|
\le |{\delta}_0(t,w_0)| + |w-w_0| \le  10{\varrho}
H_0^{-1/2}  \le  10{\varrho} h^{-1/2}/c \qquad t=t',\ t''
$$ 
by the uniform Lipschitz continuity of $w \mapsto {\delta}_0(t,w)$. 
Thus when
${\varrho}  \ll 1$  we find that there exists
$(t'',{\beta}_+(z'), z') \in X_0$ for some $|z'| \le 20{\varrho}
H_0^{-1/2}$
such that $|w - ({\beta}_+(z'), z')| = |{\delta}_0(t'',w)|$.

In the following, we let $m_{{\varrho}}= m_{{\varrho}}(t_0,w_0)$.
Now $|{\beta}_+(0)| \le C{\varrho}^{-2}m_{{\varrho}} < C <
2C{\varrho}^2H_0^{-1/2}$, $|{\beta}'_+(0)| \le
C{\varrho}^{-1}H_0^{1/4}m_{{\varrho}}^{1/2}$ and $\mn{{\beta}_+''} \le
C_2H_0^{1/2}$ by ~\eqref{beta+est}. Thus, when $|w| \le
10{\varrho} H_0^{-1/2}$ we may use
Lemma~\ref{minlem} with ${\gamma} = 
CC_2^{-1}{\varrho}^{-1}H_0^{-1/4}m_{{\varrho}}^{1/2}$ and
${\lambda} = C_2\min(2C{\varrho}^2,10{\varrho}) \le 1/6$ when
${\varrho}^2 < 10 {\varrho} \ll 1$ (since $H_0^{1/2}$ is replaced by
$C_2H_0^{1/2}$).   
 Thus, we obtain that
\begin{equation*}
|z'-w'| \le CC_2^{-1}{\varrho}^{-1}H_0^{-1/4}m_{{\varrho}}^{1/2}
+ |w'| 
\end{equation*}
which gives that $|z'| \le
C C_2^{-1}{\varrho}^{-1}H_0^{-1/4}m_{{\varrho}}^{1/2} + 2|w'|$. Since
$|{\beta}_+(0)| \le 
C{\varrho}^{-2}m_{{\varrho}}$ and $|{\beta}'_+(0)| \le
C{\varrho}^{-1}H_0^{1/4}m_{{\varrho}}^{1/2}$ we find
that for these $z'$ and $w'$ we have that
$$|{\beta}_+(z')| + |{\beta}_+(w')| \le 
C_0 ({\varrho}^{-2}m_{{\varrho}} + H_0^{1/2}|w'|^2)$$
by using Taylor's formula and the Cauchy-Schwarz inequality. 
Since $|w_1 -
{\beta}_+(z')|\le |{\delta}_0(t'',w)| \le |w_1 -{\beta}_+(w')|$ we
obtain that
\begin{equation*}
|w_1| - C_0({\varrho}^{-2}m_{{\varrho}} + H_0^{1/2}|w'|^2)
\le |{\delta}_0(t'',w)| \le |w_1| +
C_0 ({\varrho}^{-2}m_{{\varrho}} + H_0^{1/2}|w'|^2)
\end{equation*}
when $|w-w_0|\le 2{\varrho}H_0^{-1/2}$ and ${\varrho} \le
{\varrho}_0 \ll 1$. Since $|{\beta}'_-(0)| = 
{\beta}_-(0) = 0$ a similar argument (with ${\gamma} =0$ and thus
$|z'| \le 2|w'|$) gives 
$$
||{\delta}_0(t',w)|- |w_1|| \le C_1H_0^{1/2}  |w'|^2\qquad
\text{when $|w-w_0|\le 2{\varrho{}}H_0^{-1/2}$}
$$ 
for some $C_1 > 0$ when ${\varrho}\le {\varrho}_0 \ll 1$.
We obtain that 
\begin{equation}\label{dpmvar}
{\delta}_0(t'',w) - {\delta}_0(t',w) \le C_2
({\varrho}^{-2}m_{{\varrho}} + H_0^{1/2}|w'|^2) 
\qquad \text{when $|w-w_0|
\le 2{\varrho{}}H_0^{-1/2}$}.
\end{equation}
In fact, if ${\delta}_0(t'',w)$ and ${\delta}_0(t',w)$ have the same
sign then we find $|{\delta}_0(t'',w) - {\delta}_0(t',w) |=
 ||{\delta}_0(t',w)| - |{\delta}_0(t'',w)||$. Moreover,
${\delta}_0(t'',w)$ and ${\delta}_0(t',w)$ may only have
different signs when ${\beta}_-(w') \le w_1 \le {\beta}_+(w')$ or
${\beta}_+(w') \le w_1 \le {\beta}_-(w')$, so we find in this case that
$$|{\delta}_0(t,w)| \le |{\beta}_+(w') -{\beta}_-(w') | \le
C ({\varrho}^{-2}m_{{\varrho}} + H_0^{1/2}|w'|^2)\qquad t = t',\ t''.$$ 
We obtain from
~\eqref{dpmvar} and the monotonicity of $t \mapsto {\delta}_0(t,w)$ 
that 
\begin{multline}\label{dpmvar1}
{\varrho}^2|{\delta}_0(t,w) - {\delta}_0(t_0,w)| \le
{\varrho}^2({\delta}_0(t'',w) - {\delta}_0(t',w))\\ 
\le C_2 ( m_{{\varrho}} + H_0^{1/2}{\varrho}^2|w'|^2)\le 2C_2
m_{{\varrho}}  (1 + {\varrho}^2|w-w_0|^2)
\end{multline}
when $t = t'$, $t''$ and $|w-w_0| \le {\varrho}m_{\varrho}^{-1} <
2{\varrho{}}H_0^{-1/2}$, since 
$|w'| \le |w-w_0|$.  Now $G_1$ is slowly varying, thus we find for
small enough ${\varrho{}} >0$ that
\begin{equation*}
H_1^{1/2}(t,w) \le C_3H_1^{1/2}(t,w_0)
 \qquad\text{when  $|w-w_0| \le
   2{\varrho{}}H_0^{-1/2}  \le 
   2{\varrho{}}H_1^{-1/2}(t,w_0)$}
\end{equation*}
for  $t = t'$, $t''$.
By the uniform Lipschitz continuity we find that
\begin{equation}\label{6.11}
\w{{\varrho}{\delta}_0(t,w) }\le
\w{{\varrho}{\delta}_0(t,w_0) }(1 + {\varrho}|w-w_0|)\qquad \forall\,t
\end{equation} 
which implies that for  $t = t'$, $t''$, we have
\begin{equation}\label{h1est}
H_1^{1/2}(t,w) \w{{\varrho}{\delta}_0(t,w) }  \le
C_3H_1^{1/2}(t,w_0) \w{{\varrho}{\delta}_0(t,w_0) } (1 +
{\varrho}|w-w_0|)
\end{equation}
when $|w-w_0| \le 2{\varrho{}}H_0^{-1/2}.$ By using ~
\eqref{6.7}, \eqref{dpmvar1}, \eqref{h1est} and taking the infimum we
obtain
\begin{multline*}
m_{\pm,{\varrho}}(t_0,w) \le C_4 m_{{\varrho}}(t_0,w_0)(1 +
{\varrho}|w-w_0|)^2 \\
\text{when $|w-w_0| \le
  {\varrho{}} m_{{\varrho}}^{-1} (t_0,w_0) \le
2{\varrho{}}H_0^{-1/2}$}
\end{multline*}
uniformly for small ${\varrho}$. By taking the maximum and then the minimum, we
obtain~\eqref{6.25} and thus Proposition~\ref{h0slow}.
\end{proof}

In section~\ref{norm}, we shall choose a fixed ${\varrho} = {\varrho}_0
\ll 1$ in order to get invertible operators and suitable norms.  In
the following, we shall for simplicity only consider $m_{1}$,
since all the $m_{{\varrho}}$ are equivalent when ${\varrho} \ge c
>0$, this is really no restriction: the following results also
holds for any $m_{{\varrho}}$, $0 <{\varrho} \le 1$ but with
constants depending on ${\varrho}$.
The next Proposition shows that if $m_1 \ll 1$ then the conditions in
Proposition~\ref{symbest} are satisfied.

\begin{prop}\label{h0sign}
There exists $0 < {\kappa}_0 < 1$ and $c_0 >0$ such that if  $m_{1}
\le {\kappa}_0$ at 
$(t_0,w_0) \in \br \times T^*\br ^n$, then there exist $g^\sh$
~~orthonormal coordinates so that $w_0 = (z_1,0)$, $|z_1| <
|{\delta}_0(t_0,w_0)| +  1$ and
\begin{equation}\label{signref}
\sgn(w_1)f(t_0,w) \ge 0\quad\text{when $|w_1| \ge (1 +
H_0^{1/2}|w'|^2)/c_0$ and $|w| \le c_0 H_0^{-1/2}$}
\end{equation}
where $c_0h^{1/2} \le H_0^{1/2} <
4 m_{1}(t_0,w_0)\w{{\kappa}_0^{1/2}{\delta}_0(t_0,w_0)}^{-1} \le
4 {\kappa}_0$. 
Here ${\kappa}_0$ and $c_0$ only
depend on the seminorms of $f$ in $S(h^{-1},hg^\sh)$.
\end{prop}

\begin{proof}  
Let ${\varrho}_1 = {\varrho}_0/8$ for the fixed ${\varrho}_0 \le 1$ in
Lemma~\ref{betalemma} and assume that $m_{1}(t_0,w_0) \le {\kappa}_0 <
{\varrho}_1^2$. Since $m_{{\varrho}_1} \le m_{1}$ we can use
Proposition~\ref{weightprop} to find $t'\le t_0 \le t''$ such that
${\Delta} = \max_{t=t',\, t''}|{\delta}_0(t,w_0) -
{\delta}_0(t_0,w_0)| < {\varrho}_1^{-2}m_{{\varrho}_1}(t_0,w_0) <
1$, $|{\delta}_0(t,w_0)| \le 2 {\varrho}_1H_1^{-1/2}(t,w_0)$ for
$t=t'$, $t''$ and 
\begin{equation}\label{H0def} 
H_0^{1/2} =
\max(H_1^{1/2}(t',w_0),H_1^{1/2}(t'',w_0)) < 2 m_{{\varrho}_1} \le
2{\varrho}_1^2. 
\end{equation} 
By using~\eqref{6.11a} as before we find that $|{\delta}_0(t,w_0)| \le
8 {\varrho}_1H_0^{-1/2}$ for $t=t'$, $t''$.

Since $2{\varrho}_1^2 < 8{\varrho}_1 \le {\varrho}_0$ we may use
Lemma~\ref{betalemma} to obtain 
$g^\sh$ ~~orthonormal coordinates so that $w_0 = (z_1,0)$ with $|z_1| =
|{\delta}_0(t',w_0)| 
< |{\delta}_0(t_0,w_0)| + 1$
and
\begin{equation}\label{sgnfeq} 
\left\{
\begin{aligned} 
&\sgn f(t',w) = \sgn(w_1-{\beta}_-(w'))\\
&\sgn f(t'',w) = \sgn(w_1-{\beta}_+(w'))
\end{aligned}\right.\qquad \text{for $|w|\le
c_0 H_0^{-1/2}$.}
\end{equation} 
Here ${\beta}_\pm(w') \in S(H_0^{-1/2}, H_0 g^\sh)$ and
$$
ch^{1/2} \le H_0^{1/2} <
4 m_{1}(t_0,w_0)/\w{{\kappa}_0^{1/2}{\delta}_0(t_0,w_0)}  \le 4 {\kappa}_0.
$$ In fact, $H_1^{1/2}(t,w_0)\w{{\varrho}_1{\delta}_0(t,w_0)} 
< 2m_{{\varrho}_1}(t_0,w_0) \le 2m_1(t_0,w_0)$ 
and $\w{{\kappa}_0^{1/2}{\delta}_0(t_0,w_0)} \le 2
\w{{\varrho}_1{\delta}_0(t,w_0)}$ when $t = t'$ and $t''$
by~\eqref{6.7} and ~\eqref{6.11a}.  Since ${\Delta} \le 1$ we also
obtain from Lemma~\ref{betalemma} that $|{\beta}_\pm(0)|\le C$, $
|{\beta}'_\pm(0)| \le C H_{0}^{1/4}$ and $ \mn{{\beta}_\pm''} \le
C_0 H_{0}^{1/2}$.  This gives
$$|{\beta}_\pm(w')| \le C(1 + H_{0}^{1/2}|w'|^2) \qquad\text{in
$|w'| \le c_0 H_{0}^{-1/2}$}
$$ 
for some $C >0$. By condition ($\ol {\Psi}$) and~\eqref{sgnfeq}, we
obtain the result.
\end{proof}

In order to get lower bounds in terms of the weight $m_{1}$, we 
shall need the following result, which will be important for the proof. 

\begin{prop}\label{mestprop}
Let the weight $M$ be given by Definition~\ref{Mdef}. Then there
exists $C_0>0$ such that
\begin{equation}\label{Mest0}
MH_1^{3/2}\w{{\delta}_0} \le C_0m_{1}
\end{equation} 
which gives $S(MH_1^{3/2} ,G_1) \subseteq
S(m_{1}\w{{\delta}_0}^{-1},g^\sh)$. 
Here $C_0$ only depends on the seminorms of $f$ in $S(h^{-1},hg^\sh)$. 
\end{prop}

\begin{proof}[Proof of Proposition~\ref{mestprop}]
We shall omit the dependence on $t$ in the proof and put $m_{1}=
m_{1}(w_0)$. First we observe that 
if $m_{1} \ge c >0$, then $MH_1^{3/2}\w{{\delta}_0} \le C
\le Cm_{1}/c$ at $w_0$
since $\w{{\delta}_0} \le H_1^{-1/2}$ and $M \le CH_1^{-1}$ by
Proposition~\ref{mproplem}. 

Thus, we only have to consider the case $m_{1} \le {\kappa}_0 \ll
1$.  By using Proposition~\ref{h0sign} for ${\kappa}_0
\ll 1$ we obtain  coordinates so that  $|w_0| \le
|{\delta}_0(w_0)| + 1 \le H_1^{-1/2}(w_0)$, $f$ satisfies
~\eqref{signref} and thus the
conditions in Proposition~\ref{symbest}.
Since ${\kappa}_0^{1/2}\w{{\delta}_0} \le
\w{{\kappa}_0^{1/2}{\delta}_0}$ we obtain from 
Propositions~\ref{symbest} and ~\ref{h0sign}
the estimate
\begin{equation} \label{mhest}
M(0)H_1^{3/2}(0) \le C_1 H_0^{1/2}  \le 4C_1
{\kappa}_0^{-1/2}m_{1}/\w{{\delta}_0(w_0)}.  
\end{equation}
It remains to prove the estimate $M(w_0)H_1^{3/2}(w_0)\le
CM(0)H_1^{3/2}(0)$ in this case. 
By Proposition~\ref{g1prop} we have that
\begin{equation}\label{6.22}
 M(w_0) \le CM(0)(1+ H_1^{1/2}(0)|w_0|)^3
\end{equation}
and 
\begin{equation}\label{6.23}
 H(w_0) \le CH(0)(1+ H_1^{1/2}(w_0)|w_0|)^2.
\end{equation}
In the case $H_1^{1/2}(0) \le H_1^{1/2}(w_0)$ we find that $|w_0| \le
H_1^{-1/2}(w_0) \le H_1^{-1/2}(0)$ and thus $M(w_0)H_1^{3/2}(w_0) \le
64C^{5/2} M(0)H_1^{3/2}(0)$ by ~\eqref{6.22}--\eqref{6.23}. When
$H_1^{1/2}(w_0) \le H_1^{1/2}(0)$ we don't have to use ~\eqref{6.23},
instead we find from ~\eqref{6.22} that
\begin{equation*}
M(w_0)H_1^{3/2}(w_0) \le CM(0)H_1^{3/2}(0)(H_1^{1/2}(w_0)H_1^{-1/2}(0) +
 1)^3 \le 8C M(0)H_1^{3/2}(0)  
\end{equation*}
since $|w_0| \le H_1^{-1/2}(w_0)$. This completes the proof of the
Proposition.
\end{proof}

If $m_1 \cong 1$ then we find from the proof that the estimate
~\eqref{Mest0} is trivial. When $m_1 \ll 1$ we have the following
``geometrical'' interpretation of~\eqref{Mest0}.

\begin{rem} 
In the case $|{\delta}_0| \le C$ we find from
Proposition~\ref{mproplem} that 
$MH_1^{3/2}\w{{\delta}_0} \cong \mn{f''}H_1^{1/2} +
h^{1/2} \le Cm_1$ if and only if $\mn{f''} \le C_1m_1 H_1^{-1/2}$. 
By Proposition~\ref{weightprop}
there exist $t' \le t_0 \le t''$ so that 
$m(t_0,w_0) \cong H_0^{1/2} = \max(H_1^{1/2}(t',w_0),\,
H_1^{1/2}(t'',w_0))$ when $|{\delta}_0(t_0,w_0)| \le C$. 
Then~\eqref{Mcomp} holds at\/ $(t_0,w_0)$  if and only if 
\begin{equation}\label{Mest1}
\mn{f''(t_0,w_0)} \le C_2H_0^{1/2} H_1^{-1/2}(t_0,w_0). 
\end{equation} 
Now it follows from ~\eqref{dfest0} and~\eqref{Mest1}
that $H_0^{-1/2} \cong m_1^{-1}$ has the property that $F =
H_0^{-1/2}f \in S(H_1^{-3/2}, G_1)$ when
$t = t_0$ and $|w - w_0| \ll H_1^{-1/2}$. 
\end{rem}

Next, we shall prove some results about the
properties of $m_{1}$. First we shall prove the ``convexity property''
mentioned earlier.

\begin{prop}\label{qmaxpropo}
Let $m_{1}$ be given by Definition~\ref{h0def}. There exist
${\kappa}_0 >1$, $c_0 >0$ and ${\varepsilon}_0 >0$ such that if ${\kappa} \ge
{\kappa}_0$, $t' < t_0 < t''$ and
\begin{equation}\label{qmaxwd}
m_{1}(t_0,w_0) = {\kappa}\max(
m_{1}(t',w_0),m_{1}(t'',w_0))
\end{equation} 
then we have
\begin{equation}\label{nyref1}
{\delta}_0(t'',w) -{\delta}_0(t',w) \ge c_0 m_{1}(t_0,w_0) 
= c_0{\kappa}\max(m_{1}(t',w_0),m_{1}(t'',w_0))
\end{equation}
when $|w-w_0|\le {\varepsilon}_0$.
\end{prop}

\begin{proof}
Since $t_0 < t''$ we have by the triangle inequality
\begin{multline*} 
m_{+,1} (t_0,w_0) \le \inf_{t''\le t} ({\delta}_0(t,w_0)
-{\delta}_0(t_0,w_0) + H_1^{1/2}(t)\w{{\delta}_0(t,w_0)}) \\
\le {\delta}_0(t'',w_0) -{\delta}_0(t_0,w_0) +m_{+,1} (t'',w_0)
\end{multline*}
and similarly
\begin{equation*}
m_{-,1} (t_0,w_0) \le 
{\delta}_0(t_0,w_0) -{\delta}_0(t',w_0) +m_{-,1} (t',w_0).
\end{equation*} 
Since $m_{\pm,1} \le m_1$ we find that 
\begin{multline*} 
m_{1}(t_0,w_0) = \max(m_{-,1} (t_0,w_0), m_{+,1} (t_0,w_0)) \\ \le
{\delta}_0(t'',w_0)-{\delta}_0(t',w_0) + \max(
m_{1}(t',w_0),m_{1}(t'',w_0))
\end{multline*}
which gives~\eqref{nyref1} for $w= w_0$ with ${\kappa}_0 = 2$ and  $c_0 =
1/2$. 

If we choose
${\varepsilon}_0>0$ so that 
\begin{equation} \label{h0var}
1/C_0 \le m_{1}(t,w)/m_{1}(t,w_0)
\le C_0\qquad\text{when $|w-w_0| \le {\varepsilon}_0$ and $\forall\, t$} 
\end{equation}
then we obtain the result when $|w-w_0| \le {\varepsilon}_0$ for
${\kappa}_0 = 2C_0^2$ and $c_0 = (2C_0)^{-1}$.
In fact, \eqref{qmaxwd} implies that
$
m_{1}(t_0,w) = {\kappa}_1\max(
m_{1}(t',w),m_{1}(t'',w))
$
where ${\kappa}_1 \ge {\kappa}/C_0^{2}$. When ${\kappa} \ge 2C_0^2$ we
find
$
 {\delta}_0(t'',w)
-{\delta}_0(t',w) \ge \frac{1}{2}m_{1}(t_0,w) \ge
\frac{1}{2C_0}m_{1}(t_0,w_0),
$
which proves the Proposition.
\end{proof}

The following Proposition gives an estimate 
on how much $w \mapsto {\delta}_0(t,w)$ varies for different
values of $t$, 
using the monotonicity of $t \mapsto {\delta}_0(t,w)$.

\begin{prop}\label{dvar1} 
Let $m_{1}$ be given by Definition~\ref{h0def} and
let ${\Delta}(s,t,w) = {\delta}_0(t,w) -{\delta}_0(s,w) \ge 0$ for
$s \le t$. There exists ${\varepsilon_1} >0$ so that if
\begin{equation*} 
\max(m_{1}(s_0,w_0), m_{1}(t_0,w_0)) \le K\le 1
\end{equation*}
and 
\begin{equation}\label{ddcond}
|{\Delta}(s_0,t_0,w_0)| = {\lambda}K\qquad {\lambda} \ge 0
\end{equation}
for some $s_0$, $t_0 \in \br$, then
\begin{equation}\label{ddest0}
(2{\lambda}/3 - 3)K \le |{\Delta}(s_0,t_0,w)|
\le  (4{\lambda}/3 +
 3)K\qquad \text{when $|w-w_0| \le {\varepsilon}_1$.}
\end{equation}
Here ${\varepsilon}_1$ does not depend on ${\lambda}$ or $K$.
\end{prop}

\begin{proof} 
To prove ~\eqref{ddest0} it suffices to show that
\begin{equation}\label{ddest}
 |{\Delta}(s_0,t_0,w)- {\Delta}(s_0,t_0,w_0)| \le  ({\lambda}/3 +
 3)K\qquad \text{when $|w-w_0| \le {\varepsilon}$}
\end{equation}
for sufficiently small ${\varepsilon}>0$.
Observe that if
$K \ge {\kappa} >0$ then by the uniform Lipschitz continuity of $w
\mapsto  {\delta}_0(t,w)$, we obtain ~\eqref{ddest} with ${\lambda}
=0$ for small enough ${\varepsilon}$ (depending on~ ${\kappa}$).

Thus, in the following we may assume that $K \le {\kappa} \ll 1$,  and
it is no restriction to assume that $s_0 < t_0$.
Since $m_{1}(t_0,w_0) \le K \ll 1$ we obtain by Proposition~\ref{weightprop}
that there exist $t' \le t_0 \le t''$ such that  
\begin{equation} \label{d0test}
|{\delta}_0(t,w_0) -{\delta}_0(t_0,w_0)| <
m_{1}(t,w_0) \le K \qquad t=t',\ t''
\end{equation}
by ~\eqref{6.6}. We obtain from~\eqref{6.7} that $|{\delta}_0(t,w_0)|
\le 2 K H_1^{-1/2}(t,w_0)$ and $H_1^{1/2}(t,w_0) \le 2m_{1}(t,w_0) \le
2K$ for $t = t'$, $t''$. We similarly obtain~\eqref{d0test} with
~~$t_0$ replaced by ~~$s_0$ and $t'$, $t''$ by $s' \le s_0 \le s''$.
In the following we shall assume that $s = s'$, $s''$ and $t = t'$, $
t''$.  When $K$ is sufficiently small, we obtain from
Proposition~\ref{ffactprop} that ${\delta}_0
\in S(H_1^{-1/2},G_1)$ in a fixed $G_1$ neighborhood of $(t,w_0)$
and $(s,w_0)$. 
By ~\eqref{ddcond} and ~\eqref{d0test} we have that 
\begin{equation}\label{dstest}
 |{\Delta}(s,t,w_0)| \le ({\lambda} + 2)K
\end{equation}
when  $s = s'$, $s''$ and $t = t'$, $t''$.
Then we obtain that $\mn{\partial^2_w{\Delta}( s,t,w)}
\le CK$ when $|w-w_0| \le c$ for some $c >0$.  Since
$\pm{\Delta}( s,t,w) \ge 0$ for a choice of sign, we
find that
\begin{equation*}
 \pm\Big({\Delta}( s,t,w_0) + \w{ w-w_0,
   d_w{\Delta}( s,t,w_0)}\Big) + CK|w-w_0|^2/2  \ge 0
\end{equation*}
when $|w-w_0| \le c$.  By optimizing over $|w-w_0| = c$ and
using~\eqref{dstest},  
we obtain that $|d_w{\Delta}( s,t,w_0)| \le
C({\lambda} + 1)K$.
Since $\mn{\partial^2_w{\Delta}( s,t,w)}
\le CK$ when $|w-w_0| \le c$ we obtain that
\begin{equation} \label{dtest}
|{\Delta}( s,t,w) -{\Delta}( s,t,w_0)|
\le \left(\frac{\lambda}3 + 1\right)K \qquad \text{when
$|w-w_0| \le {\varepsilon} \ll 1$}
\end{equation} for $s = s'$, $s''$ and $t = t'$, $t''$,
where ${\varepsilon}$ does not depend on ${\lambda}$.
Since ${\Delta}(s'',t',w) \le {\Delta}(s_0,t_0,w) \le
{\Delta}(s',t'',w) $ for $s_0 < t_0$,
we find by ~\eqref{d0test} that
\begin{equation*}
{\Delta}(s_0,t_0,w) - {\Delta}(s_0,t_0,w_0) \le  {\Delta}(s',t'',w)
-{\Delta}(s',t'',w_0) + 2K
\end{equation*}
and 
\begin{equation*}
{\Delta}(s_0,t_0,w) - {\Delta}(s_0,t_0,w_0) \ge  {\Delta}(s'',t',w)
-{\Delta}(s'',t',w_0) - 2K.
\end{equation*}
Thus we obtain ~\eqref{ddest}
which completes the proof of the Proposition.
\end{proof}

\section{The Wick quantization}\label{norm}

In order to define the pseudo-sign we shall use the Wick quantization, 
following~\cite[Appendix~B]{de:suff} and~\cite[Section~4]{ln:coh}.
The advantage with using the Wick quantization is that positive
symbols give positive operators.  
We shall
also define the norms we shall use, following~\cite{bc:sob}. 
For $a \in
L^\infty(T^*\br^n)$ we define the Wick quantization:
\begin{equation*}
a^{Wick}(x,D_x)u(x) = \int_{T^*\br^n}a(y,{\eta})
{\Sigma}^w_{y,{\eta}}(x,D_x)u(y)\,dyd{\eta}\qquad u \in  \cal S(\br^n)
\end{equation*}
using the projections ${\Sigma}^w_{y,{\eta}}(x,D_x)$ with symbol
$${\Sigma}_{y,{\eta}}(x,{\xi}) =
{\pi}^{-n}\exp(-g^\sharp(x-y,{\xi}-{\eta})) =
{\pi}^{-n}\exp(-|x-y|^2-|{\xi}-{\eta}|^2).$$
We find that
$a^{Wick}:\ \cal S(\br^n) \mapsto \cal S'(\br^n)$
is symmetric on $\cal S(\br^n)$ if $a$ is real valued,
\begin{equation} \label{poswick}
a \ge 0 \quad \text{in $L^\infty(T^*\br^n)$}\Rightarrow
\w{a^{Wick}(x,D_x)u,u} \ge
0 \quad \text{for $u \in \cal S(\br^n)$}
\end{equation}
and
$ 
\mn{a^{Wick}(x,D_x)}_{\cal L(L^2(\br^n))} \le \mn{a}_{L^\infty(T^*\br^n)}
$ 
(see \cite[Proposition~4.2]{ln:coh}).  
We obtain from the definition that
$a^{Wick} = a_0^{w}$ where
\begin{equation}\label{gausreg}
a_0(w) = {\pi}^{-n}\int a(z)\exp(-|w-z|^2)\,dz 
\end{equation}
is the Gaussian regularization. Observe that real Wick symbols have
real Weyl symbols.

In the following, we shall
assume that $G_1 =H_1g^\sh$ is a slowly varying metric satisfying
~\eqref{tempest} and $M$ is a weight for $G_1$ 
satisfying~\eqref{mtemp}.
Also recall that $S^+(1,g^\sh)$ is given by Definition~\ref{s+def}.

\begin{prop}\label{propwick} 
Assume that $a \in L^\infty(T^*\br^n)$, then $a_0^w = a^{Wick}$ where
$a_0$ is given by ~\eqref{gausreg}. If $|a| \le CM$ then we find that
$a_0 \in S(M,g^\sh)$. If also $a \in S(M,G_1)$ in a $G_1$ ball of
fixed radius with center $w$, then $a_0 \cong a$ modulo symbols
in $S(H_1M,G_1)$ in a fixed $G_1$ neighborhood of $w$.  If $a \ge M$
we obtain $a_0 \ge cM$, and if $a \ge M$ in a $G_1$ ball of fixed radius
with center $w$ then $a_0 \ge cM -CH_1M$ in a fixed $G_1$ neighborhood
of $w$, for some constants $c$, $C >0$.
If $|da| \le C$ almost everywhere, then  $a_0\in
S^+(1, g^\sh)$. 
\end{prop} 

These results are well known, but for completeness we give a proof.
Observe that the results are uniform in the metrics and weights.

\begin{proof} 
Since $a$ is measurable satisfying $|a| \le C M$, we find that
$a^{Wick} = a_0^{w}$ where $a_0$ is given by ~\eqref{gausreg}.  Since
$M(z) \le C M(w)(1 + |z-w|)^{3}$ by ~\eqref{mtemp}, we obtain that $a_0(w)
= \cal O(M(w))$. By differentiating on the exponential factor, we find
$a_0 \in S(M,g^\sh)$, and similarly we find that $a_0 \ge M/C$ if $a
\ge M$.

If $a \in S(M,H_1\,g^\sh)$ in a $G_1$
ball of radius $c>0$ and center at $w$, then we write
\begin{multline*}
a_0(w) = {\pi}^{-n}\int_{T^*\br^n} a(z)\exp(-|w-z|^2)\,dz =
{\pi}^{-n}\int_{|w-z| \le cH_1^{-1/2}(w)/2} a(z)\exp(-|w-z|^2)\,dz\\ +
{\pi}^{-n}\int_{|w-z| \ge cH_1^{-1/2}(w)/2} a(z)\exp(-|w-z|^2)\,dz
\end{multline*}
where the last term is $\cal O(H_1^N(w)M(w))$ for any $N$.

Thus, after multiplying with a cut-off function, we may assume that $a
\in S(M,G_1)$ everywhere. Taylor's formula gives
\begin{multline*}
a_0(w) = {\pi}^{-n}\int_{T^*\br^n} a(w +z)\exp(-| z|^2)\,dz \\ = a(w) +
 {\pi}^{-n} \int_0^1\int_{T^*\br^n} (1-{\theta})\w{a''(w +
{\theta}z)z,z} e^{-|z|^2}\,dzd{\theta}
\end{multline*}
where $a'' \in S(MH_1,G_1)$ since $G_1 = H_1g^\sh$.  Since
differentiation commutes with convolution, $M(w + {\theta}z) \le C
M(w)(1 + |z|)^{3}$ and $H_1(w + {\theta}z) \le CH_1(w)(1 + |z|)^2$
when $|{\theta}| \le 1$, we find that $a_0(w) \cong a(w)$ modulo
symbols in $S(H_1M,G_1)$.  Similarly, we obtain that $a_0 \ge cM$
modulo $S(H_1M,g^\sh)$ for some $c>0$ if $a \ge M$ in a fixed $G_1$
ball.  Since
$da_0(w) = {\pi}^{-n}\int_{T^*\br^n} da(z)\exp(-|w-z|^2)\,dz $, we obtain
the last  statement. 
\end{proof}

\begin{rem}\label{distrwick}
Observe that if $a(t,w)$ and $g(t,w) \in L^\infty(\br\times T^*\br^n)$
and $\ddt a(t,w) \ge g(t,w)$ in $\cal D'(\br)$ for almost all $w
\in T^*\br^n$, then we find $\w{\ddt (a^{Wick})u,u} \ge \w{g^{Wick}u,u}$
in $\cal D'(\br)$ when $u \in \cal S(\br^n)$.
\end{rem}

 In fact, the condition
means that 
\begin{equation*}
-\int a(t,w) {\phi}'(t)\,dt \ge \int g(t,w){\phi}(t)\,dt \qquad 0 \le
{\phi} \in C_0^\infty(\br)
\end{equation*}
for almost all $w\in T^*\br^n$, which by~\eqref{poswick} gives
\begin{equation*}
-\int \w{a^{Wick}(t,x,D_x)u,u} {\phi}'(t)\,dt \ge \int
\w{g^{Wick}(t,x,D_x)u,u}{\phi}(t)\,dt \qquad 0 \le {\phi} \in
  C_0^\infty(\br)
\end{equation*}
for $u \in  \cal S(\br^n)$.

We are going to use the symbol classes $S(m_{\varrho}^k, g_{\varrho})$
where  $g_{\varrho} = {\varrho}^2g^\sh$ and $m_{{\varrho}}$ is given by
Definition~\ref{h0def}. Observe that $S(m_{\varrho}^k, g_{\varrho}) =
S(m_1^k, g^\sharp)$ for all $0 < {\varrho} \le 1$. 
In fact, $g_{\varrho} \cong g^\sh$ and $m_{\varrho} \le m_1
\le {\varrho}^{-2}m_{\varrho}$ by ~\eqref{h0eq}.
By ~\cite[Corollary~6.7]{bc:sob} we can define Sobolev spaces
$H(m_{\varrho}^k,g_{\varrho})$ with the following properties: $\cal S$
is dense in $H(m_{\varrho}^k,g_{\varrho})$,
the dual of $H(m_{\varrho}^k,g_{\varrho})$ is naturally identified
with $H(m_{\varrho}^{-k},g_{\varrho})$, and 
\begin{equation}\label{7.7}
u \in  H(m_{\varrho}^k,g_{\varrho}) \iff
a^wu \in L^2 = H(1,g_{\varrho}) \qquad \forall\, a \in
S(m_{\varrho}^k,g_{\varrho})
\end{equation}
and then $u = a_0^w v$ for some $a_0 \in
  S(m_{\varrho}^{-k},g_{\varrho})$ and $v \in L^2$.
 Observe that $
H(m_{\varrho}^k,g_{\varrho}) = H(m_1^k,g^\sharp)$ for all $0 <
{\varrho} \le 1$, but not uniformly. 
We also find from ~~\cite[Corollary ~4.4]{bc:sob} that $a^w$ is
bounded as an operator: 
\begin{equation}\label{nest0}
 u \in H(m_{\varrho}^j,g_{\varrho}) \mapsto a^wu \in
H(m_{\varrho}^{j-k},g_{\varrho})\qquad\text{when $a \in
S(m_{\varrho}^k,g_{\varrho})$,}
\end{equation} and the bound only depends on the
seminorms of ~$a$ in $S(m_{\varrho}^k,g_{\varrho})$. 
Let
${\mu}_{\varrho}^w = m_{\varrho}^{Wick}$, i.e.,
\begin{equation}\label{mudef}
{\mu}_{\varrho}(t,w) = {\pi}^{-n}\int_{T^*\br^n}
 m_{\varrho}(t,z)\exp(-|w-z|^2)\,dz
\end{equation}
Since $m_{\varrho}$ satisfies~\eqref{wtemp} we find from
Proposition~\ref{propwick} that $ m_{\varrho}/c_0 \le
{\mu}_{\varrho} \in L^\infty(\br, S(m_{\varrho},g_{\varrho}))$
uniformly for $0 < 
{\varrho} \le 1$ for some $c_0 >0$. For small enough ${\varrho}>0$ we get
invertible operators according to the following result.

\begin{prop}\label{mprop}
There exists $0< {\varrho}_0 \le 1$ such that when $0 < {\varrho} \le
{\varrho}_0$ we have
\begin{align}
&({\mu}_{\varrho}^{1/2})^w ({\mu}_{\varrho}^{-1/2})^w =
{\gamma}_{\varrho}^w\quad\text{is uniformly invertible in $L^2$}\label{mest1}\\
&\frac{1}{2} \le ({\mu}_{\varrho}^{-1/2})^w
{\mu}_{\varrho}^w({\mu}_{\varrho}^{-1/2})^w \le 2\qquad\text{in $L^2$.}\label{mest2}
 \end{align}
The value of ${\varrho}_0$ only depends on the
seminorms of $f$ in $S(h^{-1},hg^\sh)$.
\end{prop}

\begin{proof}
Since $g_{\varrho} = {\varrho}^2g^\sh$ is uniformly ${\sigma}$
~temperate, $g_{\varrho}/g_{\varrho}^{\sigma} = {\varrho}^4$,  and 
$m_{\varrho}$ is   uniformly ${\sigma}$, $g_{\varrho}$ ~temperate, the
calculus gives that 
$
({\mu}_{\varrho}^{1/2})^w ({\mu}_{\varrho}^{-1/2})^w = 1 + r^w_{\varrho} 
$
where $r_{\varrho} \in S({\varrho}^2,g^\sh)$ uniformly for $0 <
{\varrho} \le 1$. We obtain that the $L^2$ operator norm
$\mn{r^w_{\varrho}} \le C_0{\varrho}^2$ uniformly in ${\varrho}$. Thus,
by taking ${\varrho}_0 \le (2C_0)^{-1/2}$ we obtain~\eqref{mest1}.
Similarly, we find that 
$ 
({\mu}_{\varrho}^{-1/2})^w {\mu}_{\varrho}^w({\mu}_{\varrho}^{-1/2})^w = 1 +
s^w_{\varrho} 
$ 
where $s_{\varrho} \in S({\varrho}^2,g^\sh)$ uniformly. As before, we
may choose ${\varrho}_0$ so that $\mn{s^w_{\varrho}} \le 1/2$ when $0 <
{\varrho} \le {\varrho}_0$, which proves~\eqref{mest2}.
\end{proof}

Now we fix ${\varrho}= {\varrho}_0$ to the value given by
Proposition~\ref{mprop}, let ${\mu} = {\mu}_{{\varrho}_0}$ and
let $\mn{u}_{H(m_1^{k})}$ be the norm defining $H(m_1^k) =
H(m_1^k,g^\sh)$. The next Proposition shows that the norm in
$H(m_1^k)$ can be defined by the operator ${\mu}^w$.

\begin{prop}\label{wprop} 
Assume that ${\varrho}_0$ is given by Proposition~\ref{mprop}, and the
symbol ${\mu} = {\mu}_{\varrho_0} \in L^\infty(\br, S(m_1,g^\sh))$ is
given by ~\eqref{mudef}  with
${\varrho}= {\varrho}_0$ so that ${\mu}^w = m_{{\varrho}_0}^{Wick} \le
m_1^{Wick}$.  Then
there exist positive constants $c_0$, $c_1$ and $C_0$ such that
\begin{equation}\label{west3}
c_0 h^{1/2}\mn u^2 \le c_1 \mn{u}_{H\left(m_1^{1/2}\right)}^2 \le
\w{{\mu}^wu,u}\le 
C_0 \mn{u}_{H\left(m_1^{1/2}\right)}^2.
\end{equation}
The constants only depend on the
seminorms of $f$ in $S(h^{-1},hg^\sh)$.
\end{prop}

\begin{proof} Let $a = {\mu}^{-1/2}$.
By Proposition~\ref{mprop} we find $1 =
({\gamma}_{\varrho_0}^w)^{-1}({\mu}^{1/2})^w a^w$ with
$\mn{({\gamma}_{\varrho_0}^w)^{-1}} \le C$, which gives 
\begin{equation}\label{west1} 
 \mn u_{H(1)} \le  C \mn {({\mu}^{1/2})^w a^w u}_{H(1)} \le C' \mn
 {a^w u}_{H\left(m_1^{1/2}\right)}  \le C'' \mn u_{H(1)}
\end{equation}
by \eqref{nest0}. Thus $u \mapsto a^wu$ is an isomorphism between $H(1)$ and
$H(m_1^{1/2})$. Since the constant
metric $g^\sh$ is trivially strongly ${\sigma}$ ~temperate in the
sense of \cite[Definition~7.1]{bc:sob}, we find from
\cite[Corollary~7.7]{bc:sob} that there exists $a_0 \in
S(m_1^{1/2},g^{\sh})$ such that $a^wa_0^w = a_0^wa^w = 1$.

Since $c_0h^{1/2} \le m_1 \le {\varrho}_0^{-2}m_{{\varrho}_0} \le 1$
and ${\mu}^w = m_{\varrho_0}^{Wick}$ we find that $c_0{\varrho}_0^2h^{1/2}\mn u^2
\le \w{{\mu}^wu,u} \le C\mn{u}^2$, so we only have to prove that
$\w{{\mu}^wu,u} \cong \mn{u}_{H\left(m_1^{1/2}\right)}^2$. Since $a =
{\mu}^{-1/2}$ we find from \eqref{mest2} that
\begin{equation*}
\w{{\mu}^wu,u} = \w{{\mu}^wa^wa_0^wu, a^wa_0^wu} \ge \frac{1}{2}\mn{a_0^w
  u}^2_{H(1)}, 
\end{equation*}
and since $\mn{u}_{H\left(m_1^{1/2}\right)} = \mn{a^wa_0^w
  u}_{H\left(m_1^{1/2}\right)}  \le
C \mn{a_0^w u}_{H(1)}$, we find $\w{{\mu}^wu,u} \ge \frac{1}{2C}
\mn{u}_{H\left(m_1^{1/2}\right)}^2$. Finally, we have
\begin{equation*}
\w{{\mu}^wu,u} \le \mn{{\mu}^w
  u}_{H\left(m_1^{-1/2}\right)}\mn{u}_{H\left(m_1^{1/2}\right)} \le C 
\mn{u}_{H\left(m_1^{1/2}\right)}^2
\end{equation*}
which completes the proof of the Proposition.
\end{proof}

\section{The Pseudo-Sign}\label{aps}

In this section we shall  construct a perturbation $B(t,w)=
{\delta}_0(t,w) + {\varrho}_0(t,w)$ of ${\delta}_0$ such that ${\varrho}_0 =
\cal O(m_{1})$ and
\begin{equation} \label{r0prop}
\ddt ({\delta}_0 + {\varrho}_0) \ge c m_{1} >0
\qquad \text{in $\cal D'(\br)$ when $|t| < 1$,} 
\end{equation}
where $m_{1}$ is given by Definition~\ref{h0def} and ${\delta}_0$ by
Definition~\ref{d0deforig}.  We shall use this in Section~\ref{lower}
to prove Proposition~\ref{apsprop} with $b^w = B^{Wick}$ as
a ``pseudo-sign'' for $f$.  When $t\mapsto m_{1}(t,w)$ has a approximate
minimum at $t = t_0$ in the sense that $m(s) \le Cm(t)$
when $t \le s \le t_0$ or  $t_0 \le s \le t$, we may take
${\varrho}_0(t,w) = c\int_{t_0}^t 
m_{1}(s,w)\,ds$ since $t \mapsto {\delta}_0(t,w)$ is
non-decreasing. In general, we have to split the interval $[-1,1]$
into subintervals where $t\mapsto m_{1}(t,w)$ has approximate maximum
and minimum, and use the ``convexity property'' of $t \mapsto
{\delta}_0(t,w)$ given by Proposition~\ref{qmaxpropo} in order to
``interpolate'' ${\delta}_0$ at the approximate maxima of $t \mapsto
{\delta}_0(t,w)$.
We shall also compute the Weyl symbol $b$ for the ``pseudo-sign''
$B^{Wick} = b^w$.  All the results in this section are uniform in the
sense that they only depend on the seminorms of $f$ in
$S(h^{-1},hg^\sh)$ for $|t| \le 1$. As before, we denote by
$\lip(T^*\br^n)$ the Lipschitz continuous functions on ~~$T^*\br ^n$.

\begin{prop}\label{apsdef}
Assume that ${\delta}_0$ is given by
Definition~\ref{d0deforig} and  $m_{1}$ is given by
Definition~\ref{h0def}. Then there exist 
a positive constant $C_1$ and a real valued ${\varrho}_0(t,w)
\in L^\infty(\br\times T^*\br^n)$
such that
\begin{align}\label{r0prop0}
&|{\varrho}_0| \le C_1 m_{1}\\
&\partial_t({\delta}_0 + {\varrho}_0) \ge m_{1}/C_1 \label{r0prop1}
\end{align}
in $\cal D'(\br)$ when $|t| <1$.  We also have that $t \mapsto
{\varrho}_0(t,w)$ is a regulated function, $\forall \, w\in T^*\br^n$,
and $w \mapsto {\varrho}_0(t,w)\in \lip(T^*\br^n)$ uniformly for
almost all $|t| \le 1$.
\end{prop}

\begin{proof}
We shall make the construction of ${\varrho}_0(t,w)$ locally in $w$,
by using a partition of unity $\set{{\phi}_k(w)} \in S(1,g^\sh)$ in
$T^*\br^n$ such that $0 \le {\phi}_k \le 1$, $\sum_k {\phi}_k = 1$ and
$|w-w_k| \le {\varepsilon}$ in $\supp {\phi}_k$. Here $0
<{\varepsilon} \le \min({\varepsilon}_0,{\varepsilon}_1)$, where
${\varepsilon}_j$ are given by
Propositions~\ref{qmaxpropo}--\ref{dvar1} for $j=1$, $2$. Observe that
${\phi}_k \in \lip(T^*\br^n)$ uniformly. We also assume
that ${\varepsilon}$ is chosen small enough so that $w \mapsto
m_{1}(t,w)$ only varies with a fixed factor in $\supp {\phi}_k(w)$ for
any $t \in [-1,1]$, $\forall\, k$, and we shall keep ${\varepsilon}$
fixed in what follows.

On each $\supp {\phi}_k$ we shall construct a real valued
${\varrho}_{0k}(t,w)\in L^\infty(\br\times T^*\br^n)$ satisfying the
conditions in Proposition~\ref{apsdef} uniformly in
$\supp {\phi}_k$. By taking ${\varrho}_0(t,w) = \sum_k {\phi}_k(w)
{\varrho}_{0k}(t,w)$ we then obtain the result. 
Observe that we may ignore the values of
${\varrho}_{0k}(t,w)$ for $t$ in a zero set.

In the following, we shall keep $k$ fixed. Next, we  choose
coordinates so that $w_k = 0$, let 
${\varrho}_{0}(t,w) = {\varrho}_{0k}(t,w)$ and
\begin{equation} \label{mdef}
m(t) = m_{1}(t,0)
\end{equation}
which gives $m(t) \cong m_1(t,w)$ when $w \in \supp {\phi}_k$ and $t
\in [-1,1]$.
Thus it suffices to construct a real valued ${\varrho}_0(t,w)$ 
such that $|{\varrho}_0(t,w)| \le C_1m(t)$ and $\partial_t({\delta}_0+
{\varrho}_0) \ge m/C_1$ in $\cal D'(\br)$
when $|w| \le {\varepsilon}$ and $|t| < 1$. This is essentially a one
dimensional problem, but there are some complications at the
approximate maxima of $t \mapsto m(t)$.

Since $t \mapsto m(t)$ and $t \mapsto {\delta}_0(t,w)$
are regulated functions, we
may consider them as functions of $\wt t \in S^*\br\bigcup \br$.
Thus $m(\,\wt
t\,)$ is either $m(t)$ or $m(t\pm) = \lim_{{\varepsilon}\searrow 0}
m(t \pm {\varepsilon})$ depending on the context.
We introduce an
ordering on $S^*\br \bigcup \br$ so that
\begin{equation}\label{order}
s < s+ < t- < t\qquad\text{when $s < t$.} 
\end{equation}
In the next Lemma we shall cut the interval $[-1,1]$ where the jumps
of $t \mapsto {\delta}_0(t,w)$ are large enough for $w \in  \supp
{\phi}_k$. Then we also get a 
bound on the jumps of $t \mapsto m(t)$ in the subintervals by
Proposition~\ref{qmaxpropo}.

\begin{lem}\label{apslem}
Let $m(t)$ be given by ~\eqref{mdef}, then
there exist finitely many open disjoint
subintervals~ $I_k \subseteq\ [-1,1]$ so that $[-1,1] = \bigcup_k \ol
I_k$ and
\begin{equation}\label{jumpcond1}
{\delta}_0(t+,0)- {\delta}_0(t-,0)  \le
9 \max(m(t-),m(t+))\qquad t \in \bigcup_k I_k
\end{equation}
and 
\begin{equation}\label{jumpcond2}
 |m(t)|\le {\kappa}_1 \max(m(t-) , m(t+))\qquad t \in \bigcup_k I_k.
\end{equation}
Here ${\kappa}_1 = \max({\kappa}_0,
9/c_0) > 1$ with $c_0$ and ${\kappa}_0$ given
by Proposition~\ref{qmaxpropo}.  We also obtain that
\begin{equation}\label{jumpcond3}
{\delta}_0(t+,w)- {\delta}_0(t-,w)  >
3 \max(m(t-),m(t+))\qquad \text{if\/ $|w| \le {\varepsilon}$}
\end{equation}
when  $\pm 1 \ne t \in \bigcup_k \partial I_k$. 
\end{lem}

\begin{proof}[Proof of Lemma~\ref{apslem}] 
Since $t \mapsto {\delta}_0(t,w)$ is non-decreasing and
$|{\delta}_0| \le h^{-1/2}$, we find that for any ${\gamma} >0$
there can only be finitely may values of ~$t$ such that
${\delta}_0(t+,0)- {\delta}_0(t-,0) \ge {\gamma}$ when $|t| <
1$.  Since $m \ge ch^{1/2}$ for some $c>0$ by ~\eqref{hhhest}, there
are only finitely many values of ~$t$ for which ${\delta}_0(t+,0)-
{\delta}_0(t-,0) > 9 \max(m(t-),m(t+)) \ge 9ch^{1/2}$. Thus,
by cutting the interval $]{-1},1[$ into finitely many parts at the
discontinuity points of $t \mapsto {\delta}_0(t,0)$ we may assume
that $[-1,1] = \bigcup_k \ol I_k$, that~\eqref{jumpcond1} holds 
and
\begin{equation}\label{jumpcond0b}
{\delta}_0(t+,0)- {\delta}_0(t-,0) >
9 \max(m(t-),m(t+))\qquad \pm 1 \ne t \in \bigcup_k\partial I_k.
\end{equation} 

By letting $t'$, $t'' \to t_0 = t\in \bigcup_k I_k$ in
Proposition~\ref{qmaxpropo}, we 
obtain that $m(t) \linebreak[0] \ge {\kappa}\max(m(t-),m(t+))$ implies
$$ {\delta}_0(t+,0)- {\delta}_0(t-,0) \ge c_0{\kappa} \max(m(t-),m(t+))
$$ 
when ${\kappa} \ge
{\kappa}_0$. This contradicts  ~\eqref{jumpcond1} when
${\kappa} > \max({\kappa}_0, 9/c_0)$, and gives~\eqref{jumpcond2}. 
By letting $s_0$, $t_0 \to t\in \bigcup_k \partial I_k\setminus
\set{\pm 1}$, $s_0 < t < t_0$, and using Proposition~\ref{dvar1} with
${\lambda} >  9$ and $K = \max(m(t+), m(t-)) \le 1$ we obtain  
~\eqref{jumpcond3} since ${\varepsilon} \le {\varepsilon}_1$. This
completes the proof of the Lemma. 
\end{proof}

Thus, it suffices to construct real valued ${\varrho}_0(t,w)$ satisfying
\begin{equation}\label{locr0prop0} 
\left\{
\begin{aligned}
&|{\varrho}_0(t,w)| \le C_1 m(t)\\
&\partial_t({\delta}_0(t,w) + {\varrho}_0(t,w)) \ge m(t)/C_1 
\end{aligned}\right.
\qquad\text{in $\cal D'(\br)$ when $t
\in \bigcup I_k$ and $|w| \le {\varepsilon}$.}
\end{equation}
We shall also obtain that $t \mapsto
{\varrho}_0(t,w)$ is regulated, $w \mapsto {\varrho}_0(t,w)$ is
Lipschitz and
\begin{equation}\label{mjumpcond}
 |{\varrho}_0(t\pm,w)| \le m(t\pm) \qquad \text{for $\pm 1 \ne t \in
 \bigcup_k \partial I_k$ and $|w| \le {\varepsilon}$.}
\end{equation} 
Then we obtain from~\eqref{jumpcond3} that $t \mapsto B(t,w)= {\delta}_0(t,w)
+{\varrho}_0(t,w)$ has positive jumps $B(t+,w) - B(t-,w) > 0$
when $t \in \bigcup_k \partial I_k \setminus \set{\pm 1}$ and $|w|
\le {\varepsilon}$.
In the following, we shall consider one of the
intervals $I_k$ and let $I =I_k$ be fixed. We shall
assume that ~\eqref{jumpcond1} and~\eqref{jumpcond2} hold for $t \in I$,
and we shall split the interval $I$ into subintervals where $t\mapsto
m(t)$ has an approximate minimum. The complication is that $t \mapsto
m(t)$ is regulated, and not necessarily continuous. Because of the
``convexity property'' given by Proposition~\ref{qmaxpropo}, we cannot in
general change the value of $m(t)$ even at finitely many points.

\begin{lem}\label{cutlemma1} 
Assume that $m(t)$ is given by ~\eqref{mdef}, that ~\eqref{jumpcond1}
and~\eqref{jumpcond2} are satisfied on an open interval $I \subseteq\
[{-1},1]$, and let ${\kappa} > {\kappa}_1$ where ${\kappa}_1 >1$ is
given by Lemma~\ref{apslem}.
Then there exists an open interval $I_1 = \ \left]r_1, r_2\right[\ 
\subseteq I$, and 
$t_1 \in \ol I_1$ such that 
\begin{align}
&m(\, \wt t_1) = \inf_{t \in
  I_1}m(t)\label{unionprop1}\\ 
&m(s) \le
{\kappa}^2m(t)\qquad \text{for $t_1 < s \le t$ or $t \le s <
t_1$ when $t \in I_1$.}\label{qminprop1}
\end{align}
Let $M_1 = \sup _{r_1 < t < t_1} m(t)$ if $r_1 \ne t_1$, and
$M_2 = \sup _{t_1 < t < r_2} m(t)$ if $r_2 \ne t_1$. If\/
$r_1$ or $r_2 \notin \partial I$ we have $r_j \ne t_1$,  $j=1$ or $2$, and
obtain that
\begin{equation}\label{mMcomp1}
m(r_1-) \le M_1/ {\kappa} \qquad\text{if $r_1
  \notin \partial I\quad$ or}\qquad
m(r_2+) \le M_2/ {\kappa}\qquad\text{if $r_2
  \notin \partial I\quad$}.
\end{equation}
Thus, $t_1 \in \partial I_1$ if and  only if $t_1 \in \partial I$.
If $I_1 \ne I$ then we find
\begin{equation}\label{wtdvar}
\sup _{t \in \wt I} {\delta}_0(t,0) - \inf_{t \in \wt I} 
{\delta}_0(t,0) \ge  c_0 ch^{1/2}
\end{equation}
for any open interval $\wt I$ such that  $I_1 \subsetneq \wt
I\subseteq I$.   
\end{lem}

Since we are going to take ${\varrho}_0(t) = {\gamma}_0\int_{t_1}^t
m(s)\,ds$ near $t_1$, property ~\eqref{qminprop1} will
give $|{\varrho}_0(t)| \le 2{\gamma}_0{\kappa}^2m(t) \le m(t)$ if
${\gamma}_0 \le (2{\kappa}^2)^{-1}$. 

\begin{proof}[Proof of Lemma~\ref{cutlemma1}]
Let 
$ 
m_1 = \inf _{t\in I}m(t),
$ 
since $m(t)$ is regulated, we may choose
$ 
t_1 \in \ol I$ such that $m(\,\wt t\,\,) =m_1.
$ 
If $t_1 \in \partial I$, then we may take $\wt t_1 = t_1\pm$
depending on whether $\pm (t -t_1) >0$ for all $t \in I$. 
Next, for fixed ${\kappa} > {\kappa}_1 >1$ we define $I = ]r_1, r_2[$ with
\begin{align}\label{9.15}
&r_1 = \inf \set{t \in I:\  m(s) \le
  {\kappa}^2m(t)\quad\text{for all $t \le s <
t_1$}} \\\label{9.16}
&r_2 = \sup \set{t \in I:\ m(s) \le
  {\kappa}^2m(t)\quad\text{for all $t_1 < s \le t$}}.
\end{align}
If $t_1 \in \partial I$ then one condition is empty, depending on
whether $(-1)^j(t -t_1) < 0$ for all $t \in I$, and then we put $r_j =
t_1$. Observe that we do not use the value of $m$ at $t_1$ in
\eqref{9.15}--\eqref{9.16}. Since $t \mapsto m(t)$ is regulated and
${\kappa} > 1$, we find that $t_1 = r_j$ if and only if $r_j \in
\partial I$, and thus $r_1 < r_2$. In fact, if there exist $t_1 <
s_{\varepsilon} < t_1 + {\varepsilon}$ such that $m(s_{\varepsilon}) >
{\kappa}^2m(t_1 + {\varepsilon})$ for a sequence ${\varepsilon}
\searrow 0$, then $m(t_1+) \ge {\kappa}^2m(t+)$ which gives a
contradiction since ${\kappa} > 1$.

Next, we shall prove ~\eqref{mMcomp1}. 
If, for example, $r_2 \notin \partial I$ then since $t \mapsto m(t) \ge
ch^{1/2}$ is regulated and ${\kappa} >1$, we can find
${\varepsilon}_2>0$ so that
\begin{equation*} 
{\kappa}^{-1}m(r_2+) \le m(r_2 +
{\varepsilon}) \le {\kappa} m(r_2+)\qquad\text{when $0 < {\varepsilon}
  \le {\varepsilon}_2$.}
\end{equation*}
This implies that $m(t) \le {\kappa}m(r_2+) \le {\kappa}^2m(r_2+
{\varepsilon})$ for $r_2 < t < r_2 + 
{\varepsilon}$ and $0 < {\varepsilon} < {\varepsilon}_2$. If $M_2 <
{\kappa}m(r_2+ \wt {\varepsilon})$ for some $0 < \wt{\varepsilon} \le
{\varepsilon}_2$, then  $m(t) < {\kappa} m(r_2 +
\wt {\varepsilon})$ for $t_1 < t < r_2$, and we find from  
~\eqref{jumpcond2} that 
$$m(r_2) \le {\kappa} 
\max(m(r_2-),m(r_2+)) \le {\kappa}\max(M_2,m(r_2+)) \le
{\kappa}^2m(r_2+ \wt {\varepsilon}).$$ 
Thus, $m(t) \le
{\kappa}^2m(r_2+ \wt{\varepsilon})$  for $t_1 < t \le r_2 + \wt {\varepsilon}$,
which contradicts the definition of ~$r_2$. 

Now, if $r_1 \notin \partial I$ then $r_1 \ne t_1$ and we may choose
$r_1 < \wt t < \wt t_1$ such that $m(\,\wt t\,) = M_1 \ge
{\kappa}m(r_1-) \ge {\kappa}m(\,\wt t_1)$  by~\eqref{mMcomp1}. 
Similarly, if $r_2 \notin \partial I$ then $M_2 \ge
{\kappa}\max(m(r_2+),m(\,\wt t_1))$. 
By taking limits in
Proposition~\ref{qmaxpropo}, we find that
$
{\delta}_0(\,\wt t_1, 0) - {\delta}_0(r_1-,0) \ge
c_0 M_1  \ge C_0ch^{1/2}
$ 
if $r_1
  \notin \partial I$ and 
${\delta}_0(r_2+,0) - {\delta}_0(\,\wt t_1,0) \ge
c_0 M_2 \ge C_0ch^{1/2}$ 
if $r_2 \notin \partial I$.
Combining these estimates, we obtain ~\eqref{wtdvar}.
\end{proof}

We can now cut the interval $I$ into subintervals where $t\mapsto
m(t)$ is approximately monotone. 

\begin{prop}\label{cutlemma} 
Assume that $m(t)$ is given by ~\eqref{mdef}, that ~\eqref{jumpcond1}
and~\eqref{jumpcond2} are satisfied on an open interval $I \subseteq\
[{-1},1]$, and let ${\kappa} > {\kappa}_1$ where ${\kappa}_1 > 1$ is
given by Lemma~\ref{apslem}. Then there exists a finite family
$\set{I_k}_{k=1}^{N}$ of open and disjoint subintervals $I_k \subseteq
I$ and points $t_k \in \ol I_k$ such that $t_j \ne t_k$ if $j \ne k$
and
\begin{align}
&m(\, \wt t_k) = \inf_{t \in
  I_k}m(t)\qquad\forall\, k \label{unionprop}\\
&\ol I = \bigcup_k \ol I_k\label{8.18}
\\ 
&m(s) \le
{\kappa}^2m(t)\qquad \text{for $t_k < s \le t$ or $t \le s <
t_k$ when $t \in I_k$.}\label{qminprop}
\end{align}
If $N >1$ and $\ol I_j  \bigcap \ol I_k \ne \emptyset$ we also have
\begin{equation}
\sup_{t_j < t < t_k} m(t) > {\kappa}\max(m(\, \wt t_j), m(\, \wt t_k))  \qquad
\forall\,j,\ k.\label{qmaxprop}
\end{equation}
\end{prop}

As before, property ~\eqref{qminprop} will
give $|{\varrho}_0(t)| \le {\gamma}_0\int_{t_k}^t
m(s)\,ds \le m(t)$ on ~$I_k$. When joining the constructions we
shall find that Proposition~\ref{qmaxpropo} gives a sufficiently large
increase of $t \mapsto {\delta}_0(t,w)$ by~\eqref{qmaxprop} so that we can
interpolate between these parts.  Note that ~\eqref{qminprop} is empty
in one of the cases if $t_k \in \partial I_k$.

\begin{proof}[Proof of Proposition~\ref{cutlemma}]
We shall obtain the Proposition by repeatedly ``cutting'' the interval
$I$ using Lemma~\ref{cutlemma1}.  The Lemma first gives $I_1 =\left
]\,r_1, r_2\,\right [ \subseteq I$ satisfying
\eqref{unionprop1}--\eqref{wtdvar} which gives~~\eqref{unionprop}
and~\eqref{qminprop} for $k=1$. We also have
that $t_1 \in I_1$ or $t_1 \in \partial  I \bigcap \partial  I_1$.
If $I_1 = I$, i.e., $r_1$ and $r_2\in \partial I$ then we obtain
~\eqref{8.18}, $N=1$ and we are finished.

Else, we have ${\Xi}_2 = I \setminus \ol I_1 \ne \emptyset$.  Recursively
assume that we have chosen $I_j$ satisfying ~~\eqref{unionprop}
and~\eqref{qminprop} for $j < k$ so that ${\Xi}_k = I\setminus \left(
\bigcup_{j<k} \ol I_j \right) \ne \emptyset$.  Then we take $m_k =
\inf_{t \in {\Xi}_k} m(t) = m(\,\wt t_k)$ and use
Lemma~\ref{cutlemma1} to ``cut'' the component of ~$\ol {\Xi}_k$
containing ~$t_k$.  Then we obtain a new open interval $I_k \subseteq
{\Xi}_k$ satisfying \eqref{unionprop1}--\eqref{wtdvar} with ~$t_1$
replaced by~ $t_k$, $I_1$ by ~$I_k$ and $I$ by ~${\Xi}_k$.  Thus,
unless $t_k \in I_k$ we have that $t_k \in \partial I_k \bigcap
\partial\, {\Xi}_k$ and since $\partial\,{\Xi}_k \subseteq
\bigcup_{j<k}\partial I_j$ we then find $t_k \in \partial I_k \bigcap
\partial I_j$ for some $j < k$. Since three disjoint intervals in
~~$\br$ cannot have intersecting boundary, we find that
$\bigcap_{i=1}^3 \partial I_{j(i)} = \emptyset$ which gives that $t_k
\ne t_j$ for any $j < k$.

Note that, unless $I_k$ is equal to one of the components of
~${\Xi}_k$, we obtain from ~\eqref{wtdvar} that $t \mapsto
{\delta}_0(t,0)$ increases more than $c_0ch^{1/2}$ in any open
interval $\wt I$ such that $I_k \subsetneq \wt I \subseteq
{\Xi}_k$. If we take the intersection of these intervals we find that
$\bigcap_{I_k \subsetneq \wt I \subseteq {\Xi}_k} \wt I = \ol I_k
\bigcap {\Xi}_k$. Since $t \mapsto {\delta}_0(t,0)$ is monotone and
bounded, there is a fixed bound on the number of such intervals.
Thus, we may repeat the process only finitely many times until $\ol I
= \bigcup_{0 <k \le N} \ol I_k$. (The proof would in fact also work
with an infinite number of subintervals.) In fact, in order to get
infinitely many intervals, we must infinitely many times ``cut'' a
remaining component of the sets ~${\Xi}_k$.
In the following we
keep the original enumeration of the intervals, so that if $j < k$ then 
$I_j$ was ``cut'' before $I_k$. 

It remains to prove ~\eqref{qmaxprop}, and
it is no restriction to assume that $t_j <
t_k $. First assume that $j < k$ ($I_j$ was defined before
$I_k$) and let $\ol I_j \bigcap \ol I_k =
\set{r_{jk}}$, then $r_{jk}  \in  {\Xi}_j$. 
In fact, if $r_{jk} \in \partial \,{\Xi}_j$ then $r_{jk} \in \partial
I_i$ for some $i < j$ which is impossible.
Thus, by using ~\eqref{mMcomp1}
we find that 
\begin{equation*}
m_j \le m(r_{jk}+) \le
\sup_{t_j < t < r_{jk}}m(t)/{\kappa}
 \le \sup_{t_j < t < t_k}m(t)/{\kappa}
\end{equation*}
since $m_j$ is the infimum of $m(t)$
over ${\Xi}_j$.
Now $r_{jk}+ {\varepsilon} \in I_k$ for small enough
${\varepsilon}$, thus we find that $m(r_{jk}+) \ge m_k$.
A similar argument gives~\eqref{qmaxprop} when $k < j$ 
and completes the proof of Proposition~\ref{cutlemma}.
\end{proof}

Now we proceed with proof of Proposition ~\ref{apsdef}, i.e., the
construction of ${\varrho}_{0}(t,w)$ when $|w| \le {\varepsilon}$
and $t \in I \subseteq [-1,1]$.  Recall that $m(t) = m_1(t,0)$ is
given by~\eqref{mdef}, thus $m(t) \cong m_1(t,w)$ when $|w| \le
{\varepsilon}$. We have also assumed that ~\eqref{jumpcond1}
and~\eqref{jumpcond2} hold in $I$, and we shall construct
${\varrho}_0(t,w)$ satisfying~\eqref{locr0prop0} in~$I$ 
and~\eqref{mjumpcond} on~$\partial I$. We also need to prove that $t \mapsto
{\varrho}_0(t,w)$ is regulated in $I$, $\forall\, w \in
T^*\br^n$, and $w \mapsto {\varrho}_0(t,w)$ is uniformly
Lipschitz when $|w| \le {\varepsilon}$, $\forall\,t \in I$.  By
Proposition~\ref{cutlemma} we obtain $I_k \subseteq I$, $k= 1, \dots,
N$, such that the properties ~\eqref{unionprop}--\eqref{qmaxprop} hold
with ${\kappa} > {\kappa}_1 > 1$ and $t_j \in \ol I_j$ such that $t_j
\ne t_k$ when $j \ne k$. In the following, we shall assume that $t_j$
are ordered so that $t_j < t_k$ if and only if $j < k$. First we let
${\gamma}_0 >0$ and define
\begin{equation}\label{rho0def}
{\varrho}_0(t) =
\left\{
\begin{aligned}
&{\gamma}_0\int_{t_1}^t m(s)\,ds \qquad t \le t_1
\\
&{\gamma}_0\int_{t_N}^t m(s)\,ds \qquad t \ge t_N
\end{aligned}
\right. 
 \qquad \text{for $|w| \le {\varepsilon}$,}
\end{equation}
which is constant in $w$. Then we obtain that $t \mapsto
{\varrho}_0(t)$ is continuous when $t \le t_1$ or $t \ge t_N$, and
by~\eqref{qminprop} we have that
$$|{\varrho}_0(t)| \le 2{\gamma}_0{\kappa}^2m(t) \le m(t)\qquad\text{when $t
\le  t_1$ or $t \ge t_N$}$$ 
if ${\gamma}_0 \le (2{\kappa}^{2})^{-1}$, which
gives~\eqref{mjumpcond}.  We also
find that 
$$\ddt {\varrho}_0(t) = {\gamma}_0m(t)\qquad\text{for almost all $t
\le t_1$ or $t \ge t_N$.}$$ We shall assume that $0 < {\gamma}_0 \le
(2{\kappa}^{2})^{-1} <1/2$ in what follows, but later we shall
impose more conditions on ${\gamma}_0$.

In the following, we shall put
\begin{equation}\label{nymudef} 
{\mu}(t,w) =
{\delta}_0(t,w) + {\varrho}_0(t,w)\qquad |w| \le {\varepsilon}
\quad t \in I
\end{equation} 
and first we define ${\mu}(t_j,w) = {\delta}_0(\,\wt t_j,w)$.
It remains to construct ${\mu}(t,w)$ (or ${\varrho}(t,w)$) on
$]t_j,t_{j+1}[$ when $|w| \le 
{\varepsilon}$, $j=1,\dots, N-1$, and it is no restriction  
to consider $j=1$. 
When constructing ${\mu}(t,w)$ in $]\,t_1,t_2 \,[$ we shall ensure that
\begin{equation}\label{8.31}
 {\delta}_0(\,\wt t_1,w) \le {\mu}(t,w) \le {\delta}_0(\,\wt
 t_2,w)\qquad \text{ when $t_1 < t < t_2$ and $|w| \le
 {\varepsilon}$.}
\end{equation}
Then we obtain that $t \mapsto {\mu}(t,w)$ has non-negative jumps at $t =
t_1$ and $t_2$.  Let $r_{12} = \ol I_1 \bigcap \ol I_2$, observe that it is
possible that $r_{12} = t_1 $ or $t_2$. For $t_1 < t < t_2$, $t
\ne r_{12}$, we find from ~\eqref{qminprop} that
\begin{equation}\label{9.25}
m(s) \le {\kappa}^2m(t) \quad\text{for $t_1 < s \le t < r_{12}$ or $r_{12}
<t \le s <t_2$,}
\end{equation}
which could be empty in one case. We shall now determine where $t
\mapsto m(t)$ has an approximate maximum.
We find from~\eqref{qmaxprop} that
\begin{equation}\label{9.24}
M_{12} = \sup_{t_1 < t < t_2}m(t) > {\kappa}\max(m_1,m_2) 
\end{equation}
since $m_j = m(\,\wt t_j\,)$.
Define
\begin{align}\label{9.26}
&s_1 = \sup\set{t\in I_1:\ m(s) \le M_{12}/{\kappa}
  \quad\text{for all $t_1 <s \le t$}} \\ \label{9.27}
&s_2 = \inf\set{t\in I_2:\ m(s) \le
  M_{12}/{\kappa} \quad \text{for all $t \le s < t_2$}} .
\end{align}
Since the condition in ~\eqref{9.26} is empty when $t \le t_1$ and the
condition in ~\eqref{9.27} is empty when $t \ge t_2$, we find that
$t_1 \le s_1 \le r_{12} \le s_2 \le t_2$. Observe that we could have
$s_j = t_j$ for $j=1$ or $2$, for example if $t_j= r_{12} \in \partial
I_j$.  If $t_1 < r_{12}$ and $m(t_1+) < M_{12}/{\kappa}$ then we find 
that $s_1 > t_1$; similarly if $r_{12} < t_2$ and $m(t_2-) <
M_{12}/{\kappa}$ then we find 
$s_2 < t_2$.  We shall define ${\varrho}_0$ by~\eqref{rho0def} when
$t_1 < t < s_1$ and $s_2 < t < t_2$; when $s_1 < t < s_2$ we shall use
that $m \cong M_{12}$ has an approximate maximum according to the following

\begin{lem}\label{intlem}
Let $m(t)$ be given by ~\eqref{mdef} satisfying  ~\eqref{jumpcond2}
and~\eqref{qminprop}, let $M_{12}$ be given by ~\eqref {9.24}
and let $s_1$, $s_2$ be given by ~\eqref{9.26}--\eqref{9.27}. Then
we find that $s_1 < s_2$ and 
\begin{equation}\label{mcond2}
M_{12}/{\kappa}^3 < m(t) \le M_{12}\qquad \text{when $s_1 < t < s_2$ and $t
  \ne r_{12}  \in \ol I_1 \bigcap \ol I_2$.}
\end{equation}
\end{lem}

\begin{proof} [Proof of Lemma~\ref{intlem}]
We shall show that if $s_1 = r_{12} = s_2$ then $m(r_{12}\pm) \le
M_{12}/{\kappa} = m(r_{12})/{\kappa}$, which contradicts
~\eqref{jumpcond2} since 
${\kappa} > {\kappa}_1$.  In the case that $r_{12} = s_j \ne t_j$, $j= 1$
and $2$, we immediately obtain this from from the definition of $s_j$.
When $t_1 = r_{12} = s_2$ then $t_1 \in \partial I_1$, so that
$m(r_{12}-) = m(\,\wt t_1) =m_1 < M_{12}/{\kappa}$.  In that case $s_2
\ne t_2$ so $m(r_{12}+) \le
M_{12}/{\kappa}$.  A similar argument works when $t_2 = r_{12} = s_1$.

It remains to prove ~\eqref{mcond2}. By the definition~\eqref{9.26} we
find that if~$s_1 \ne r_{12}$ then for any $s_1 < t < r_{12}$ there
exists $t_1 < s \le t$ such that $M_{12}/{\kappa} < m(s)$ and by
~\eqref{qminprop} we have $m(s) \le {\kappa}^2m(t) \le
{\kappa}^2M_{12}$, which gives ~\eqref{mcond2} for these intervals.
We similarly obtain \eqref{mcond2} when $r_{12} < t < s_2$, which
completes the proof of the Lemma.
\end{proof}

If $s_1 > t_1$ we define 
\begin{equation}\label{rhoint1}
{\varrho}_0(t) = {\gamma}_0\int_{t_1}^t m(s)\,ds \qquad t_1 < t < s_1
  \quad \text{and} \quad \ |w| \le {\varepsilon}
\end{equation}
and if $s_2 < t_2$  we define 
\begin{equation}\label{rhoint2}
{\varrho}_0(t) = {\gamma}_0\int_{t_2}^t m(s)\,ds \qquad s_2 < t < t_2
\quad  \text{and} \quad   \ |w| \le {\varepsilon}
\end{equation}
which is constant in $w$.
Then we obtain that $t \mapsto {\varrho}_0(t)$ is continuous in
$]\,t_1,s_1\,[\, \bigcup\, ]\,s_2,t_2\,[$ and 
\begin{equation} \label{rho0size}
|{\varrho}_0(t)| \le
2{\gamma}_0{\kappa}^2m(t) \le m(t)\qquad \forall\, t \in
]\,t_1,s_1\,[\, \bigcup\, ]\,s_2,t_2\,[ 
\end{equation} 
by ~\eqref{qminprop} since ${\gamma}_0 \le (2{\kappa}^2)^{-1}$. We also
find that 
$$\ddt {\varrho}_0(t) = {\gamma}_0m(t)\qquad\text{for almost all $t
\in ]\,t_1,s_1\,[\, \bigcup\, ]\,s_2,t_2\,[$.}$$ Since ${\mu}(t,w) =
{\delta}_0(t,w) + {\varrho}_0(t)$ we find that ${\mu}(t,w) >
{\delta}_0(t,w) \ge {\delta}_0(t_1+,w)$ when $t_1 < t < s_1$, and
${\mu}(t,w) < {\delta}_0( t,w) \le {\delta}_0(t_2-,w)$ when $s_2 < t <
t_2$ which gives ~\eqref{8.31} for these intervals. Thus it only
remains to construct ${\mu}(t,w)$ when $s_1 \le t \le s_2$ and
$|w| \le {\varepsilon}$.

First we define ${\mu}(s_1,w) = {\mu}(s_1-,w)$ and ${\mu}(s_2,w) =
{\mu}(s_2+,w)$ when $|w| \le {\varepsilon}$ and $s_j \ne t_j$.
Then we find for $|w| \le {\varepsilon}$ that
\begin{align}\label{mu1comp}
&{\delta}_0(s_1-,w) \le {\mu}(s_1,w) \le  {\delta}_0(s_1-,w) +
2 {\gamma}_0 M_{12}\\\label{mu2comp}
&{\delta}_0(s_2+,w) \ge {\mu}(s_2,w) \ge  {\delta}_0(s_2+,w) -
2 {\gamma}_0 M_{12}.
\end{align}
Recall that we have defined ${\mu}(s_j,w) = {\delta}_0(\,\wt t_j,w)$
in the case $s_j= t_j$.  Now, we shall define ${\mu}(t,w)$ for $t \in
]\,s_1,s_2\,[$ by linear interpolation using that $m \cong M_{12}$ is
essentially constant and we have a lower bound on the variation:
${\delta}_0(s_2,w)- {\delta}_0(s_1,w) \ge c M_{12}$.

\begin{lem}\label{interpollemma}
Let $m(t)$ be given by ~\eqref{mdef} satisfying ~\eqref{jumpcond1},
~\eqref{unionprop} and ~\eqref{qmaxprop}, let $M_{12}$ be defined by ~\eqref
{9.24} and let $s_1$, $s_2$ be defined by ~\eqref{9.26}--\eqref{9.27}.
Then there exist $s_1=r_1 < r_2 < \dots < r_{N-1} < r_{N}= s_2$, $0 < c_1 <
C_1$ and ${\mu}_j(w) \in \lip(T^*\br^n)$ uniformly  such that
${\mu}(s_1,w) \le {\mu}_1(w) < \dots < {\mu}_N(w) \le {\mu}(s_2,w)$
and 
\begin{align}\label{interpolcond1} 
&c_1M_{12} \le {\mu}_{j+1}(w) - {\mu}_j(w) \le C_1M_{12}\\
&|{\mu}_{j}(w) - {\delta}_0(t,w)| \le C_1M_{12} \qquad\text{for $\quad
  r_j < t < 
r_{j+1}$}\label{interpolcond2} 
\end{align}
when $|w| \le {\varepsilon}$ and $j = 1,
\dots, N-1$ . 
\end{lem}

Observe that Lemma~\ref{interpollemma} also holds when $s_j = t_j$ for
$j=1$ and/or $j =2$. 
We postpone the proof of Lemma~\ref{interpollemma} until later, and
define ${\mu}(t,w)$ by linear interpolation:
\begin{equation*}
{\mu}(t,w) = \frac{t-r_j}{r_{j+1}- r_j}{\mu}_{j+1}(w) +
\frac{t-r_{j+1}}{r_{j} - r_{j+1}}{\mu}_{j}(w)
\qquad r_j < t < r_{j+1}
\end{equation*}
when $|w|\le {\varepsilon}$, $j = 1, \dots, N-1$.
We obtain from~\eqref{interpolcond1} that 
\begin{equation*}
\ddt  {\mu}(t,w) = (r_{j+1}- r_j)^{-1}({\mu}_{j+1}(w)- {\mu}_j(w)) \ge
c_1M_{12}/2 \qquad r_j < t < r_{j+1}
\end{equation*}
when $|w|\le {\varepsilon}$, $j = 1, \dots, N-1$, since $|r_{j+1}-
r_j| \le 2$. We find from ~~\eqref{interpolcond1}
and~\eqref{interpolcond2} that 
\begin{equation*}
|{\varrho}_0(t,w)| = |{\mu}(t,w)- {\delta}_0(t,w)|\le
|{\mu}_j(w)-{\delta}_0(t,w)| + C_1M_{12} \le 2C_1 M_{12}
\qquad r_j < t < r_{j+1}
\end{equation*}
when $|w|\le {\varepsilon}$, $j = 1, \dots, N-1$. By
~\eqref{mcond2} this implies that 
$$|{\varrho}_0(t,w)| \le 2C_1{\kappa}^3m(t)\qquad r_j < t < r_{j+1}
\qquad t \ne r_{12} \qquad |w|\le {\varepsilon}
$$ for $j = 1, \dots, N-1$.  We find that
${\delta}_0(\,\wt t_1,w) \le {\mu}(s_1,w) \le {\mu}(t,w) \le
{\mu}(s_2,w) \le {\delta}_0(\,\wt t_2,w)$ when $s_1 \le t \le s_2$,
which gives ~\eqref{8.31}.
Since $w \mapsto {\mu}_j(w)$ is uniformly
Lipschitz, we obtain that $w \mapsto
{\varrho}_0(t,w) \in \lip(T^*\br^n)$ uniformly when $s_1 < t <
s_2$. Since $t \mapsto {\mu}(t,w)$ is regulated, we
find that $t \mapsto {\varrho}_0(t,w)$ is regulated on $]\,s_1,s_2\,[$
when $|w| \le {\varepsilon}$.
This completes the proof of the Proposition~\ref{apsdef}.
\end{proof}

\begin{proof} [Proof of Lemma~\ref{interpollemma}] 
In order to have control of the variation of ${\delta}_0(t,w)$ in the
$w$ variables, we shall estimate the variation of $ {\Delta}(\,\wt
s,\wt t,w) = {\delta}_0(\,\wt t,w) - {\delta}_0(\,\wt s,w)$  when
$|w| \le {\varepsilon}$. Since
${\varepsilon} \le {\varepsilon}_1$ and $m(\,\wt t\,) \le M_{12} \le 1$ for
$t_1 <t < t_2$, we find by taking limits in Proposition~\ref{dvar1}
that if $t_1 < s < t < t_2$ and
\begin{equation} \label{lambdacond}
|{\Delta}(\,\wt s,\wt t,0)| = {\lambda}M_{12}\qquad {\lambda} \ge 0
\end{equation} 
then 
\begin{equation}\label{ddest1}
 \left(\frac{2{\lambda}}{3} -3\right )M_{12} \le |{\Delta}(\,\wt s,\wt
 t,w)| \le \left(\frac{4{\lambda}}{3}+ 3\right )M_{12}\quad\text{when
 $|w|\le {\varepsilon}$.}
\end{equation}
This also holds if we
replace ~$\wt s$ or $\wt t$ with ~$\wt t_1$ or ~$\wt t_2$ given
by~\eqref{unionprop}, since $m(\,\wt t_j) < M_{12}$ by ~\eqref{qmaxprop}.

Since ${\varepsilon}\le {\varepsilon}_0$ and $M_{12} \ge
{\kappa}\max(m(s_1-), m(s_2+))$ by ~\eqref{9.26}--\eqref{9.27} when
$s_j \ne t_j$, we
obtain in this case that
\begin{equation}\label{dd0est}
{\delta}_0(s_2+,w) - {\delta}_0(s_1-,w) \ge c_0M_{12}
>0\qquad\text{for $ |w| \le {\varepsilon}$}
\end{equation}
by taking the limits in Proposition~\ref{qmaxpropo}. Since
${\kappa}m(\,\wt t_j\,) < M_{12}$ by ~\eqref{qmaxprop}, 
this also holds in the case $s_j = t_j$
if we substitute $\wt t_1$ for $s_1-$ and/or $\wt t_2$ for $s_2+$.
If we have a uniform upper bound then we immediately obtain the result
with $N=2$. If not, we have to divide the interval and use the bounded
variation given by ~\eqref{ddest1}. The complications are the
unbounded jumps of $t \mapsto {\delta}_0(t,w)$ when $s_j = t_j$.

By ~\eqref{jumpcond1} we obtain a bound on the jumps  
\begin{equation}\label{jumpc}
\begin{alignedat}{2} 
&{\delta}_0(t+,0)- {\delta}_0(t-,0) \le 9 M_{12}&\qquad&\text{for $t_1 < t
  < t_2$}\\
&{\delta}_0(t+,w)- {\delta}_0(t-,w) \le 15 M_{12}&\qquad&\text{for $t_1 < t
  < t_2$ and $|w|\le {\varepsilon}$}
\end{alignedat}
\end{equation}
by using~\eqref{ddest1} with ${\lambda} = 9$.
In the case $t = t_j$, we have no bounds on the jumps, but we
shall use a trick to handle that case.  
First we shall define ${\mu}_j(w)$ for $j=1$, 2, such that   
${\mu}(s_1,w) \le {\mu}_1(w) <  {\mu}_2(w)\le {\mu}(s_2,w)$ 
and ${\mu}_2(w) - {\mu}_1(w) \ge c_0/3$ when
$|w|\le {\varepsilon}$. We shall also obtain that
\begin{equation}\label{8.44}
\left\{
\begin{aligned} 
&{\delta}_0(s_1+,0)- 9 M_{12} \le
{\mu}_1(0)\\
&{\mu}_2(0) \le {\delta}_0(s_2-,0)  + 9 M_{12}
\end{aligned}\right.
\end{equation}
and 
\begin{equation}\label{8.45}
\left\{
\begin{aligned} 
&{\delta}_0(s_1+,w) - 15 M_{12} \le {\mu}_1(w) \le {\delta}_0(s_1+,w)
 + M_{12}\\
& {\delta}_0(s_2-,w) - M_{12} \le {\mu}_2(w) \le  {\delta}_0(s_2-,w)
 + 15 M_{12}
\end{aligned}\right.\qquad \text{for $|w|\le {\varepsilon}$}.
\end{equation}

If $s_1 > t_1$ then we put ${\mu}_1(w) =
{\delta}_0(s_1-,w) + {\gamma}_1M_{12}$ 
where ${\gamma}_1 = \min(1,c_0/3)$. 
By~\eqref{mu1comp} we find that ${\mu}_1(w) \ge  {\mu}(s_1,w)$
if ${\gamma}_0 \le c_0/6$. 
By ~\eqref{jumpc}
we obtain \eqref{8.44}--\eqref{8.45} in this case. 

If $s_1 = t_1$ and ${\delta}_0(t_1+,0) - {\delta}_0(\,\wt t_1,0) \le
9M_{12}$ then we obtain as before from ~\eqref{ddest1} that
${\delta}_0(t_1+,w) - {\delta}_0(\,\wt t_1,w) \le 15 M_{12}$ for
$|w|\le {\varepsilon}$. If we put ${\mu}_1(w) = {\delta}_0(\,\wt
t_1,w) = {\mu}(t_1,w)$, we obtain \eqref{8.44}--\eqref{8.45} in this
case.

In the last case when $s_1 = t_1$ and
${\delta}_0(t_1+,0) - {\delta}_0(\,\wt t_1,0) > 9M_{12}$, we find
by using~\eqref{ddest1} with ${\lambda} = 9$ that
$ {\delta}_0(t_1+,w) - {\delta}_0(\,\wt t_1,w) > 3
M_{12}$ when $|w|\le {\varepsilon}$.
In that case we let 
${\mu}_1(w) = {\delta}_0(t_1+,w) - 3M_{12} \ge  {\delta}_0(\,\wt
t_1,w) = {\mu}(t_1,w)
$ for $|w|\le {\varepsilon}$, which gives
\eqref{8.44}--\eqref{8.45} for $j =1$.

Similarly, if $s_2 \ne t_2$ then we put ${\mu}_2(w) =
{\delta}_0(s_2+,w) - {\gamma}_1M_{12} \le {\mu}(s_2,w)$
by~\eqref{mu2comp}. If $s_2 = 
t_2$ and ${\delta}_0(\,\wt t_2,0) - {\delta}_0(t_2-,0) \le
9M_{12}$ then we put ${\mu}_2(w) = {\delta}_0(\,\wt t_2,w) =
{\mu}(t_2,w)$. Finally, when $s_2 = t_2$ and ${\delta}_0(\,\wt
t_2,0) - {\delta}_0(t_2-,0) > 9M_{12}$, then we let ${\mu}_2(w) =
{\delta}_0(t_2-,w) + 3M_{12}\le {\mu}(t_2,w)$, and we obtain
as before \eqref{8.44}--\eqref{8.45} for $j=2$. By the definition of
${\mu}_j(w)$ we obtain that  ${\mu}_2(w) - {\mu}_1(w) \ge \min(3,c_0/3)$.

We are going to consider the value of 
\begin{equation}\label{kdef}
K =  \left({\mu}_2(0) - {\mu}_1(0)\right)/M_{12} \ge
\min(3,c_0/3). 
\end{equation}
Now we have no fixed upper bound on~
$K$, and therefore we shall consider the cases when $K \gtrless 45$.

In the case $K \le 45$ we find 
that 
${\delta}_0(s_2-,0) - {\delta}_0(s_1+, 0) \le 47 M_{12}$ 
by ~\eqref{8.45}. We obtain that 
${\delta}_0(s_2-,w) - {\delta}_0(s_1+, w) \le 197 M_{12}/3$ 
when $|w|\le{\varepsilon}$ by taking ${\lambda} = 47$ in~\eqref{ddest1}.
This gives ${\mu}_2(w)-{\mu}_1(w) \le (30 + 197/3)M_{12} < 96M_{12}$
for $|w| \le {\varepsilon}$.
By ~\eqref{8.45} we find that 
$|{\mu}_{1}(w) - {\delta}_0(t,w)| \le 97 M_{12}$ when $s_1 < t <
s_2$ and $|w| \le {\varepsilon}$.
Thus, we obtain the result in this case with  $N = 2$,
$c_1 = \min(3, c_0/3)$ and $C_1 = 97$.

Next, we consider the case $K > 45$. Then we obtain that $K_1 =
({\delta}_0(s_2-,0) - {\delta}_0(s_1+,0))/M_{12} \ge K -18 > 27$ by
~\eqref{8.44}. 
By the jump condition~~\eqref{jumpc} we find that $]s_1, s_2[\ \ni t \mapsto
({\delta}_0(t,0) -  {\delta}_0(s_1+,0))/M_{12}$ takes values in any closed
interval of length $9$ in $[0, K_1]$. Thus, we can find $r_2 \in\ ]s_1,
s_2[ $ such that $({\delta}_0(r_2,0)
- {\delta}_0(s_1+,0))/M_{12} \in [9, \,18]$, and we find that 
$K_2 =
({\delta}_0(s_2-,0) - {\delta}_0(r_2,0))/M_{12} \ge K_1 -18 > 9$. If
recursively  $K_j =
({\delta}_0(s_2-,0) - {\delta}_0(r_j,0))/M_{12} > 27$
then we choose $r_{j+1} \in\ ]r_j, 
s_2[ $ such that $({\delta}_0(r_{j+1},0)
- {\delta}_0(r_j,0))/M_{12} \in [9, \,18]$, until $9<
K_{N-1} \le 27$.
Then by using ~\eqref{ddest1} with $9 \le {\lambda} \le 18$
we find that $({\delta}_0(r_2,w)
- {\delta}_0(s_1+,w))/M_{12}$ and $({\delta}_0(r_{j+1},w)-
{\delta}_0(r_j,w))/M_{12} \in [3,27]$ for $1 < j < N-1$, and similarly
we find that $({\delta}_0(s_2-,w) - {\delta}_0(r_{N-1},w))/M_{12} \in [3, 39]$.
Putting ${\mu}_N(w) = {\mu}_2(w)$, redefining ${\mu}_2(w)
= {\delta}(r_2,w)$ and letting ${\mu}_j(w) = {\delta}_0(r_j,w)$ for $2
< j < N$, we find from~\eqref{8.45} that
$ 
2M_{12} \le {\mu}_{j+1}(w) - {\mu}_j(w)\le
54 M_{12}$ when $|w| \le {\varepsilon}$ and  $1 \le  j \le N-1
$. 
By the construction and ~\eqref{8.45} we find that 
$ 
|{\mu}_{j}(w) - {\delta}_0(t,w)| \le 42 M_{12}$ when $r_j <
t < r_{j+1}$ and $|w| \le {\varepsilon}
$ 
for $ 1 \le j \le N-1$.  Thus we obtain the result in this case with $c_1 =
2$ and $C_1 = 54$, which completes the proof of the Lemma.
\end{proof}

Now we shall compute the Weyl symbol for the Wick operator given by
the ``pseudo-sign'' ${\delta}_0 +
{\varrho}_0$.   
As before, we shall use the symbol classes $S^+(1,g^\sh)$ given by
Definition ~\ref{s+def}.

\begin{prop}\label{wickweyl}
Let $B = {\delta}_0 + {\varrho}_0$, where ${\delta}_0$ is given by
Definition~~\ref{d0deforig} and ${\varrho}_0(t,w)$ is the
real valued symbol given by Proposition~\ref{apsdef},
satisfying $|{\varrho}_0(t,w)| \le Cm_{1}(t,w)$ for almost all $|t|\le
1$. Then we find
\begin{equation*}
B^{Wick} = b^w \qquad |t| \le 1
\end{equation*}
where $b= {\delta}_1 + {\varrho}_1 \in S(H_1^{-1/2},g^\sh)\bigcap
S^+(1,g^\sh)$ is real valued and regulated in $t$, and ${\varrho}_1
\in S(m_{1}, g^\sh) \subseteq S(H_1^{1/2}\w{{\delta}_0}, g^\sh)$ for
almost all $|t| \le 1$.  There also exists a positive constant
${\kappa}_2$ with the following properties.  For any ${\lambda} > 0$,
there exists $c_{\lambda} > 0$ such that if $|{\delta}_0| \ge
{\lambda} H_1^{-1/2}$ and $H_1^{1/2} \le c_{\lambda}$ then $|b| \ge
{\kappa}_2{\lambda}H_1^{-1/2}$.  If $H_1^{1/2}(t,w_0) \le {\kappa}_2$
and $|{\delta}_0(t,w_0)| \le {\kappa}_2H_1^{-1/2}(t,w_0)$ then we have
$S(H_1^{-1/2},G_1) \ni {\delta}_1(t,w) = {\delta}_0(t,w) +
{\varrho}_2(t,w)$ when $|w-w_0| \le {\kappa}_2H_1^{-1/2}(t,w_0)$ with
real valued $ {\varrho}_2(t,w) \in S(H_1^{1/2},G_1)$.
\end{prop}

\begin{proof}
Let ${\delta}_0^{Wick} = {\delta}_1^w$ and ${\varrho}_0^{Wick} =
{\varrho}_1^w$.  Since $|{\delta}_0|\le CH^{-1/2}_1$, $|{\varrho}_0|
\le C m_{1}$ and the symbols are real valued, we obtain from
Proposition~\ref{propwick} and ~\eqref{hhhest} that ${\delta}_1 \in
S(H_1^{-1/2},g^\sh)$ and ${\varrho}_1 \in S(m_{1}, g^\sh) \subseteq
S(H_1^{1/2}\w{{\delta}_0}, g^\sh)$ are real valued for almost all $|t|
\le 1$. Observe that $m_{1} \le 1$, and since $|{\delta}_0'| \le 1$
almost everywhere we find that $b \in S^+(1,g^\sh)$ for almost all
$|t| \le 1$ by Proposition~\ref{propwick}. Since ${\delta}_0(t,w)$ and
${\varrho}_0(t,w)$ are regulated in~ $t$, we find from~\eqref{gausreg}
that the same holds for ${\delta}_1(t,w)$ and~${\varrho}_1(t,w)$.

When $|{\delta}_0| \ge {\lambda}H_1^{-1/2}$ at $(t,w)$, ${\lambda} >
0$, then by the Lipschitz continuity and slow variation we find that
$|{\delta}_0| \ge {\lambda}H_1^{-1/2}/C_0$ in a $G_1$ neighborhood
${\omega}$ of $(t,w)$ (depending on ${\lambda}$). Since $|{\varrho}_0|
\le CH_1^{1/2}\w{{\delta}_0}$ we find that $|{\delta}_0 + {\varrho}_0|
\ge {\lambda}H_1^{-1/2}/2C_0$ when $H_1^{1/2}$ is small enough in
~${\omega}$. By the slow variation, it suffices that $H^{1/2}_1(t,w)$
is small enough if the neighborhood ${\omega}$ is sufficently small.
Proposition~\ref{propwick} gives $|b| \ge c
{\lambda}H_1^{-1/2}/2C_0 - C{\lambda}H_1^{1/2}/2C_0 \ge
c{\lambda}H_1^{-1/2}/3C_0$ at $(t,w)$ when $H_1^{1/2}(t,w)$ is small
enough.

If $|{\delta}_0| \le {\kappa}_2H_1^{-1/2}$ and $H_1^{1/2} \le {\kappa}_2$ for
sufficiently small ${\kappa}_2 >0$, then $|{\delta}_0| \le
C_0{\kappa}_2H_1^{-1/2}$ 
and $H_1^{1/2} \le C_0{\kappa}_2$ in a fixed $G_1$ neighborhood. Thus,
for ${\kappa}_2 \ll 1$ we obtain that ${\delta}_0 \in S(H_1^{-1/2},
G_1)$ in a fixed $G_1$ neighborhood. 
Then we obtain the last statement from Proposition
~\ref{propwick}, which completes the proof.
\end{proof}

\section{The Lower Bounds}\label{lower}

In this section we shall show that $B^{Wick} = b^w$ given by
Proposition~\ref{wickweyl} satisfies the conditions in
Proposition~\ref{apsprop}, finally proving that
Proposition and completing the proof of the
Nirenberg-Treves conjecture. Recall that
$f \in C(\br, S(h^{-1}, hg^\sh))$, where both $0 <
h \le 1$ and $g^\sh = (g^\sh)^{\sigma}$ are constant.
First, we shall obtain lower
bounds on $\re b^wf^w$.

\begin{prop}\label{lowersign}
Assume that
$b = {\delta}_1 + {\varrho}_1$ 
is given by Proposition~\ref{wickweyl}.
Then we have
\begin{equation}\label{lowerbound}
\re\w{(b^wf^w)\restr{t}u,u} \ge \w{C_t^wu,u} \qquad\forall\  u \in
C_0^\infty(\br^n)\qquad\text{for almost all $|t| \le 1$}
\end{equation}
where $C_t \in S(m_{1}(t), g^\sh)$  has uniformly bounded
seminorms, which only depend on the seminorms of $f$
in $S(h^{-1},hg^\sh)$ for $|t| \le 1$.
\end{prop}

\begin{proof}
We shall prove the Proposition by localizing with respect to the
metric $G_1$ for fixed ~$t$. Observe that we may ignore terms
in $\op S(MH_1^{3/2}\w{{\delta}_0},g^\sharp) \subseteq
\op S(m_{1}, g^\sharp)$ by Proposition~\ref{mestprop}.
We fix $|t| \le 1$ such that $b = {\delta}_1 + {\varrho}_1 \in
S(H_1^{-1/2},g^\sh) + S(m_1,g^\sh)$ for this ~$t$.   
In the following we shall omit the $t$ variable and put $H_1(w)
= H_1(t,w)$, $m_{1}(w) =m_{1}(t,w)$ and $M(w) =
M(t,w)$. We shall localize with respect to the metric $G_1 =
H_1g^\sh$, and as before we shall assume the coordinates chosen so
that $g^\sh(w) = |w|^2$. In the following, we shall use the
neighborhoods 
$${\omega}_{w_0}({\varepsilon}) = \set{w:\ |w-w_0| \le
{\varepsilon}H_1^{-1/2}(w_0)}\qquad w_0 \in T^*\br^n.$$

By the slow variation of $G_1$ and the uniform Lipschitz continuity of
$w \mapsto {\delta}_0(w)$ we find that there exists ${\kappa}_0 >0$
with the following property. If $0 < {\kappa}\le {\kappa}_0$ then there
exist positive constants $c_{\kappa}$ and ${\varepsilon}_{\kappa}$ so
that for any $w_0 \in T^*\br^n$ we have
\begin{alignat}{2}
&|{\delta}_0(w)| \le {\kappa}H_1^{-1/2}(w)&\qquad &w \in
{\omega}_{w_0}({\varepsilon}_{\kappa})  \qquad\text{or}\label{d0case}\\
&|{\delta}_0(w)| \ge c_{\kappa} H_1^{-1/2}(w)&\qquad &w \in
{\omega}_{w_0}({\varepsilon}_{\kappa}). \label{d1case}
\end{alignat}
We may also assume that 
${\varepsilon}_{\kappa}$ is small enough so that 
$w \mapsto H_{1}(w)$ and $w \mapsto M(w)$ only vary with a fixed factor
in~ ${\omega}_{w_0}({\varepsilon}_{\kappa})$.  In fact, we have by the 
Lipschitz continuity that $w \mapsto {\delta}_0(w)$ varies with at most
$2{\varepsilon_{\kappa}}H_1^{-1/2}(w_0)$ in
${\omega}_{w_0}({\varepsilon}_{\kappa})$, thus if
${\varepsilon}_{\kappa} \ll {\kappa}$ we obtain 
that ~\eqref{d0case} holds when $|{\delta}_0(w_0)| \ll
{\kappa}H_1^{-1/2}(w_0)$ and ~\eqref{d1case} holds when
$|{\delta}_0(w_0)| \ge c{\kappa}H_1^{-1/2}(w_0)$.

Now we let ${\kappa}_1$ be given by Proposition~\ref{ffactprop},
${\kappa}_2$ by Proposition~\ref{wickweyl}, choose ${\kappa} = \min
({\kappa}_0,{\kappa}_1, {\kappa}_2)$ and
let ${\varepsilon}_{\kappa}$ and $c_{\kappa}$ be given by
\eqref{d0case}--\eqref{d1case}. 
Since $H_1$ only varies with a 
fixed factor in $ {\omega}_{w_0}({\varepsilon}_{\kappa})$,
Proposition~\ref{wickweyl} (with ${\lambda} = c_{\kappa}$) gives
${\kappa}_3 > 0$ such that 
\begin{equation} \label{kappa3ref}
|b|\ge {\kappa}_2c_{\kappa}H_1^{-1/2}\qquad\text{in ~$
{\omega}_{w_0}({\varepsilon}_{\kappa})$} \quad\text{if
$H_1^{1/2}(w_0) \le {\kappa}_3$ and ~\eqref{d1case} holds in ~$
{\omega}_{w_0}({\varepsilon}_{\kappa}) $.}
\end{equation}  
Take a partition of unity $\set{{\psi}_k(w)}_k$, $\set{{\Psi}_k(w)}_k$
and $\set{{\Phi}_k(w)}_k \in S(1,G_{1})$ with values in ~~$\ell^2$,
such that $0 \le {\psi}_k \le 1$, $0 \le {\Psi}_k \le 1$, $0 \le
{\Phi}_k \le 1$, $\sum_{k} {\psi}_k^2 \equiv 1$, ${\psi}_k{\Psi}_k =
{\psi}_k$, ${\Psi}_k{\Phi}_k = {\Psi}_k$, ${\Psi}_k = {\phi}_k^2$ for
some $\set{{\phi}_k(w)}_k \in S(1,G_{1})$ and $\supp {\Phi}_k
\subseteq {\omega}_k = {\omega}_{w_k}({\varepsilon}_{\kappa})$.
Observe that $\set{{\omega}_k}_k$ cover $T^*\br^n$ since $\sum_j {\psi}_j^2
\equiv 1$. Let $A = bf$ and $A_k = {\Psi}_kA \in
S(MH_1^{-1/2},g^\sharp)\bigcap S^+(M,g^\sharp)$
uniformly, which are real valued symbols.  Next, we shall localize
$\re b^wf^w$ using the following Lemma.

\begin{lem}\label{bflem}
Assume that $A = bf$ where $b = {\delta}_1 + {\varrho}_1$ is given by
Proposition~\ref{wickweyl}, and let $A_k = {\Psi}_kA$. Then we find
\begin{equation}\label{bfcong}
\re (b^wf^w) \cong A^w \cong \sum_{k} {\psi}_k^w A_k^w {\psi}_k^w
\end{equation}
modulo $\op S(MH_{1}^{3/2}\w{{\delta}_0},g^\sharp) \subseteq \op
S(m_{1},g^\sharp)$.
\end{lem}

\begin{proof} 
We shall consider $\set{{\psi}_k}_k$ as having values in $\ell^2$ and
$\set{A_k}_k =\set{A_k{\delta}_{jk}}_{jk}$ as a diagonal matrix in
$\cal L(\ell^2,\ell^2)$. Then we can use the calculus as for scalar
valued symbols (see ~\cite[p.\ 169]{ho:yellow}), and we shall compute
the symbols of~\eqref{bfcong}.  Observe that in the domains
~${\omega}_k$ where $H_1^{1/2} \ge c > 0$, we find that $A_k
\in S(MH_1^{3/2}, g^\sharp)$ giving the result in this case. Thus in
the following, 
we shall assume that $H_1^{1/2} \le
{\kappa}_2$ in~${\omega}_k$, with ${\kappa}_2$ given by
Proposition~\ref{wickweyl}, and we shall consider the cases
where~\eqref{d0case} or~\eqref{d1case} hold.

First we consider the case when~\eqref{d1case} holds. Then we find that
$\w{{\delta}_0} \cong H_1^{-1/2}$ in~${\omega}_k$ so $S(MH_1,g^\sharp) =
S(MH_1^{3/2}\w{{\delta}_0},g^\sharp) \subseteq S(m_{1}, g^\sharp)$ 
in~${\omega}_k$. Since $g^\sh/G_1^{\sigma} = H_1$, we find from
Theorem~18.5.5 in ~\cite{ho:yellow} that the 
symbol of $b^wf^w \in \op S(MH_1^{-1/2},g^\sharp)$ has an expansion in
$S(MH_1^{(j-1)/2},g^\sharp)$, $j \ge 0$. Moreover, since $b \in
S^+(1,g^\sharp)$ we find from ~\thetag{18.4.8} in~\cite{ho:yellow} that
the symbol of $
 b^wf^w $ is equal to $
bf + \frac{1}{2i}\pb{b,f}
$
modulo $S(MH_1, g^\sharp)$.  Since $\re(a^w) = \frac{1}{2}(a^w +
\overline a^w) = (\re a)^w$ we find that the symbol of $\re (b^wf^w)
$  is equal to  $bf$ in~${\omega}_k$
modulo $S(MH_1,g^\sharp)$.  Similarly, since $\sum_{k} {\psi}_k^w
A_k^w {\psi}_k^w$ is symmetric and $A_k \in S^+(M,g^\sh)$, we find
that the symbol of $\sum_{k} {\psi}_k^w A_k^w {\psi}_k^w $ is equal to
$ A$ in~${\omega}_k$ modulo
$S(MH_1,g^\sh)$, which proves the result in this case.

Finally, we consider the case when ~\eqref{d0case} holds and
$H_1^{1/2} \le {\kappa}_2$ in~${\omega}_k$. Then $b = {\delta}_1 +
{\varrho}_1 \in S(H_1^{-1/2},G_1)+ S(m_{1}, g^\sharp)$ in~${\omega}_k$
by Proposition ~\ref{wickweyl}. Thus, $b^wf^w = {\delta}_1^wf^w
+{\varrho}_1^wf^w$ where the symbol of ${\delta}_1^wf^w$ has an
expansion in $S(MH_1^{-1/2}H_1^{j},G_1)$ and as before the symbol of
${\varrho}_1^wf^w$ has an expansion in ~$S(m_1MH_1^{j/2},g^\sh)$
in~${\omega}_k$. By taking the symmetric part we obtain only even $j$,
and since $M \le CH_1^{-1}$ we obtain that the symbol of $\re
(b^wf^w)$ is in $S(m_1,g^\sh)$ in ~${\omega}_k$. Similarly, since $A_k
\in S(MH_1^{-1/2},G_1)+ S(Mm_{1},g^\sharp)$ is real, we find that the
symbol of $\sum_{k} {\psi}_k^w A_k^w {\psi}_k^w$ is equal to $A$
modulo $ S(m,g^\sh)$ in ~${\omega}_k$, which proves ~\eqref{bfcong}
and the Lemma.
\end{proof}

By Lemma~\ref{bflem} it suffices to get lower bounds on $A_k^w$.
We obtain from Lemma ~\ref{lowerlem} below that $\w{A_k^w{\psi}_k^w
u,{\psi}_k^w u} \ge \w{C^w_{k}{\psi}_k^w u, {\psi}_k^wu}$, where
$C_k\in S(m_{1}, g^\sharp)$ uniformly in $k$. Thus we obtain 
from Lemma~\ref{bflem} that 
$$\re\w{b^wf^w u,u} \ge
\sum_{k} \w{{\psi}_k^wC^w_{k}{\psi}_k^w u,u} +\w{R^wu,u}\qquad u \in
C_0^{\infty}({\mathbf R}^n)
$$ where $\sum_{k}{\psi}_k^wC^w_{k}{\psi}_k^w$ and $R^w \in \op
S(m_{1},g^\sharp)$. This completes the proof of
Proposition~\ref{lowersign}.
\end{proof}

\begin{lem}\label{lowerlem}
Let $A = bf$ where $b = {\delta}_1 + {\varrho}_1$ is given by
Proposition~\ref{wickweyl}, and let $A_k = {\Psi}_kA$.
There exists $C_{k} \in S(m_{1}, g^\sharp)$ uniformly,
such that
\begin{equation*}
\w{A_k^w u,u} \ge \w{C^w_{k} u,u}\qquad
u \in C_0^{\infty}({\mathbf R}^n)\qquad \forall\, k.
\end{equation*} 
\end{lem}

\begin{proof}
We shall keep $k$ fixed and as before we are going to consider the
cases when $H_1^{1/2} \cong 1$ or $H_1^{1/2} \ll 1$, and when
~\eqref{d0case} or ~\eqref{d1case} holds in ${\omega}_k$.

Assume that $H_1^{1/2}(w_k) \ge c > 0$, then $A_k \in
S(MH_1^{-1/2},g^\sharp) \subseteq S(MH_1^{3/2},g^\sharp) \subseteq
S(m_{1},g^\sharp)$ uniformly by 
Proposition~\ref{mestprop}.  Thus, we obtain the Lemma with $C_k =
A_k$ in this case. Thus we may assume that  $H_1^{1/2}(w_k) \le {\kappa}_4
= \min({\kappa}_0,{\kappa}_1, {\kappa}_2, {\kappa}_3)$ in
what follows.

Next, we consider the case when ~\eqref{d0case} holds and
$H_1^{1/2}(w_k) \le {\kappa}_4 $ in ${\omega}_k$. Then we obtain from
Proposition~\ref{wickweyl} that $b = {\delta}_0 + {\varrho}_1 +
{\varrho}_2 = {\delta}_0 + {\varrho}_3$ is real valued, where
${\varrho}_1 \in S(m_{1}, g^\sharp)$ and ${\varrho}_2 \in S(H_1^{1/2},
G_1)$ in ~~${\omega}_k$. We find from Proposition~\ref{ffactprop} that
we may choose $g^\sharp$~orthonormal coordinates so that $|w|\le
cH_1^{-1/2}(w_k)$ and
\begin{equation*}
\left\{
\begin{aligned}
&f(w) = {\alpha}_1(w)(w_1- {\beta}(w'))\\
&{\delta}_0(w) = {\alpha}_0(w)(w_1- {\beta}(w'))
\end{aligned}
\right.
\qquad\text{in ~${\omega}_k$. }
\end{equation*}
Here $ c\le {\alpha}_0 \in S(1, G_1)$
and $c|f'| \le {\alpha}_1\in S(|f'|, G_1)$ for some
$c>0$, and ${\beta} \in S(H_1^{-1/2}, G_1)$.
Since ${\Psi}_k =
{\phi}_k^2$, we find that 
\begin{equation*} 
A_k(w) = {\Psi}_kbf = {\phi}_k^2(w)\big({\alpha}_0(w){\alpha}_1(w) (w_1-
{\beta}(w'))^2 + {\alpha}_1(w)(w_1- {\beta}(w')){\varrho}_3(w)\big),
\end{equation*}
and we shall construct an approximate square root to $A_k$ by
completing the square. Let
\begin{multline*} 
{\gamma}_k(w) = {\phi}_k(w)\sqrt{{\alpha}_0(w) {\alpha}_1(w)}\big(w_1
-{\beta}(w') + {\varrho}_3(w)/2{\alpha}_0(w)\big) \\ \in
S(|f'|^{1/2}H_1^{-1/2}, G_1) + S(|f'|^{1/2}m_{1},g^\sharp),
\end{multline*}
uniformly, which is real valued since ${\varrho}_3 = {\varrho}_1 +
{\varrho}_2 \in S(H_1^{1/2}, G_1) + S(m_{1}, g^\sharp)$ is real.  Then
we find
\begin{equation*}
 {\gamma}_k^2 = A_k + {\Psi}_k{\alpha}_1{\varrho}_3^2/4{\alpha}_0 
\end{equation*}
where $ {\Psi}_k{\alpha}_1{\varrho}_3^2/4{\alpha}_0 \in S(|f'|
m^2_{1},g^\sharp) + S(m_{1},g^\sharp)$,
since we have $|f'|H_1 \le MH_1^{3/2} \le Cm_{1}$ and $|f'| \le
CH_1^{-1/2}$ by~\eqref{dfest0}. As before,
~\cite[Theorem~18.5.5]{ho:yellow} gives 
\begin{equation}\label{9.6a} 
A_k^w \cong {\gamma}_k^w{\gamma}_k^w + R_k^w \qquad \text{modulo $\op
  S(m_{1},g^\sharp)$}
\end{equation}
where the term $R_k \in S(|f'| m^2_{1},g^\sharp) \subseteq
S(|f'|H_1\w{{\delta}_0}^2, g^\sharp)$ by ~\eqref{hhhest} is real
valued. In fact, the composition
of operators in $ \op S(|f'|^{1/2}H_1^{-1/2}, G_1)$ and $\op
S(|f'|^{1/2}m_{1},g^\sharp)$ gives an expansion in $S(|f'|m_{1}
H_1^{(j-1)/2}, g^\sharp)$, with even ~$j$ in this case since the
operator is symmetric. In order to treat the error term $R_k$, we shall
localize where $|{\delta}_0| \gtrless 1$, for example with
${\phi}({\delta}_0) \in S(1,g^\sh)$ where ${\phi} \in C_0^{\infty}(\br)$.
Since $|R_k| \le C|f'|H_1 \le C_0m_1$ when $|{\delta}_0| \le c$ and
$ c_0|f'| \le {\alpha}_1{\alpha}_0^{-1} \in S(|f'|, G_1)$ in
${\omega}_k$, we find that 
$  R_k \cong {\chi}_k{\Psi}_k
 {\alpha}_1{\alpha}_0^{-1}{\delta}_0^2$ modulo 
   $S(m_{1},g^\sharp) 
$ 
where ${\chi}_k \in S(H_1,g^\sharp)$ is real valued and supported
where $|{\delta}_0| \ge 1$. 
Now we have ${\gamma}_k^2 \cong A_k \cong
{\Psi}_k{\alpha}_1{\alpha}_0^{-1}{\delta}_0^2$ modulo
$S(|f'|\w{{\delta}_0},g^\sharp)$. Thus we obtain that 
\begin{equation}\label{9.6b}
 {\chi}_k{\gamma}_k^2 \cong  R_k\qquad \text{modulo
   $S(|f'|H_1\w{{\delta}_0},g^\sharp) \subseteq S(m_1,g^\sh)$} 
\end{equation}
since $|f'|H_1\w{{\delta}_0} \le MH_1^{3/2}\w{{\delta}_0} \le Cm_1$.
Recall that the symbol ${\gamma}_k  \in
S(|f'|^{1/2}H_1^{-1/2}, G_1) + S(|f'|^{1/2}m_{1},g^\sharp)$.
By taking the real valued
\begin{equation*}
\widetilde  {\gamma}_k = (1 + {\chi}_k/2){\gamma}_k  \in
S(|f'|^{1/2}H_1^{-1/2}, G_1) + S(|f'|^{1/2}H_1^{1/2}\w{{\delta}_0},g^\sharp)
\end{equation*}
we obtain that $\wt {\gamma}_k \cong {\gamma}_k$ modulo
${\chi}_k{\gamma}_k/2 \in S(|f'|^{1/2}H_1^{1/2},g^\sh)$. Thus we find
that
\begin{equation*}
  \widetilde {\gamma}_k^w \widetilde {\gamma}_k^w \cong {\gamma}_k^w 
{\gamma}_k^w + \frac{1}{2}(({\chi}_k{\gamma}_k)^w{\gamma}_k^w +
{\gamma}_k^w ({\chi}_k{\gamma}_k)^w)\qquad\text{modulo $S(m_1,g^\sh)$.}
\end{equation*}
In order to handle the last terms, we observe that
the composition of operators in $\op S(|f'|^{1/2}H_1^{1/2}, g^\sh)$ and
$\op S(|f'|^{1/2}H_1^{1/2}\w{{\delta}_0},g^\sharp)$ gives operators in $
\op S(m_1,g^\sharp)$.
Also, the composition of operators in $\op S(|f'|^{1/2}H_1^{-1/2}, G_1)$ and
$\op S(|f'|^{1/2}H_1^{1/2},g^\sharp)$ gives symbol expansions in $S(|f'|
H_1^{j/2}, g^\sharp)$ with even ~$j$ for symmetric operators.
Thus, we obtain that $\frac{1}{2}(({\chi}_k{\gamma}_k)^w{\gamma}_k^w +
{\gamma}_k^w ({\chi}_k{\gamma}_k)^w) \cong ({\chi}_k{\gamma}_k^2)^w$
modulo $\op S(m_1,g^\sharp)$. By \eqref{9.6a}--\eqref{9.6b} we find that
\begin{equation*}
A_k^w \cong {\gamma}_k^w{\gamma}_k^w + ({\chi}_k{\gamma}^2_k)^w \cong
\widetilde {\gamma}_k^w \widetilde {\gamma}_k^w \cong (\widetilde
{\gamma}_k^w)^* \widetilde {\gamma}_k^w \ge 0 \qquad \text{modulo $\op
S(m_{1},g^\sharp)$},
\end{equation*}
which proves the Lemma in this case.

Finally, we consider the case when $H_1^{1/2}(w_k) \le {\kappa}_4$ and
~\eqref{d1case} holds in ${\omega}_k$.  Then, we shall use the uniform
Fefferman-Phong estimate, but we also have to consider the
perturbation ${\varrho}_1(w)$, which is in a bad symbol class.  Since
$|{\delta}_0(w)| \ge c_{\kappa} H_1^{-1/2}(w)$, we find $|{\delta}_0|
\cong H_1^{-1/2}$ in ~${\omega}_k$.  Thus, we may ignore terms in
$S(MH_1,g^\sharp) \subseteq S(MH_1^{3/2}\w{{\delta}_0},g^\sharp)$
supported in ${\omega}_k$.  Since $H_1^{1/2}(w_k) \le {\kappa}_3$ and
$|{\delta}_0| \ge c_{\kappa}H_1^{-1/2}$, \eqref{kappa3ref} gives that
$|b|\ge {\kappa}_2c_{\kappa}H_1^{-1/2}$ in 
${\omega}_k$.  Since $b \in S^+(1,g^\sh)$, we find by the chain rule
that $|b|^{{\lambda}} \in S(H_1^{-{\lambda}/2},g^\sharp)\bigcap
S^+(H_1^{(1-{\lambda})/2},g^\sharp)$ in ~~${\omega}_k$. In fact, we have
$\partial_w |b|^{{\lambda}} = \sgn(b){\lambda}|b|^{{\lambda}-1}
\partial_w b \in S(H_1^{(1-{\lambda})/2},g^\sharp)$ in~~${\omega}_k$
since $\partial_w 
b \in S(1,g^\sh)$. Let $0 \le {a}_k =
{\Psi}_k|f| \in S(M,G_1)$ 
and ${\beta}_k = {\Phi}_k|b|^{1/2} \in S(H_1^{-1/4},g^\sharp)\bigcap
S^+(H_1^{1/4},g^\sharp)$ since $\partial {\Phi}_k = \cal
O(H_1^{1/2})$, then we obtain that $A_k = {\Psi}_kbf = {a}_k{\beta}_k^2$.  By
using Lemma~\ref{calclem} below, we find that
\begin{equation*}
A_k^w \cong  {\beta}_k^w {a}_k^w{\beta}_k^w - ({a}_kr_k)^w
\qquad \text{}
\end{equation*}
modulo $S(MH_1,g^\sharp) \subseteq S(m_1,g^\sh)$,
with real $r_k \in S(H_1^{1/2},g^\sharp)$. 
By taking the real symbol
$${\lambda}_k = {\beta}_k + {\Phi}_kr_k|b|^{-1/2}/2 \in
S(H_1^{-1/4},g^\sharp)\bigcap S^+(H_1^{1/4},g^\sharp) +
S(H_1^{3/4},g^\sharp)
$$
we find that  
$ 
 A_k^w \cong {\lambda}_k^w{a}_k^w{\lambda}_k^w
$ 
modulo $\op S(MH_1,g^\sharp)$. In fact, since ${a}_k \in S(M,G_1)$ and 
${\Phi}_k r_k|b|^{-1/2} \in
S(H_1^{3/4},g^\sharp)$ we obtain that 
that ${a}_k^w({\Phi}_k r_k|b|^{-1/2})^w \cong ({a}_k
r_k|b|^{-1/2})^w \in \op S(MH_1^{3/4},g^\sh)$
modulo $\op S(MH_1^{5/4},g^\sh)$. Since ${\beta}_k = {\Phi}_k|b|^{1/2} \in
S^+(H_1^{1/4},g^\sharp)$ and  ${\Phi}_k = 1$ on $\supp {\Psi}_k$, 
we find that ${\beta}_k^w {a}_k^w
({\Phi}_k r_k|b|^{-1/2})^w \cong ({a}_k r_k)^w$ modulo
$\op S(MH_1, g^\sharp)$ by using  ~\thetag{18.4.8} in ~\cite{ho:yellow}.  
By the uniform Fefferman-Phong estimate
\cite[Lemma~18.6.10]{ho:yellow}, there exists $C>0$ so that
\begin{equation*}
\w{{a}_k^w u,u} \ge -CM(w_k)H_1^2(w_k)\mn{u}^2 \qquad \forall\ u \in
C_0^{\infty}({\mathbf R}^n). 
\end{equation*}
Since ${\lambda}_k$ is real valued this gives
\begin{equation*}
  \w{{\lambda}_k^w{a}_k^w{\lambda}_k^w u,u} \ge -C M(w_k)
  H_1^{2}(w_k)\mn{{\lambda}_k^w u}^2 = \w{c_k^w u,u}
\end{equation*}
where $c_k^w = -CM(w_k)H_1^2(w_k){\lambda}_k^w{\lambda}_k^w \in \op
S(MH_1^{3/2},g^\sharp)$ uniformly. This completes the proof 
of Lemma~\ref{lowerlem}. 
\end{proof}

\begin{lem}\label{calclem}
Let ${\beta}_k \in S(H_1^{-1/4},g^\sharp)\bigcap
S^+(H_1^{1/4},g^\sharp)$ and ${a}_k \in S(M,G_1)$ be real valued
symbols. Then there exists a real valued symbol $r_k \in
S(H_1^{1/2},g^\sharp)$ such that
\begin{equation}\label{9.8}
  {\beta}_k^w {a}_k^w{\beta}_k^w \cong
  \left({a}_k({\beta}_k^2 + r_k)\right)^w
\end{equation}
 modulo $\op S(MH_1,g^\sharp)$.
\end{lem}

\begin{proof}
We have ${\beta}_k^w {a}_k^w{\beta}_k^w  = \re {\beta}_k^w
 {a}_k^w{\beta}_k^w$ since ${a}_k$ and ${\beta}_k$ are real,
thus we find that  
\begin{equation*} 
 {\beta}_k^w {a}_k^w{\beta}_k^w =
 \re\left([{\beta}_k^w,{a}_k^w]{\beta}_k^w +
   {a}_k^w B_k^w  \right)
= \frac{1}{2}\left[[{\beta}_k^w,{a}_k^w], {\beta}_k^w\right] +
\frac{1}{2}\left({a}_k^wB_k^w +  B_k^w{a}_k^w\right)
\end{equation*}
where $B_k^w = {\beta}_k^w{\beta}_k^w \in \op S(H_1^{-1/2},g^\sharp)$
is symmetric. We find from ~\thetag{18.4.8} in ~\cite{ho:yellow} that
$B_k = {\beta}_k^2 + r_k$ with ${\beta}_k^2 \in
S(H_1^{-1/2},g^\sharp)\bigcap S^+(1,g^\sharp)$ and $r_k \in
S(H_1^{1/2},g^\sharp)$ since $\partial {\beta}_k^2 = 2 {\beta}_k \partial
{\beta}_k $ where $\partial {\beta}_k \in 
S(H_1^{1/4},g^\sharp)$.  By using the expansion ~\thetag{18.4.8}
in ~\cite{ho:yellow} and  that
$B_k \in S^+(1,g^\sharp)$, we find that
$ 
 \frac{1}{2}\left({a}_k^wB_k^w +  B_k^w{a}_k^w\right) \cong
 \left({a}_kB_k \right)^w = (a_k({\beta}_k^2 + r_k))^w$ modulo $\op
 S(MH_1^{},g^\sharp)$ since the operator is symmetric. 
We also obtain that
$[{\beta}_k^w,{a}_k^w] \cong \frac{1}{i}\{\,{\beta}_k,
{a}_k\,\}^w$ modulo $\op S(MH_1^{5/4},g^\sharp)$ where $\{\,{\beta}_k,
{a}_k\,\} \in S(MH_1^{3/4},g^\sharp)$. Thus we find that
$\left[[{\beta}_k^w,{a}_k^w], {\beta}_k^w\right] \cong
-\{\,\{\,{\beta}_k, {a}_k\,\},{\beta}_k\,\}^w \cong 0$ modulo
$\op S(MH_1^{}, g^\sharp)$, which gives ~\eqref{9.8}.
Since the operator in~\eqref{9.8} is symmetric, $r_k$ is real.
This proves the Lemma.
\end{proof}

We shall finish the paper by giving a proof of Proposition~\ref{apsprop}.

\begin{proof}[Proof of Proposition~\ref{apsprop}]
Let $B^{Wick} = ({\delta}_0 + {\varrho}_0)^{Wick}$, where ${\delta}_0
+ {\varrho}_0$ is the ``pseudo-sign'' given by
Proposition~\ref{apsdef}. We find that $B^{Wick} = b^w  = ({\delta}_1 +
{\varrho}_1)^w$ where $b(t,w) \in  L^\infty(\br, S(H_1^{-1/2}, g^\sh) +
S^+(1, g^\sh))$ is given by
Proposition~\ref{wickweyl} for $|t| \le 1$.  Now
${\partial_t } ({\delta}_0 + {\varrho}_0) \ge m_{1}/C_1$ in $\mathcal
D'({\mathbf R})$ when $|t| < 1$ by Proposition~\ref{apsdef}. Let
${\mu}^w  \in  L^\infty(\br, \op S(m_{1},g^\sharp))$ be given by
Proposition~\ref{wprop}, then $ m_1^{Wick} \ge {\mu}^w$. 
Thus we find by
Remark~\ref{distrwick} that 
\begin{equation}\label{dbest}
 {\partial_t }\w{b^w u,u} = \w{{\partial_t } B^{Wick} u,u} \ge C_1^{-1}
 \w{{\mu}^wu,u}\qquad\text{in  $\mathcal D'({\mathbf R})$}
\end{equation}
when $u \in C_0^\infty({\mathbf R}^n)$.  We obtain from
Proposition~\ref{wprop} that there exist positive constants $c_0$ and
$c_1$ so that
\begin{equation}\label{mlow}
\w{{\mu}^wu,u} \ge c_1\mn{u}_{H\left(m_1^{1/2}\right)}^2 \ge
c_0h^{1/2}\mn{u}^2\qquad 
u \in C_0^\infty({\mathbf R}^n). 
\end{equation}
Here $\mn{u}_{H\left(m_1^{1/2}\right)}$ is the norm of the Sobolev space
$H(m_1^{1/2},g^\sharp) =H(m_1^{1/2})$ given by ~\eqref{7.7}
with ${\varrho}=1$ and $k = 1/2$. By
Proposition~\ref{lowersign} we find for almost all 
$|t| \le 1$ that
\begin{equation}
\re \w{(B^{Wick} f^w)\restr{t}u,u} = \re\w{(b^w f^w)\restr{t}u,u} \ge
\w{C_t^wu,u} \quad u \in C_0^\infty({\mathbf R}^n)
\end{equation}
with $C_t \in S(m_{1}(t),g^\sharp)$ uniformly. Thus
we obtain from~\eqref{nest0} and duality that there exists a positive
constant $c_2$ such that
\begin{equation}\label{cest}
|\w{C_t^wu,u}| \le \mn{u}_{H\left(m_1^{1/2}\right)} \mn{C_t^w
    u}_{H\left(m_1^{-1/2}\right)} \le c_2
\mn{u}_{H\left(m_1^{1/2}\right)}^2 \le c_2\w{{\mu}^wu,u}/c_1
\end{equation}
for $u \in C_0^\infty({\mathbf R}^n)$ and almost all $|t| \le 1$. We obtain
Proposition~\ref{apsprop} from~\eqref{dbest}--\eqref{cest}, which
completes the proof of Theorem~\ref{mainthm} and the Nirenberg-Treves
conjecture. 
\end{proof}

\bibliographystyle{amsplain}

\begin{thebibliography}{10}

\bibitem{bf}
Richard Beals and Charles Fefferman, \emph{On local solvability of linear
  partial differential equations}, Ann. of Math. \textbf{97} (1973), 482--498.

\bibitem{bc:sob}
Jean-Michel Bony and Jean-Yves Chemin, \emph{Espace fonctionnels associ\'es au
  calcul de {Weyl}-{H\"{o}r\-mander}}, Bull. Soc. Math. France \textbf{122}
  (1994), 77--118.

\bibitem{de:thesis}
Nils Dencker, \emph{On the propagation of singularities for pseudo-differential
  operators of principal type}, Ark. Mat. \textbf{20} (1982), 23--60.

\bibitem{de:ln}
\bysame, \emph{The solvability of non ${L}^2$ solvable operators}, Journees
  ''\'Equations aux D\'eriv\'ees Partielles'', St.\ Jean de Monts, France,
  1996.

\bibitem{de:suff}
\bysame, \emph{A sufficient condition for solvability}, International
  Mathematics Research Notices \textbf{1999:12} (1999), 627--659.

\bibitem{de:psi}
\bysame, \emph{On the sufficiency of condition $({\Psi})$}, Preprint.

\bibitem{dieu}
Jean Dieudonn\'e, \emph{Foundations of modern analysis}, Academic Press, New
  York and London, 1960.

\bibitem{ho:weyl}
Lars {H\"{o}rmander}, \emph{The {Weyl} calculus of pseudo-differential
  operators}, Comm. Partial Differential Equations \textbf{32} (1979),
  359--443.

\bibitem{ho:yellow}
\bysame, \emph{The analysis of linear partial differential operators}, vol.
  I--IV, Springer Verlag, Berlin, 1983--1985.

\bibitem{ho:conv}
\bysame, \emph{Notions of convexity}, {Birkh\"{a}user}, Boston, 1994.

\bibitem{ho:solv} \bysame, \emph{On the solvability of
pseudodifferential equations}, Structure of solutions of differential
equations (M.~Morimoto and T.~Kawai, eds.), World Scientific, New
Jersey, 1996, pp.~183--213.


\bibitem{ln:2d}
Nicolas Lerner, \emph{Sufficiency of condition $({\Psi})$ for local
  solvability in two 
  dimensions}, Ann. of Math. \textbf{128} (1988), 243--258.

\bibitem{ln:ex}
\bysame, \emph{Nonsolvability in ${L}^2$ for a first order operator satisfying
  condition $({\Psi})$}, Ann.\ of Math. \textbf{139} (1994), 363--393.

\bibitem{ln:coh}
\bysame, \emph{Energy methods via coherent states and advanced
  pseudo-differential calculus}, Multidimensional complex analysis and partial
  differential equations (P.~D. Cordaro, H.~Jacobowitz, and S.~Gidikin, eds.),
  Amer. Math. Soc., Providence, R.I., USA, 1997, pp.~177--201.

\bibitem{ln:per}
\bysame, \emph{Perturbation and energy estimates}, Ann. Sci. École Norm. Sup.
  \textbf{31} (1998), 843--886.

\bibitem{ln:fact}
\bysame, \emph{Factorization and solvability}, Preprint.

\bibitem{ln:private}
\bysame, {Private communication}.

\bibitem{nt}
Louis Nirenberg and Fran\c{c}ois Treves, \emph{On local solvability of linear
  partial differential equations. {P}art {I}: Necessary conditions}, Comm.
  Partial Differential Equations \textbf{23} (1970), 1--38, {\em Part {II}:
  Sufficient conditions}, Comm. Pure Appl. Math. {\bf 23} (1970), 459--509;
  {\em Correction}, Comm. Pure Appl. Math. {\bf 24} (1971), 279--288.

\bibitem{trep}
Jean-Marie Tr\'epreau, \emph{Sur la r\'esolubilit\'e analytique microlocale des
  op\'erateurs pseudodiff\'erentiels de type principal}, Ph.D. thesis,
  Universit\'e de Reims, 1984.

\end{thebibliography}

\providecommand{\bysame}{\leavevmode\hbox to3em{\hrulefill}\thinspace}

\end{document}